\documentclass[11pt]{article}

\usepackage[T1]{fontenc}

\usepackage{xcolor}
\definecolor{revisioncolor}{rgb}{.65,0.05,0}

\usepackage{amssymb}
\usepackage{amsmath}
\usepackage{array}
\usepackage{hhline}
\usepackage{epsfig}
\usepackage{psfrag}
\usepackage{afterpage}
\usepackage{theorem}
\usepackage{stmaryrd}
\usepackage{mathdots}
\newcommand{\pentagram}{\star}
\usepackage{color}
\usepackage{framed}
\makeatletter
\def\udots{\mathinner{\mkern1mu\raise\p@
\vbox{\kern7\p@\hbox{.}}\mkern2mu
\raise4\p@\hbox{.}\mkern2mu\raise7\p@\hbox{.}\mkern1mu}}
\makeatother

\usepackage{eufrak}
\newcommand{\goth}[1]{\EuFrak{#1}}
\newlength{\listeespacee}
\newlength{\tableautasse}
\newcommand{\sss}{{\goth{s}}}

\newcommand{\HH}{\mbox{$\mathbb H$}}
\newcommand{\sign}{\operatorname{sign}}
\newcommand{\SO}{\operatorname{SO}}

\newcommand{\T}{{\rm T}}

\newcommand{\ch}{{\check{h}}}

\newcommand{\cH}{{\check{H}}}

\newcommand{\cOmega}{{\check{\Omega}}}

\newcommand{\dd}{\,\text{\rm d}}

\newcommand{\MM}{{\cal M}}

\newcommand{\MMbar}{{\overline{{\cal M}}}}

\newcommand{\UU}{{\cal U}}

\newcommand{\WW}{{\cal W}}

\newcommand{\II}{{\cal I}}

\newcommand{\KK}{{\cal K}}
\newcommand{\LL}{{\cal L}}

\newcommand{\TT}{{\cal T}}

\newcommand{\FF}{{\cal F}}
\newcommand{\GG}{{\cal G}}

\newcommand{\BB}{{\cal B}}
\newcommand{\HHH}{{\cal H}}
\newcommand{\OO}{{\cal O}}

\newcommand{\cKK}{{\check{\cal K}}}

\newcommand{\cf}{{\check{f}}}
\newcommand{\cu}{{\check{u}}}
\newcommand{\ceta}{{\check{\eta}}}

\newcommand{\Dc}{{\widehat{D}}}
\newcommand{\Abb}{\mathbb A}
\newcommand{\R}{\mathbb R}
\newcommand{\N}{\mathbb N}

\newcommand{\C}{\mathbb C}
\newcommand{\K}{\mathbb K}

\newcommand{\Rad}{\operatorname{Rad}}

\newcommand{\Id}{\operatorname{Id}}

\newcommand{\Zbar}{\overline{Z}}

\newcommand{\vbar}{\overline{v}}

\newcommand{\xbar}{\overline{x}}
\newcommand{\ybar}{\overline{y}}

\newcommand{\partialbar}{\overline{\partial}}
\newcommand{\tr}{\operatorname{tr}}

\newcommand{\Diff}{\operatorname{Diff}}
\newcommand{\End}{\operatorname{End}}

\newcommand{\Mat}{\mbox{\rm Mat}}
\newcommand{\vect}{\operatorname{vect}}

\newcommand{\M}{\mathbb M}

\newcommand{\im}{\operatorname{Im}}
\newcommand{\diag}{\operatorname{diag}}
\newcommand{\Span}{\operatorname{span}}

\newlength{\aux}

\newlength{\arraycolsepsauvegardegenerale}
\newlength{\arraycolsepsauvegardetemporaire}

\newlength{\tabcolsepsauvegardegenerale}
\newlength{\tabcolsepsauvegardetemporaire}

{\setlength{\tabcolsepsauvegardetemporaire}{\tabcolsep}\setlength{\tabcolsep}{0.5\tabcolsepsauvegardegenerale}\begin{tabular}}
{\end{tabular}\setlength{\tabcolsep}{\tabcolsepsauvegardetemporaire}}

\newenvironment{nnarrowtabular}%
{\setlength{\tabcolsepsauvegardetemporaire}{\tabcolsep}\setlength{\tabcolsep}{0.33\tabcolsepsauvegardegenerale}\begin{tabular}}
{\end{tabular}\setlength{\tabcolsep}{\tabcolsepsauvegardetemporaire}}

{\setlength{\arraycolsepsauvegardetemporaire}{\arraycolsep}\setlength{\arraycolsep}{0.5\arraycolsepsauvegardegenerale}\begin{array}}
{\end{array}\setlength{\arraycolsep}{\arraycolsepsauvegardetemporaire}}

\newenvironment{nnarray}%
{\setlength{\arraycolsepsauvegardetemporaire}{\arraycolsep}\setlength{\arraycolsep}{0cm}\begin{array}}
{\end{array}\setlength{\arraycolsep}{\arraycolsepsauvegardetemporaire}}

{\begin{list}{--}{\setlength{\leftmargin}{0cm}\setlength{\itemindent}{1.05cm}\setlength{\labelwidth}{1cm}\setlength{\labelsep}{0.5em}\setlength{\listparindent}{\parindent}}}%
{\end{list}}
{\begin{list}{}{\setlength{\leftmargin}{0cm}\setlength{\itemindent}{0cm}\setlength{\labelwidth}{1.5cm}\setlength{\labelsep}{0.6cm}\setlength{\listparindent}{\parindent}}}%
{\end{list}}

\theoremstyle{change}
\newtheorem{enonce}{}[section]
\newtheorem{de}[enonce]{Definition}
\newtheorem{prop}[enonce]{Proposition}
\newtheorem{prop-de}[enonce]{Proposition/Definition}
\newtheorem{de-prop}[enonce]{Definition/Proposition}
\newtheorem{lem}[enonce]{Lemma}
\newtheorem{te}[enonce]{Theorem}
\newtheorem{cor}[enonce]{Corollary}

\theorembodyfont{\rmfamily}
\newtheorem{rem-notation}[enonce]{Remark/Notation}
\newtheorem{lem-notation}[enonce]{Lemma/Notation}
\newtheorem{def-notation}[enonce]{Definition/Notation}
\newtheorem{notation}[enonce]{Notation}
\newtheorem{rem}[enonce]{Remark}
\newtheorem{rem-de}[enonce]{Remark/Definition}
\newtheorem{importantrem}[enonce]{Important Remark}

\newtheorem{terminology-notation}[enonce]{Terminology/Notation}

\newtheorem{consequence-notation}[enonce]{Consequence/Notation}
\newtheorem{reminder}[enonce]{Reminder}
\newtheorem{reminder-rem}[enonce]{Reminder/Remark}

\newtheorem{example}[enonce]{Example}
\newtheorem{comment}[enonce]{Comment}

\newcounter{claimcounter}
\makeatletter

\renewcommand{\section}{\@startsection{section}{1}{0mm}{\baselineskip}{\medskipamount}{\centering
\bf \mathversion{bold}}}
\makeatother 

\psfrag{1}{{\small\bf (1)}}
\psfrag{1c}{{\small\bf (1$^\C$)}}
\psfrag{2}{{\small\bf (2)}}
\psfrag{2'}{{\small\bf (2')}}
\psfrag{2c}{{\small\bf (2$^\C$)}}
\psfrag{3}{{\small\bf (3)}}
\psfrag{3'}{{\small\bf (3')}}
\psfrag{3c}{{\small\bf (3$^\C$)}}

\author{Charles Boubel}
\begin{document}

\thispagestyle{empty}
\begin{center}
%\maketitle
{\bf %\Large The algebra of the parallel endomorphisms of a germ of pseudo-Riemannian metric
ON THE ALGEBRA OF PARALLEL ENDOMORPHISMS OF A PSEUDO-RIEMANNIAN METRIC\medskip\\}
{\sc Charles Boubel}\smallskip\\July 30th.\@ 2012, revised, March\@ 17th.\@ 2014
\end{center}

\begin{center}\setlength{\aux}{\textwidth}\addtolength{\aux}{-1.7cm}
\parbox{\aux}{\small {\bf Abstract.} On a (pseudo\nobreakdash-)Riemannian manifold $(\MM,g)$, some fields of endomorphisms {\em i.e.\@} sections of $\End(T\MM)$ may be parallel for $g$. They form an associative algebra $\goth e$, which is also the commutant of the holonomy group of $g$. As any associative algebra, $\goth e$ is the sum of its radical and of a semi-simple algebra $\goth s$. This $\goth s$ may be of eight different types, see \cite{boubel2013a}. Then, for any self adjoint nilpotent element $N$ of the commutant of such an $\goth s$ in $\End(T\MM)$, the set of germs of metrics such that $\goth{e}\supset\goth{s}\cup\{N\}$ is non-empty.
We parametrize it. Generically, the holonomy algebra of those metrics is the full commutant $\goth{o}(g)^{\goth{s}\cup\{N\}}$ and then, apart from some ``degenerate'' cases, $\goth{e}=\goth s \oplus (N)$, where $(N)$ is the ideal spanned by $N$. To prove it, we introduce an analogy with complex Differential Calculus, the ring $\R[X]/(X^n)$ replacing the field $\C$. This treats the case where the radical of $\goth e$ is principal and consists of self adjoint elements. We add a glimpse on the case where this radical is not principal.\medskip

\medskip

\noindent{\bf M.S.C.\@ 2010:} 53B30, 53C29, secondary 53B35, 53C10, 53C12.}
\end{center}

We build germs of pseudo-Riemannian metric admitting a parallel nilpotent endomorphism $N\in\Gamma(\End(T\MM))$, together with a semi-simple associative algebra of parallel endomorphisms $\goth s\not\ni N$ (possibly $\goth s=\R\Id$).\\

\noindent{\bf Motivation.} Let $(\MM,g)$ be a (pseudo\nobreakdash-)Riem\-an\-nian manifold, $H$ its holonomy group and $m\in\MM$. The metric $g$ is said to be K\"ahler if it admits an almost complex structure $J$ which is parallel: $DJ=0$ with $D$ the Levi-Civita connection of $g$. A natural question is to ask whether other fields of endomorphisms, {\em i.e.\@} sections of $\End(T\MM)$, may be parallel for a Riemannian metric. The answer is nearly immediate. First, one restricts the study to metrics that do not split into a non trivial Riemannian product, called here ``indecomposable''. Otherwise, any parallel endomorphism field is the direct sum of parallel such fields on each factor (considering as a unique factor the possible flat factor). Then a brief reasoning ensures that only three cases occur: $g$ may be generic {\em i.e.\@} admit only the homotheties as parallel endomorphisms, be K\"ahler, or be hyperk\"ahler {\em i.e.\@} admit two (hence three) anticommuting parallel complex structures. The brevity of this list is due to a simple fact: the action of the holonomy group $H$ of an indecomposable Riemannian metric is irreducible {\em i.e.\@} does not stabilize any proper subspace. In particular, this compels any parallel endomorphism field to be of the form $\lambda\Id+\mu J$ with $J$ some parallel, skew adjoint almost complex structure. Now, such irreducibility fails in general for an indecomposable pseudo-Riemannian metric, so that a miscellany of other parallel endomorphism fields may appear. This gives rise to the following question:
\begin{center}
Which (algebra of) parallel endomorphism fields may a pseudo-Riemannian metric admit ?
\end{center}
Its first natural step, treated here, is local {\em i.e.\@} concerns germs of metrics.

The interest of this question lies also in the following. When studying the holonomy of indecomposable pseudo-Riemannian metrics, the irreducible case may be exhaustively treated: the full list of possible groups, together with the corresponding spaces of germs of metrics (and possibly compact examples) may be provided. After a long story that we do not recall here, this has been done, even for germs of arbitrary torsion free affine connections, see {\em e.g.\@} the surveys \cite{bryant1996,schwachhoefer}. Yet, in the general case, the representation of $H$ may be non-semi-simple ---~see the survey \cite{galaev-leistner2010} of this field, and \cite{galaev-leistner2008} for the Lorentzian case~--- and such an exhaustive answer is out of reach, except perhaps in very low dimension, see {\em e.g.\@} the already long list of possible groups in dimension four in \cite{BB-ikemakhen1997,ghanam-thompson2001}. Thus, intermediate questions are needed: not aiming at the full classification, but still significant. Investigating the commutant $\End(T_m\MM)^H$ of $H$ at some point $m$ of $\MM$, instead of $H$ itself ---~that is to say studying the algebra of parallel endomorphisms~--- is such a question. It has been partially treated, namely for an individual self adjoint endomorphism with a minimal polynomial of degree 2, by G.\@ Kru\v{c}kovi\v{c} and A.\@ Solodovnikov \cite{kruckovic-solodovnikov}. I thank V.\@ S.\@ Matveev for this reference. One may also notice that determining all the parallel tensors, not only the endomorphisms, would mean determining the algebraic closure of the holonomy group $H$. So this work is a step towards this.

Finally, the metrics sharing the same Levi-Civita connection as a given metric $g$ are exactly the $g(\,\cdot\,,U\,\cdot\,)$ with $U$ self adjoint, invertible and parallel. So describing $\End(T_m\MM)^H$ enables to describe those metrics, which is also a useful work. The skew adjoint, invertible and parallel endomorphisms are similarly linked with the parallel symplectic forms.\medskip

Now, as any associative algebra, after Wedderburn--Mal$\check{\rm c}$ev theorem:\medskip\\
\hspace*{\fill} $\End(T_m\MM)^H=\goth{s}\oplus\goth{n}$\hspace*{\fill}\medskip\\
with $\goth n:=\Rad(\End(T_m\MM)^H)$ a nilpotent ideal, its radical, and $\goth s\simeq\End(T_m\MM)^H)/\goth n$ a semi-simple subalgebra.
See \cite{boubel2013a} for details.

Once again, only the semi-simple part $\goth s$ allows an exhaustive treatment. We provided it in \cite{boubel2013a}: $\goth s$ may be of eight different types, including the generic, K\"ahler and hyperk\"ahler types $\goth s\simeq\R$, $\goth s\simeq\C$ and $\goth s\simeq{\mathbb H}$, plus five other ones appearing only for {\em pseudo}-Riemannian metrics; see here Theorem \ref{structure_s} p.\@ \pageref{structure_s}. On the contrary, no list of possible forms for $\goth n$ may presently be given. Recall that, purely algebraically, the classification of nilpotent associative algebras is today out of reach. Even the case of {\em pairs} of {\em commuting} nilpotent matrices is an active subject; we did not find any explicit review of it, but {\em e.g.\@} \cite{basili-iarrobino-katami2010} and its bibliography may be consulted. Besides, our geometric context does not seem to simplify significantly the algebraic nature of $\goth n$. So we treat here a natural first step, the case where $\End(T_m\MM)^H$ contains:\medskip

-- one of the eight semi-simple algebras $\goth s\subset\End(T_m\MM)$ listed in \cite{boubel2013a},\medskip

-- and some given nilpotent endomorphism $N$ not belonging to $\goth s$.\medskip

\noindent It turns out that each such algebra is produced by a non empty set of metrics, which we parametrize. More precisely, we show the following.\medskip

\noindent{\bf Theorem} {\slshape Take $\goth s$ one of the eight algebras listed in  Theorem 1.10 of \cite{boubel2013a}, reproduced here in section \ref{preliminaire} p.\@ \pageref{structure_s}, and $N$ any self adjoint nilpotent element of the commutant of $\goth s$ in $\End(T_m\MM)$. The set of germs of metrics such that $\End(T_m\MM)^H$ contains $\goth s\cup\{N\}$ {\em i.e.\@} such that the elements of $\goth s\cup\{N\}$ extend as parallel endomorphisms fields, admits a parametrization (explicit, or obtained {\em via} Cartan-K\"ahler theory).

On this set, generically, equality $H^0=({\rm O}(g)^{\goth{s}\cup\{N\}})^0$ holds.

Apart from some  exceptional cases,  indicated in Corollary \ref{cor_forme_e} p.\@ \pageref{cor_forme_e}, that behave differently for reasons of Linear Algebra, the algebra of parallel endomorphisms of those metrics is $\goth{s}\oplus(N)$, where $(N)$ is the ideal spanned by $N$ in the algebra $\langle \goth s\cup\{N\}\rangle$.\medskip}

This shows a last motivation: producing metrics $g$ with a parallel field of nilpotent endomorphisms $N$ means producing commutants ${\rm O}(g)^N$ as holonomy groups. Examples of such metrics have been recently built by A.\@ Bolsinov and D.\@ Tsonev \cite{bolsinov-tsonev}. Involving additionaly $\goth s$ in the theorem means that we do the same work with classical holonomy groups as ${\rm U}(p,q)$, ${\rm Sp}(p,q)$ {\em etc.\@} instead of ${\rm O}(g)\simeq{\rm O}(p,q)$.\medskip

\noindent{\em Remark.} Non null parallel endomorphisms yield simple and quite strong consequences on the Ricci curvature. See Section 3 of \cite{boubel2013a}.\medskip

\noindent{\bf Contents and structure of the article.} It consists of five parts. Part 1 recalls the eight possible types for $\goth s$ given in \cite{boubel2013a}, and standard facts about nilpotent endomorphisms. Part 2 introduces an analogy between manifolds with a complex structure and manifolds with a ``nilpotent structure'' {\em i.e.\@} an integrable field of nilpotent endomorphisms: an analogue of complex differential calculus arises, $N$ and $\R[X]/(X^n)$ replacing $J$ and $\C=\R[X]/(X^2+1)$. Counterparts of holomorphic functions, of their power series expansion, appear. See in particular Definition \ref{defnilomorphic} and Theorem \ref{development}. Part 3 uses Part 2 to parametrize the set $\mathcal G$ of germs of metrics $g$ admitting a parallel nilpotent endomorphism field $N$. More exactly, we deal with the case where $N$ is $g$-self adjoint. Indeed if some $N\in\goth n$ is parallel, so are its self- and skew-adjoint parts $\frac12(N\pm N^\ast)$, so it is natural to study first the cases $N^\ast=\pm N$. The case $N^\ast=-N$ demands some more theory ---~introducing counterparts of $\partial$, $\partialbar$, of the Dolbeault lemma {\em etc.} We hope to publish it later. Part 4 uses \cite{boubel2013a} and Part 3 to show the main result, Theorem \ref{realisation_nilomorphe} p.\@ \pageref{realisation_nilomorphe}. If $g$ is hyperk\"ahler or of a similar type, this is done by solving an exterior differenial system exacty as is done by R.\@ Bryant in \cite{bryant1996}, but in the framework of ``$\R[X]/(X^n)$-differential calculus'' introduced in Part 2. Corollary \ref{cor_forme_e} gives the form of $\End(T_m\MM)^H$ in each case. To show it we need to compute the general matrix of the elements of all the algebras involved here: commutants, bicommutants {\em etc.} They are also of practical interest, so we gathered them in Lemma \ref{lemme_commutants}, based itself on Notation \ref{notation_forme_e1}--\ref{bases_privilegiees}. Part 5, a lot shorter, uses \cite{boubel2013a} and Part 2 to give a glimpse, through a simple example, on the case where the holonomy group is the commutant O$(g)^{\{N,N'\}}$ of {\em two} algebraically independent nilpotent endomorphisms.\medskip

\noindent{\bf General setting and some general notation.} $\MM$ is a simply connected manifold of dimension $d$ and $g$ a Riemannian or pseudo-Riem\-an\-nian metric on it, whose holonomy representation does not stabilize any nondegenerate subspace {\em i.e.\@} does not split in an orthogonal sum of subrepresentations. In particular, $g$ is not a Riemannian product. We set $H\subset\SO^0(T_m\MM,g_{|m})$ the holonomy group of $g$ at $m$ and $\goth h$ its Lie algebra. As $\MM$ is simply connected, dealing with $H$ or $\goth h$ is indifferent. Let $\goth{e}$ be the algebra $\End(T_m\MM)^{\goth{h}}$  of the parallel endomorphisms of $g$ ---~to commute with $\goth h$ amounts to extend as a parallel field~--- and let ${\goth e}={\goth s}+{\goth n}$ be its Wedderburn--Mal$\check{\rm c}$ev decomposition. If $A$ is an algebra and $B\subset A$, we denote by $\langle B\rangle$, $(B)$, and $A^B$ the algebra, respectively the ideal, spanned by $B$, and the commutant of $B$ in $A$. When lower case letters: $x_i, y_i$ {\em etc.\@} denote local coordinates, the corresponding upper case letters: $X_i, Y_i$  {\em etc.\@} denote the corresponding coordinate vector fields. Viewing vector fields $X$ as derivations, we denote Lie derivatives ${\mathcal L}_Xu$ also by $X.u$.

The matrix $\diag(I_p,-I_q)\in\mathrm{M}_{p+q}(\R)$ is denoted by \addtolength{\arraycolsep}{-.5ex}$I_{p,q}$, $\mbox{\footnotesize$\left(\begin{array}{cc}0&-I_{p}\\I_{p}&0\end{array}\right)$}\in\mathrm{M}_{2p}(\R)$\addtolength{\arraycolsep}{.5ex} by $J_p$ and $\mbox{\footnotesize$\left(\begin{array}{cc}0&I_{p}\\I_{p}&0\end{array}\right)$}\in\mathrm{M}_{2p}(\R)$ by $L_p$. If $V$ is a vector space of even dimension $d$, we recall that an $L\in\End(V)$ is called {\em paracomplex} if $L^2=\Id$ with $\dim\ker(L-\Id)=\dim\ker(L+\Id)=\frac d2$.

Finally, take $A\in\Gamma(\End(T\MM))$, paracomplex or nilpotent. If it is integrable {\em i.e.\@} if its matrix is constant in well-chosen local coordinates, we call it a ``paracomplex structure'' or a ``nilpotent structure'', like a complex structure, as opposed to an almost complex one.\medskip

\noindent{\bf Acknoledgements.} I thank W.\@ Bertram and V.\@ S.\@ Matveev  for the references they indicated to me, L.\@ B\'erard Bergery and S.\@ Gallot for two indications, M.\@ Audin, P.\@ Mounoud and P.\@ Py for their comments on the writing of certain parts of the manuscript, and the referee for pointing out a mistake in what is now Corollary \ref{cor_forme_e}, and for his careful reading.

\section{Preliminaries: the algebra $\goth s$; standard facts about nilpotent endomorphisms.}\label{preliminaire}

 See Theorem 1.10 of \cite{boubel2013a} for the proof and details on each case.
 
\begin{te}\label{structure_s}
The semi simple part $\sss$ of $\End(T_m\MM)^H$ is of one of the following types, where $\underline J$, $J$,  and $L$ denote respectively self adjoint complex structures and skew adjoint complex and paracomplex structures.\medskip

\noindent{\bf  (1) generic}, $\sss=\vect(\Id)\simeq\R$\medskip

\noindent{\bf  (1$^{\C}$) ``complex Riemannian''}, $\sss=\vect(\Id,\underline J)\simeq\C$. Here $2|d$, $d\geqslant4$ and $\sign(g)=(\frac{d}{2},\frac{d}{2})$.\medskip

\noindent{\bf  (2)(pseudo\nobreakdash-)K\"ahler}, $\sss=\vect(\Id,J)\simeq\C$. Here $2|d$.\medskip

\noindent{\bf  (2') parak\"ahler}, $\sss=\vect(\Id,L)\simeq\R\oplus\R$. Here $2|d$, $\sign(g)=(\frac{d}{2},\frac{d}{2})$.\medskip

\noindent{\bf  (2$^{\C}$) ``complex K\"ahler''}, $\sss=\vect(\Id,\underline J,L,J)\simeq\C\oplus\C$. Here $4|d$ and $\sign(g)=(\frac{d}{2},\frac{d}{2})$.\medskip

\noindent{\bf  (3) (pseudo\nobreakdash-)hyperk\"ahler}, $\sss=\vect(\Id,J_1,J_2,J_3)\simeq\HH$. Here $4|d$.\medskip

\noindent{\bf  (3') ``para-hyperk\"ahler''}, $\sss=\vect(\Id,J,L_1,L_2)\simeq{\rm M}_2(\R)$. Here $4|d$ and $\sign(g)=(\frac{d}{2},\frac{d}{2})$.\medskip

\noindent{\bf  (3$^{\C}$) ``complex hyperk\"ahler''}, $\sss=\vect(\Id,\underline J,J,L_1,L_2,\underline JJ,\underline JL_1,\underline JL_2)\simeq{\rm M}_2(\C)$. Here $8|d$ and $\sign(g)=(\frac{d}{2},\frac{d}{2})$.\medskip

Each type is produced by a non-empty set of germs of metrics. On a dense open subset of them, for the $C^2$ topology, the holonomy group of the metric is the commutant ${\rm SO}^0(g)^{\goth s}$ of ${\goth s}$ in ${\rm SO}^0(g)$.
\end{te}

\begin{rem}\label{groupes_dholonomie}If $G$ is a subgroup of GL$_d(\K)$, let ${\bf V}$ be its standard representation in $\K^d$, ${\bf V}^\ast:g\mapsto(\lambda\mapsto\lambda\circ g^{-1})$ be the contragredient one and, if $\K=\C$, $\overline{\bf V}^\ast$ its complex conjugate. The possible signatures of $g$ and the generic holonomy group and representation are the following.\medskip

{\small\noindent\hspace*{\fill}{\begin{nnarrowtabular}{cccccccc}
{\bf (1)}&{\bf (1$^\C$)}&{\bf (2)}&{\bf (2')}&{\bf (2$^\C$)}&{\bf (3)}&{\bf (3')}&{\bf (3$^\C$)}\\\hline
$(p,q)$&$(p,p)$&$(2p,2q)$&$(p,p)$&$(2p,2p)$&$(4p,4q)$&$(2p,2p)$&$(4p,4p)$\\\hline
SO$^0(p,q)$&SO$(p,\C)$&U$(p,q)$&GL$^0(p,\R)$&GL$(p,\C)$&Sp$(p,q)$&Sp$(2p,\R)$&
Sp$(2p,\C)$
\\\hline
${\bf V}$&${\bf V}$&${\bf V}$&${\bf V}\oplus{\bf V}^\ast$&${\bf V}\oplus\overline{\bf V}^\ast$&${\bf V}$&${\bf V}\oplus{\bf V}^\ast$&${\bf V}\oplus\overline{\bf V}^\ast$
\end{nnarrowtabular}}\hspace*{\fill}}

\end{rem}

Then we sum up standard facts, in a presentation of our own. This makes the article self-contained and enables to introduce some coherent notation used all along. Let $N$ be in $\End(\R^d)$, nilpotent of index $n$.\medskip

We set $F^{a,b}=\im N^{a}\cap\ker N^b$, for $a\in\llbracket1,n-1\rrbracket$ and $b\in\llbracket1,n-a\rrbracket$. Though we do not use the $F^{a,b}$ explicitly here, we had them throughout in mind. They are ordered by inclusion as shown in Table \ref{table_F} p.\@ \pageref{table_F}.
\begin{table}[t!]
{\small\hspace*{\fill}$\displaystyle\begin{array}{cccccccccc}
&F^{\ast,0}=\!\!\!\!&\multicolumn{3}{l}{F^{n,\ast}=\im N^{n}=\{0\}}\\[-.4ex]
&\parallel&\cap\\
&\ker N^0&F^{n-1,1} &\multicolumn{3}{l}{=\im N^{n -1}}                  &         & \\[-.4ex]
&&      \cap        &         &                            &         &                           &         & \\ 
&&    F^{n-2,1}       & \subset & F^{n-2,2}&\multicolumn{3}{l}{=\im N^{n -2}}     & \\[-.4ex]
&&     \cap         &         &       \cap                 &         &                           &         & \\ 
&&    F^{n-3,1}       & \subset &    F^{n-3,2}                 & \subset &F^{n-3,3}&\multicolumn{3}{l}{=\im N^{n -3}}\\[-.4ex]
&&     \cap         &         &        \cap                &         &            &         & \\[-.7ex]
&&     \vdots       &         &        \vdots              &         &                           & \ddots  & \\
&&     \cap       &         &        \cap        \\
&&F^{0,1}&\subset & F^{0, 2}      & \subset &
\cdots         & \subset & F^{0, n}&=\im N^0\\
&&\parallel&&\parallel&&&&\parallel&\parallel\\
&&\ker N&&\ker N^2&&&&\ker N^n&\!\!\!=\R^d\phantom =\\
\end{array}$\hspace*{\fill}}
\caption{\label{table_F}The $F^{a,b}=\im N^{a}\cap\ker N^b$, defined for $(a,b)\in\N^2$. Those for $a\not\in\llbracket n-b,n\rrbracket$, seemingly absent, are equal to some $F^{a',b'}$ present here.}
\end{table}

\begin{notation}\label{notationNpointwise}
{\bf (i)} The invariant factors of $N$ are:
\begin{center}$(\underbrace{X,\ldots,X}_{\text{\tiny$d_1$ times}},\underbrace{X^2,\ldots,X^2}_{\text{\tiny$d_2$ times}},\ldots,\underbrace{X^{n},\ldots,X^{n}}_{\text{\tiny$d_n$ times}})=\bigl((X^{a})_{k=1}^{d_{a}}\bigr)_{a=1}^{n}$,
\end{center}
for some $n$-tuple $(d_a)_{a=1}^n$. We call the $(d_a)_a$ the {\em characteristic dimensions} of $N$, see {\bf (ii)} for a justification. We set $D_a:=\sum_{k=1}^ad_k$ and $D_0:=0$.

{\bf (ii)} We denote by $\pi$ the projection $\R^d\rightarrow\R^d/\im N$. Then for each $a\in\llbracket1,n\rrbracket$, $d_a=\dim(\pi(\ker N^{a})/\pi(\ker N^{a-1}))$. For any $(a,b)$, $F^{a,b}/(F^{a,b-1}+F^{a+1,b})$ is canonically isomorphic, through $N^a$, to \linebreak[4]
$\pi(\ker N^{b-a})/\pi(\ker N^{b-a-1})$.
\end{notation}

\begin{rem-de}\label{remNpointwise} Let $\R[\nu]$ the real algebra generated by $\nu$ satisfying the unique relation $\nu^n=0$, {\em i.e.\@} $\R[\nu]=\R[X]/(X^n)\simeq\R[N]$. Setting $\nu V:=N(V)$ for $V\in\R^d$ turns $\R^d$ into an $\R[\nu]$-module. As such, $\R^d\simeq\prod_{a=1}^{n}(\nu^{n-a}\R[\nu])^{d_{a}}$ {\em i.e.\@} $d_1$ factors on which $\nu$ acts trivially, $d_2$ factors on which $\nu$ is 2-step nilpotent {\em etc.} Notice that this isomorphism is not canonical, even up to an automorphism of each of the factors. We set $D:=D_n=\sum_{a=1}^nd_a$. We define a $D$-tuple of vectors $\beta=(X_i)_{i=1}^D$ to be an {\em adapted spanning family} of $\R^d$ as an $\R[\nu]$-module if each $(X_i)_{i=1+D_a}^{D_{a+1}}$, pushed on the quotient, is a basis of $\pi(\ker N^{a+1})/\pi(\ker N^{a})$, see Notation \ref{notationNpointwise}. In other terms, $\beta$ spans $\R^d$ as an $\R[\nu]$-module and the only relation the $X_i$ satisfy is $\nu^aX_i=0$ for $D_{a-1}<i\leqslant D_a$; or: each $(X_i)_{i=1+D_{a-1}}^{D_a}$ spans the factor $(\nu^{n-a}\R[\nu])^{d_{a}}$ as an $\R[\nu]$-module. It is a basis if and only if the $\R[\nu]$-module $\R^d$ admits bases {\em i.e.\@} it is free {\em i.e.\@} $d_a=0$ for $a<n$.

{\em We denote by $n(i)$ the nilpotence index of $N$ on each submodule $\langle X_i\rangle$}; so $n(i)=a$ for $D_{a-1}<i\leqslant D_a$. We denote by $(X_i,(Y_{i,a})_{a=1}^{n(i)-1})_{i=1}^D$ the basis $(X_i,(N^aX_{i})_{a=1}^{n(i)-1})_{i=1}^D$ of $\R^d$ as an $\R$-vector space.
\end{rem-de}

Now let $g$ be a symmetric bilinear form on $\R^d$ such that $g(N\,\cdot\,,\,\cdot\,)=g(\,\cdot\,,N\,\cdot\,)$.

\begin{rem-de}\label{formes_quotient}For each $a\in\llbracket1,n\rrbracket$, the symmetric bilinear form $g(\,\cdot\,,N^{a-1}\,\cdot\,)$ is well defined on the quotient $\pi(\ker N^{a})/\pi(\ker N^{a-1})\simeq\R^{d_a}$. Indeed, if $X\in\ker N^{a}$ and $Y\in\im N$, so $Y=NZ$, then: $g(X, N^{a-1}Y)=g(X,N^{a}Z)=g(N^{a}X,Z)=0.$

We denote by $(r_a,s_a)$ its signature; $r_a+s_a\leqslant d_a$. It is standard that the couple $(N,g)$ is characterized up to conjugation by this family of dimensions and signatures, called here its {\em characteristic signatures}. See {\em e.g.\@} the elementary exposition \cite{leep-schueller}. The form $g$ is non degenerate if and only if each $g(\,\cdot\,,N^{a-1}\,\cdot\,)$ is and then, $[\,\cdot\,]$ being the floor function:
$$\sign(g)=\bigl(\sum_{a=1}^n\left[{\textstyle \frac a2}\right]d_a+\sum_{\text{$a$ odd}}r_a,\sum_{a=1}^n\left[{\textstyle\frac a2}\right]d_a+\sum_{\text{$a$ odd}}s_a\bigr).$$
\end{rem-de}

Assuming, by Proposition \ref{defbilinearassociated},  the ``$\R[\nu]$-module'' viewpoint, Propositions \ref{matrice_h_en_general} and \ref{defsignature} follow.

\begin{prop-de}\label{defbilinearassociated} Set $E=\bigl(\R^d,N\bigr)\simeq\prod_{a=1}^{n}(\nu^{n-a}\R[\nu])^{d_{a}}.$ On $E$, the $\R[\nu]$-bilinear forms $h$ are in bijection with the real forms $g$ satisfying $g(N\,\cdot\,,\,\cdot\,)=g(\,\cdot\,,N\,\cdot\,)$, through:
$$h=\sum_{a=1}^n\nu^{a-1}g(\,\cdot\,,N^{n-a}\,\cdot\,).$$
We call such an $h$ {\em the $\R[\nu]$-bilinear form associated with $g$}, and $g$ {\em the real form associated with $h$}. Be careful that $g$ is not the real part of $h$, but the coefficient of the highest power $\nu^{n-1}$ of $\nu$ in $h$.
\end{prop-de}

\begin{prop-de}\label{matrice_h_en_general}Let $h$ be an $\R[\nu]$-bilinear form on $\R^d$ as in \ref{defbilinearassociated}. If $\beta=(X_i)_{i=1}^D$ is an adapted spanning family (see \ref{remNpointwise}) of $\R^d$, $\Mat_\beta(h)=\sum_{a=0}^{n-1}\nu^aH_a\in{\rm M}_D(\R[\nu])$ where:\medskip

-- $H_a=\text{\small $\left(\begin{array}{cc}0&0\\0&\cH^a\end{array}\right)$}$, the upper left null square block, of size $D_{n-1-a}$, corresponding to $\Span_{\R[\nu]}\bigl\{X_i;N^{n-1-a}X_i=0\bigr\}$,\medskip

-- the upper left square block $\cH^a_0$ of $\cH^a$ of size $d_{n-a}$, %{\em i.e.\@}
corresponding to $\Span_{\R[\nu]}\bigl\{X_i;N^{n-1-a}X_i\neq N^{n-a}X_i=0\bigr\}$ is of signature $(r_{n-a},s_{n-a})$ introduced in \ref{formes_quotient}, hence of rank $r_{n-a}+s_{n-a}$.\medskip

\noindent We call the $(r_a,s_a)_{a=1}^D$ the {\em signatures of $h$}. So if $S\oplus \im N=\R^d$, $(r_a,s_a)$ is the signature of the (well defined ) form $h_{n-a}$ on the quotient $(S\cap\ker N^a)/(S\cap\ker N^{a-1})$.
\end{prop-de}

\noindent{\bf Proof.} For the first point: $\nu^ah(X,\,\cdot\,)=0$ as soon as $X\in\ker N^a$.\hfill{\rm q.e.d.}

\begin{prop}\label{defsignature} In Proposition \ref{matrice_h_en_general}, choosing an adequate $\beta$, we may take, for all $a$, $\cH^a$ null except $\cH^a_0=I_{d_a,r_a,s_a}:=\diag(I_{r_a},\linebreak[2]-I_{s_a},\linebreak[2]0_{d_a-r_a-s_a})\linebreak[2]\in{\rm M}_{d_a}(\R[\nu])$. So with such a $\beta$:
\begin{align*}
\Mat_\beta(h)&=\diag\left(\nu^{n-a}I_{d_{a},r_{a},s_{a}}\right)_{a=1}^n\\%JDG
&=\left(\begin{array}{ccc}\nu^{n-1}I_{d_{1},r_{1},s_{1}}\\&\ddots\\&&\nu^{0}I_{d_{n},r_{n},s_{n}}\end{array}\right)\in{\rm M}_D(\R[\nu]).
\end{align*}

\noindent Each block $\nu^{n-a}I_{d_a,r_{a},s_{a}}$ corresponds to the factor $(\nu^{n-a}\R[\nu])^{d_{a}}$ {\em i.e.\@} to:
$$\Span_{\R[\nu]}\bigl\{(X_i)_{D_{a-1}<i\leqslant D_a}\bigr\}=\Span_{\R[\nu]}\bigl\{X_i;N^{a-1}X_i\neq N^{a}X_i=0\bigr\}\subset\R^d.$$
The $\R[\nu]$-conjugation class of $h$ is given by the signatures of $h$, and $h$ is non degenerate if and only if $r_a+s_a=d_a$ for each $a$.
\end{prop}

Finally, take $N$ a nilpotent structure on $\MM$.

\begin{notation}\label{notationNfield}
The distributions $\im N^a$ and $\ker N^a$ are integrable, we denote their respective integral foliations by $\II^a$ and $\KK^a$ and set $\II:=\II^1$. {\bf\mathversion{bold}From now on, $\UU$ is an open neighbourhood of $m$ on which those foliations are trivial.} We still denote by $\pi$ the projection $\UU\rightarrow\UU/\II$.

Saying that $N$ is integrable is saying that there is a coordinate system $(x_i,(y_{i,a})_{a=1}^{n(i)-1})_{i=1}^D$ of $\UU$ such that at each point, the basis
$$\left(X_i,\left(Y_{i,a}\right)_{a=1}^{n(i)-1}\right)_{i=1}^D:=\left(\frac{\partial}{\partial x_i},\left(\frac{\partial}{\partial y_{i,a}}\right)_{a=1}^{n(i)-1}\right)_{i=1}^D$$
is of the type given in Remark \ref{remNpointwise}.
\end{notation}

\section{Introducing a special class of functions}\label{preliminaires_nilomorphes}

We introduce here some material which is a bit more general than our strict subject. Just afterwards, back to our germs of pseudo-Riemannian metrics with a parallel field of nilpotent endomorphisms, it will simplify a lot the statements and the proofs and above all make them natural.\medskip

We still denote $\R[X]/(X^n)$ by $\R[\nu]$. We now mimic the definition of a holomorphic function $f$ from a manifold $\MM$ with a complex structure $J$, to $\C=\R[{\rm i}]$. The latter is such that $\dd f\circ J={\rm i}\dd f$. Here $\MM$ is endowed with a ``nilpotent structure'': an integrable field of endomorphisms such that $N^{n-1}\neq N^n=0$. This leads to the following definition.

\begin{de}\label{defnilomorphic} If $f: (\MM,N)\rightarrow\R[\nu]$, is differentiable, we call it here {\em nilomorphic} (for the nilpotent structure $N$) if $\dd f\circ N=\nu\dd f$.
\end{de}

\begin{notation}If $\eta$ is a function or more generally a tensor with values in $\R[\nu]$, we denote by $\eta_a\in\R$ its coefficient of degree $a$ in its expansion in powers of $\nu$, so that: $\eta=\sum_{a=0}^{n-1}\eta_a\nu^a$.
\end{notation}

\begin{example}\label{coordonnees_nilomorphes}The simplest example of such functions are ``nilomorphic coordinates'', built once again similarly as complex coordinates $z_j:=x_j+{\rm i}y_j$ on a complex manifold. Take the $N$-integral coordinates $(x_i,(y_{i,a})_a)_i$ introduced on $\UU$ in Notation \ref{notationNfield} and set:
$$z_i:=x_i+\nu y_{i,1}+\nu^2y_{i,2}+\ldots+\nu^{n(i)-1}y_{i,n(i)-1}\in\R[\nu].$$
Then each $\nu^{n-n(i)}z_i$ is nilomorphic. Indeed, take $X_i$ any coordinate vector transverse to $\im N$ and $a\in\N$. Then $(N^aX_j).(\nu^{n-n(i)}z_i)=0$ if $i\neq j$ and $(N^aX_i).(\nu^{n-n(i)}z_i)=\nu^{n-n(i)}\nu^a$ (it is immediate if $N^aX_i\neq 0$; besides $N^aX_i=0$ if and only if $a\geqslant n(i)$ so that both sides of the equality vanish simultaneously). In particular $(N^aX_i).z_i=\nu^a(X_i.z_i)$.
\end{example}

\begin{def-notation}\label{def_coordonnees_nilomorphes} We now call the $z_i$ of Example \ref{coordonnees_nilomorphes} themselves ``nilomorphic coordinates'' even if only the $\nu^{n-n(i)}z_i$ are nilomorphic functions. The reason appears in Remark \ref{formes_elementaires}. We also introduce a notation that will much alleviate the use of nilomorphic coordinates:$$(\nu y_i):=\sum_{a=1}^{n(i)-1}\nu^a y_{i,a},\ \text{so that } z_i=x_i+(\nu y_{i}),\ \text{and: }(\nu y):=\bigl((\nu y_i)_{i=1}^D\bigr).$$
\end{def-notation}

\begin{rem}\label{generalisationdefnilomorphe}Definition \ref{defnilomorphic} may be stated for functions with value in any $\R[\nu]$-module.
\end{rem}

\begin{rem} A system of holomorphic coordinates provides an isomorphism between a neighbourhood of any point $m$ of $(\MM,J)$ onto a neighbourhood of the origin in $\C^K$. So do the nilomorphic $(\nu^{n-n(i)}z_i)_i$, from a neighbourhood of any point $m$ of $(\MM,N)$ onto a neighbourhood of the origin in some $\R[\nu]$-module ${\mathbb M}$. There is a small difference: ${\mathbb M}$ is not free {\em i.e.\@} ${\mathbb M}\not\simeq(\R[\nu])^K$ in general, but ${\mathbb M}\simeq\prod_a\left(\nu^a\R[\nu]\right)^{K(a)}$. See the previous section. This is linked to the $\nu^{n-n(i)}$ factoring the coordinates $z_i$.
\end{rem}

\begin{reminder}\label{basic}A tensor $\theta$ on a foliated manifold $(\MM,\FF)$ is said to be {\em basic} for $\FF$, or $\FF$-basic, if it is everywhere, locally, the pull back by $p:\MM\rightarrow\MM/\FF$ of some tensor $\overline \theta$ of $\MM/\FF$.
\end{reminder}

We will need also the following auxiliary definition.

\begin{de}\label{defadapted} Let $\M$ be an $\R[\nu]$-module. A function $\check f:(\UU/\II)\rightarrow\M$ is said to be {\em adapted} (to $N$) if for each $a\in\llbracket0,n-1\rrbracket$, $\nu^a\check f$ is $\pi(\KK^{a})$-basic {\em i.e.\@} constant along the leaves of $\pi(\KK^{a})$. If $\M=\R[\nu]$, this means that each coefficient $\check f_a$ is $\pi(\KK^{n-1-a})$-basic.

Similarly, a (multi)linear form $\ceta$ defined on $\UU/\II$ with values in $\M$ is called {\em adapted} if each $\nu^a\ceta$ is $\cKK^{a}$-basic. If $\M=\R[\nu]$, this means that each coefficient $\ceta_a$ is $\cKK^{n-1-a}$-basic.
\end{de}

Here is the main property of nilomorphic functions we will use. The proof is simple, but this is a key statement so we called it a theorem.

\begin{te}\label{development}Let $\M$ be an $\R[\nu]$-module and $f\in C^{n-1}(\UU,\M)$. Then $f$ is nilomorphic for $N$ if and only if, in any nilomorphic coordinates system $(z_i)_{i=1}^D=(x_i+(\nu y_i))_{i=1}^D$, it reads:$$f=\sum_\alpha\frac1{\alpha!}\frac{\partial^{|\alpha|}\cf}{\partial x^\alpha}(\nu y)^\alpha,$$
where $\cf$ is some adapted function (see Def. \ref{defadapted}) of the coordinates $(x_i)_{i=1}^D$ and where, classically:\vspace{-1ex}
\begin{center}
$\alpha$ is a multi-index $\displaystyle(\alpha_i)_{i=1}^D$, $\quad\displaystyle|\alpha|:=\sum_{i=1}^D\alpha_i$,  $ \quad\displaystyle\alpha!:=\prod_{i=1}^D\alpha_i!\ $,\\$\displaystyle\frac{\partial^{|\alpha|}\cf}{\partial x^\alpha}:=\left(\frac{\partial^{|\alpha_1|}}{\partial x_1^{\alpha_1}}\ldots\frac{\partial^{\alpha_D}}{\partial x_D^{\alpha_D}}\right)\cf$,\ \ and  $\  \displaystyle(\nu y)^\alpha:=\prod_{i=1}^D(\nu y_i)^{\alpha_i}$.
\end{center}
\end{te}

\begin{rem} Theorem \ref{development} is very similar to the fact that a function $(\MM,J)\rightarrow\C$ is holomorphic if and only if it is equal to a power series in the neighbourhood $\UU$ of any point. In the complex case, we may consider that the coordinates $x_i$ and $y_i$ parametrize the integral leaves of $\im J$ (complicated manner to mean the whole $\UU$), and that a single point $m$ is a manifold transverse to this leaf. Then $f$ is holomorphic if and only if it reads, in any holomorphic coordinates system:
$$f=\sum_\alpha\frac1{\alpha!}\left(\frac{\partial^{|\alpha|}f}{\partial z^\alpha}\right)_{|m}z^\alpha.$$
In the formula of Theorem \ref{development} appear:\medskip

-- instead of the $z^\alpha$, the $(\nu y)^\alpha$, which are the powers of the coordinates parametrising the integral leaves $\II$ of $\im N$ (those are not the whole $\UU$),\medskip

-- instead of the value and derivatives of $f$ at the single point $m$ (a ``transversal'' to $\UU$), the values and derivatives of $\cf$, which is $f$ along the level $\TT:=\{(\nu y)=0\}=\{\forall i,(\nu y_i)=0\}$, a transversal to the leaves of $\II$.\medskip

So in the complex case, you choose the value and derivatives of a holomorphic function at some point (ensuring  a convergence condition), the rest of the function is given by a power series. In the ``nilomorphic'' case, you choose the value of $f$ along some transversal $\TT$ to $\II$ (ensuring the ``adaptation'' condition \ref{defadapted}), the rest of the function is given by a power series. As $\nu^n=0$, this series in powers of $(\nu y)$ is even a {\em polynomial}, of degree $n-1$. So $f$ is polynomial along the leaves of $\II$ ---~thus this notion makes sense, in any $N$-nilomorphic coordinates system; see Example \ref{jets} for an explanatory point of view. Transversely to those leaves however, $f$ may be only of class $C^{n-1}$.
\end{rem}

\begin{rem}One of the interests of the development formula of Theorem \ref{development} is that it holds even if the invariant factors of $N$ have different degrees {\em i.e.\@} the $\R[\nu]$-module $(T\MM,N)$ is not free. This will enable to build metrics making $N$ parallel in such cases, {\em e.g.\@} as in Example \ref{exemple_trois}
\end{rem}

\noindent{\bf Proof of Theorem \ref{development}.} {\bf (i) The ``if'' part.} Take $f$ of the form given in the proposition, $X_i$ any coordinate vector transverse to $\im N$ and $a\in\N^\ast$. Let us check that $(N^aX_i).f=\nu^a(X_i.f)$.  If $\alpha=(\alpha_i)_{i=1}^D$, $\alpha\pm1_i$ stands for $(\alpha_1,\ldots,\alpha_{i-1},\alpha_i\pm1,\alpha_{i+1},\ldots,\alpha_D)$.
\begin{gather*}\nu^aX_i.\left(\frac1{\alpha!}\frac{\partial^{|\alpha|}\cf}{\partial x^\alpha}(\nu y)^\alpha\right)=\nu^a\frac{\partial}{\partial x_i}\frac{\partial^{|\alpha|}\cf}{\partial x^{\alpha}}(\nu y)^\alpha=\nu^a\frac1{\alpha!}\frac{\partial^{|\alpha|+1}\cf}{\partial x^{\alpha+1_i}}(\nu y)^\alpha.
\end{gather*}
As $\cf$ is adapted (Def.\@ \ref{defadapted}), $\nu^{n(i)}\cf$ is constant along the leaves of $\KK^{n(i)}$. So, as $X_i\in\ker N^{n(i)}$, and setting $\alpha':=\alpha-1_i$:
\begin{align*}(N^aX_i).\left(\frac1{\alpha!}\frac{\partial^{|\alpha|}\cf}{\partial x^\alpha}(\nu y)^\alpha\right)&=\chi_{\{a<n(i)\}}\frac1{\alpha!}\frac{\partial^{|\alpha|}\cf}{\partial x^{\alpha}}\nu^a\alpha_i(\nu y)^{\alpha-1_i}\\
&=\chi_{\{a<n(i)\}}\nu^a\frac1{\alpha'!}\frac{\partial^{|\alpha'|+1}\cf}{\partial x^{\alpha'+1_i}}(\nu y)^{\alpha'}
\end{align*}
As $(N^aX_i).\left(\frac1{\alpha!}\frac{\partial^{|\alpha|}\cf}{\partial x^\alpha}(\nu y)^\alpha\right)=0$ if $\alpha_i=0$, we get the following equality, which concludes:
$$(N^aX_i).f=\chi_{\{a<n(i)\}}\nu^a\sum_{\alpha'}\frac1{\alpha'!}\frac{\partial^{|\alpha'|+1}\cf}{\partial x^{\alpha'+1_i}}(\nu y)^{\alpha'}=\nu^aX_i.f.$$

{\bf (ii) The ``only if'' part.} If $f$ is nilomorphic, then its restriction $\cf$ to the level $\TT=\{\forall i,(\nu y_i)=0\}$, as a function of the $x_i$, is adapted: if $N^aX=0$, $X.(\nu^a\cf)=(N^aX).\cf=0$. Therefore, denoting by $\cf'$ the restriction of $f$ to $\TT$, viewed as a function of the $x_i$, the function $f'$ defined as $f':=\sum_\alpha\frac1{\alpha!}\frac{\partial^{|\alpha|}\cf'}{\partial x^\alpha}(\nu y)^\alpha$ is nilomorphic by {\bf (i)}.\medskip

\noindent{\em Claim.} A nilomorphic function $g:\UU\rightarrow\M$ null on $\TT$ is null.\medskip

Applying the claim to $g=f-f'$ gives $f=f'$, so that $f$ is of the wanted form. Now we have to prove the claim. For any $i$ and $a>0$, as $g$ is nilomorphic, $(N^aX_i).g=\nu^a(X_i.g)$. In the quotient $\M/\nu\M$, this reads $(N^aX_i).[g]_{\M/\nu\M}=0$ so $[g]_{\M/\nu\M}\equiv0$. This gives rise to an induction: $X_i.g\in\nu\M$ so in $\M/\nu^2\M$, $(N^aX_i).[g]_{\M/\nu^2\M}=0$, hence  $[g]_{\M/\nu^2\M}\equiv0$. By induction we get $[g]_{\M/\nu^b\M}\equiv0$ for all $b$ and finally $g=0$.\hfill{\rm q.e.d.}\medskip

The tangent spaces $(T_m\MM,N)$ are $\R[\nu]$-modules, through: $\nu.X:=NX$. Seeing the action of $N$ as that of a scalar leads naturally to introduce the ``nilomorphic'' version of the tensors: same theory, $\R[\nu]$-linearity replacing $\R$-linearity. Let us take local nilomorphic coordinates $z_i=x_i+(\nu y_i)$ and, after Notation \ref{notationNfield}, $X_i:=\frac{\partial}{\partial x_i}$, $Y_{i,a}:=\frac{\partial}{\partial y_{i,a}}$. Then $(X_i)_{i=1}^D$ is an adapted spanning family of each $T_m\MM$, see Definition \ref{remNpointwise}.

\begin{rem}\label{fchecknoncanonique}The adapted function $\cf$ in Theorem \ref{development} depends in general on the choice of the transversal $\TT=\{(\nu y)=0\}$. We do not study this dependence here. The interested reader may look at the link with the expansion of functions in jet bundles given in Example \ref{jets} to understand a meaning of it. Yet notice that $\cf_0$ is canonical {\em i.e.\@} does not depend on the choice of $\TT$. More generally, the value of $\cf_a$ along each leaf of $\pi(\KK^{n-a})$ does not depend on it either, up to an additive constant.
\end{rem}

\begin{de-prop}\label{defnilomorphicvectorfield}A vector field $V$ on $(\MM,N)$ is called {\em nilomorphic} if $\LL_VN=0$. Equivalently: in nilomorphic coordinates $z_i=(x_i+(\nu y_i))_i$, $V=\sum_iv_iX_i$ with nilomorphic functions $v_i$.\end{de-prop}

\noindent{\bf Proof.} Any vector field reads $V:=\sum_iv_iX_i$ with $v_i:\MM\rightarrow\R[\nu]$. Now $\LL_VN=0$ if and only if, for any $a$ and $j$, $[V,N^{a+1}X_j]=N[V,N^aX_i]$. The $(N^aX_i)_{a,i}$ commute, so $[V,N^{a+1}X_j]=-\sum_i(\LL_{N^{a+1}X_j}v_i)X_i$ and $N[V,N^aX_i]=-\sum_i(\LL_{N^aX_j}v_i)NX_i=\sum_i(\LL_{N^aX_j}v_i)\nu.X_i$. Hence $\LL_VN=0$ $\Leftrightarrow$ $\forall$ $i$, $j$, $a$, $\LL_{N^{a+1}X_j}v_i=\nu\LL_{N^aX_j}v_i$, the result.\hfill{\rm q.e.d.}

\begin{de-prop}\label{defnilomorphicform}Let $\eta$ be some (multi)linear form on $(\MM,N)$, with values in $\R[\nu]$. We say here that $\eta$ is {\em nilomorphic} if:\medskip

{\bf (i)} at each point $m$, it is $\R[\nu]$-(multi)linear {\em i.e.}, if $(V_i)_{i=1}^k\in T_m\MM$:
$$
\forall (a_i)_{i=1}^k\in\N^k, \eta(N^{a_1}V_1,\ldots,N^{a_k}V_k)=\nu^{\sum_ia_i}\eta(V_1,\ldots,V_k),$$

{\bf (ii)} $\LL_{NV}\eta=\nu\LL_{V}\eta$ for all nilomorphic vector field $V$.\medskip

If {\bf (i)} is verified, then {\bf (ii)} means that in nilomorphic coordinates $z_i=x_i+(\nu y_i)$, the coefficients $\eta(X_{i_1},\ldots,X_{i_k})$ of $\eta$ are nilomorphic functions (left to the reader).
\end{de-prop}

\begin{rem}\label{remnilomorphicform}Point {\bf (i)} above implies that $\nu^a\eta(V_1,\ldots,V_k)=0$, {\em i.e.} $\nu^{n-a}|\eta(V_1,\ldots,V_k)$, as soon as some $V_i$ is in $\ker N^a$.  So in particular, setting $\eta=\sum_{a=0}^{n-1}\nu^a\eta_a$ with real $\eta_a$, at each point $m$, each $\eta_a$ is the pull back of a (multi)linear application $(T_m\MM/\ker N^{n-1-a})^k\rightarrow\R[\nu]$.
\end{rem}

Applying Theorem \ref{development} to the coefficients of any nilomorphic (multi)\-linear form $\eta$ gives:

\begin{prop}\label{development_form}Let $\eta$ be an $\R[\nu]$-(multi)linear form on $\UU$. Then $\eta$ is nilomorphic for $N$ if and only if, introducing $(z_i)_{i=1}^D$ as in \ref{def_coordonnees_nilomorphes} it reads:
$$\eta=\sum_{(i_1,\ldots,i_k)}\sum_\alpha\frac1{\alpha!}\frac{\partial^{|\alpha|}\ceta_{i_1,\ldots,i_k}}{\partial x^\alpha}(\nu y)^\alpha \dd z_{i_1}\otimes\ldots\otimes\dd z_{i_k}$$
where the $\ceta_{i_1,\ldots,i_k}$ are adapted functions of $(x_i)_{i=1}^D$ with value in \linebreak[4]
$\nu^{n-\min_{l=1}^k n(i_l)}\R[\nu]$. 
\end{prop}

\begin{rem}\label{formes_elementaires}
So, in the coordinates, the ``elementary'' nilomorphic multilinear forms are the $\nu^{n-\min(n(i_1),\ldots,n(i_k))}\dd z_{i_1}\otimes\ldots\otimes\dd z_{i_1}$. If $k>1$, the $z_i$ are needed to write them, not only the $\nu^{n-n(i)}z_i$;  thus we chose in \ref{def_coordonnees_nilomorphes} to define the former as the ``nilomorphic coordinates''. 
Using them:
$$\eta=\sum_{(i_1,\ldots,i_k)}\sum_\alpha\frac1{\alpha!}\frac{\partial^{|\alpha|}\ceta_{i_1,\ldots,i_k}}{\partial x^\alpha}(\nu y)^\alpha\bigl( \nu^{n-\min(n(i_1),\ldots,n(i_k))}\dd z_{i_1}\otimes\ldots\otimes\dd z_{i_k}\bigr)$$
but here each $\ceta_{i_1,\ldots,i_k}$, valued in $\R[\nu]/(\nu^{\min_{l=1}^k n(i_l)})$, must be such that $\nu^{n-\min_{l}n(i_l)}\ceta_{i_1,\ldots,i_k}$ is adapted. So the expression of Prop.\@ \ref{development_form} is simpler.
\end{rem}

The following result, which is now immediate, characterizes the nilomorphic forms in terms of real ones.

\begin{de-prop}\label{defprenilomorphic} Let $\theta\in\Gamma(\otimes^kT^\ast\MM)$ be a real $k$-linear form on $(\MM,N)$. We call it pre-nilomorphic if:\medskip

{\bf (i)} for any $(V_j)_{j=1}^k$, the $\theta\bigl((V_j)_{j=1}^{i-1},NV_i,(V_j)_{j=i+1}^k\bigr)$ are equal to each other, for all $i$,\medskip

{\bf (ii)} $\LL_{NV}\theta=\LL_{V}\theta(N\,\cdot,\,\cdot\,,\ldots,\,\cdot\,)$ for all nilomorphic vector field $V$.\medskip

\noindent Then the following $\R[\nu]$-valued $k$-linear form is nilomorphic:
$$\Theta:=\sum_{a=0}^{n-1}\nu^a\theta\bigl(N^{n-1-a}\,\cdot,\,\cdot\,,\ldots,\,\cdot\,\bigr).$$
We call it the {\em nilomorphic form associated with $\theta$}. In this sense, any nilomorphic $k$-linear form $\Theta=\sum_{a=0}^{n-1}\Theta_a\nu^a$, with real $\Theta_a$, is associated with its coefficient $\Theta_{n-1}$ ---~which is, necessarily, pre-nilomorphic.
\end{de-prop}

The example and comments below are unnecessary for the following. They are given as they are natural, and give another point of view on nilomorphic functions ---~making a link with another work \cite{bertram2011,bertram-souvay2012}.

\begin{example}\label{jets} Natural manifolds with a nilpotent structure are the jet bundles ${\rm J}^n\WW$ over some differentiable manifold $\WW$. The fibre at some point $m\in\WW$ is $\{{f:}]{-\varepsilon},\varepsilon[\rightarrow\WW$ ; $\varepsilon>0$ and $f(0)=m\}/{\sim}$, where $f\sim g$ if in some neighbourhood of $0$, then in all of them, $\|f(t)-g(t)\|=o(t^n)$ when $t\rightarrow0$. So ${\rm J}^0\WW=\WW$ and ${\rm J}^1\WW=T\WW$. With each local chart $\varphi=(x_i)_{i=1}^d:\OO\rightarrow\R^d$ on some open set $\OO$ of $\WW$ is functorially associated a natural chart $\widetilde\varphi=((x_{i,a})_{a=0}^n)_{i=1}^d:{\rm J}^n\OO\rightarrow\R^d$, defined by:
$$\widetilde\varphi([f])=((x_{i,a})_{a=0}^n)_{i=1}^d=:(\xbar_{a})_{a=0}^n\text{ if }f(t)=\xbar_{0}+t\xbar_{1}+\ldots+t^n\xbar_{n}+o(t^n).$$
A change of chart $\theta$ on $\WW$ induces a change of chart $\widetilde\theta$:  let the successive differentials of $\theta$, up to order $n$, act on the $\xbar_{a}$.
The projections $$\WW\leftarrow{\rm J}^1\WW\leftarrow{\rm J}^2\WW\leftarrow\ldots\leftarrow{\rm J}^n\WW\leftarrow\ldots$$ endow each ${\rm J}^n\WW$ with a flag of foliations $\KK^{1}\subset\ldots\subset\KK^n$: in any chart of the type $\widetilde\varphi$, $\KK^{a}$ is given by the levels of $(\xbar_{b})_{b=0}^{n-a}$. Now $\R[X]/(X^{n+1})$ acts naturally on $T{\rm J}^n\WW$, through the endomorphism $N$ defined as follows. If the path $(f+sg)_{s\in\R}$ represents, in some chart and at $s=0$, some tangent vector $v$ to ${\rm J}^n\WW$ at the point $[f]$, $N(v):=(f+s\widehat g)_{s\in\R}$ with $\widehat g(t)=tg(t)$. This definition is consistent, all this is classical. At $v=(\vbar_a)_{a=0}^n\in T_{[f]}{\rm J}^n\WW$, in some chart of the type $\widetilde\varphi$, $N$ reads:
$$N\bigl((\vbar_a)_{a=0}^n\bigr)=\left(0,(\vbar_a)_{a=0}^{n-1}\right).$$
In this chart, $\Mat(N)$ is constant, block-Jordan, so $N$ is integrable; its invariant factors are the $d$-tuple $(X^{n+1},\ldots,X^{n+1})$. So $({\rm J}^n\WW,N)$ is a manifold with a nilpotent structure. Each $\KK^a$ is the integral foliation of $\ker N^a=\im N^{n+1-a}$. Moreover, the submanifold $\TT_0$ of the jets of constant functions is a privileged transversal to $\II=\KK^n$. {\em Conversely, any manifold $\MM$ with a nilpotent structure $N$ with invariant factors of the same degree $n+1$, and endowed with some fixed transversal $\TT$ to $\II$, is locally modelled on the $n^\text{th}$ jet bundle of the (local) quotient $\MM/\II$.} Explicitly, any nilomorphic coordinate system of $\MM$ is exactly given by some transversal $\TT$ to $\II$ and some chart $\varphi$ of $\MM/\II$. If $\TT$ is fixed on $\MM$, the correspondence $\Psi$, given on the left in nilomorphic coordinates induced by some chart $\varphi$ of $\MM/\II$, and on the right in the chart $\widetilde \varphi$:
$$\begin{array}{rccc}
&(\MM,N)&\!\!\!\rightarrow\!\!\!&{\rm J}^n(\MM/\II)\\
\Psi\!:\!\!\!&
\!\!\!\left(x_i+(\nu y_i)\right)_{i=1}^D=\left(x_i+\sum_{a=1}^{n}\nu^ay_{i,a}\right)_{i=1}^D&\!\!\!\mapsto\!\!\!& \left[\left(x_i+\sum_{a=1}^{n} t^ay_{i,a}\right)_{i=1}^D\right]\end{array}$$
is independent of the choice of $\varphi$ ---~because of Theorem \ref{development}. This way, $(\MM,\Diff(\MM,N,\TT))$ identifies with ${\rm J}^n(\MM/\II)$ with the standard action of $\Diff(\MM/\II)$ on it. In this sense, a manifold with a nilpotent structure with invariant factors $(X^{n+1},\ldots,X^{n+1})$ is, locally, a jet bundle of order $n$ where you ``forgot'' what the submanifold of constant jets is.

Introducing functions $\widetilde u:{\rm J}^n\WW\rightarrow\R[\nu]=\R[X]/(X^{n+1})\simeq{\rm J}^n_{|0}\R$ is therefore natural; they are $(N,\nu)$-nilomorphic if and only if they represent, modulo some constant, the $n$-jet of functions $u:\WW\rightarrow\R$. That is to say, $\widetilde u$ is nilomorphic if and only if there is some $u$ and some $f_0$ in $C^\infty(\MM,\R)$ such that, for all $[f]\in{\rm J}^n\WW$:
$$\widetilde u([f])=[u\circ f] + [f_0],\ \text{where $[\,\cdot\,]$ stands for ``$n$-jet of''.}$$
This condition is independent of the choice of a privileged transversal $\TT$, which amounts to a change of $[f_0]$. So through $\Psi$, it works on any $(\MM,N)$ with the invariant factors of $N$ all of the same degree.\medskip

{\em Note: the case $n=1$ is elementary.} The fibre of ${\rm J}^1\WW=T\WW\overset{\pi}{\longrightarrow}\WW$ is a vector space, so for any $(m,v)\in T\WW$, $T_{(m,v)}T\WW$  identifies with $T_m\WW$. So $N(V):=\dd\pi_{|(m,v)}(V)$ may be viewed as an element of $T_{(m,v)}T\WW$. By construction, this $N\in\End(T{\rm J}^1\WW)$ is 2-step nilpotent. We let the reader check it is the same as the $N$ built above.\medskip

{\em Slight adaptation for the case where $N$ has any invariant factors.} Let $\widehat\KK$ be the integral flag $\pi(\KK^1)\subset\ldots\subset\pi(\KK^{n})\subsetneq\pi(\MM)$ of the $\pi(\ker N^a)$ on the local quotient $\pi(\MM)=\MM/\II$. We define a space ${\rm J}^{\widehat\KK}(\MM/\II)$ of ``$\widehat\KK$-jets'' of functions from $\R$ to $\pi(\MM)$ by: $f\sim g$ if, for each $a$, $f$ and $g$ have the same jet in ${\rm J}^{a}\bigl(\pi(\MM)/\pi(\KK^a)\bigr)$. In coordinates adapted to $\widehat \KK$, a $\widehat\KK$-jet $[f]$ is the data of the coordinates of $f$ up to the order $a$, as soon as they are transverse to $\KK^a$. To define $\Psi$ as above we use, on the left side, coordinates adapted to $\widehat \KK$ and factor those spanning each $\KK^{a+1}\smallsetminus\KK^{a}$ by $\nu^{n-a}$ {\em i.e.\@} we use the $\nu^{n-n(i)}z_i$ introduced in Example \ref{coordonnees_nilomorphes}. Then a function $\widetilde u:\MM\rightarrow\R[\nu]$ is nilomorphic if and only if it is, locally, the $\widehat\KK$-jet of some function $u:\pi(\MM)\rightarrow\R$, plus some other fixed $\widehat\KK$-jet $[f_0]$. 
\end{example}

\begin{rem}\label{lien_bertram}To build a Differential Calculus he calls ``simplicial'' on general topological spaces, see {\em e.g.\@} \cite{bertram2011}, W.\@ Bertram studied those jet bundles (I thank him for the references). With A.\@ Souvay \cite{bertram-souvay2012}, he used truncated polynomial rings $\K[X]/(X^n)$ with $\K$ any topological ring. The obtained formulas are equivalent, in the case where the $\R[\nu]$-module is free {\em i.e.\@} all the Jordan blocks of $N$ have the same size, to some of this section, notably to Theorem \ref{development}. It is not immediately explicit in the statements of \cite{bertram-souvay2012}, but observe the ``radial expansion'' in Th.\@ 2.8, and its consequences {\em e.g.\@} the expansion of a $\K[X]/(X^n)$-valued function on the bottom of p.\@ 14 in the proof of Th.\@ 3.6, or that given in the proof of Th.\@ 2.11. Besides, as non free $\K[X]/(X^n)$-modules are direct sums of free ones, Theorem 4.5 of \cite{bertram-souvay2012} generalizes the principle of these expansions.\end{rem}

\begin{rem}\label{connexion} Theorem \ref{development} has a noticeable consequence: the fact, for a function $\MM\rightarrow\R$, to be {\em polynomial along the leaves of $\II$}, in the form that appears in its statement, makes sense, regardless of the chosen $N$-integral local coordinates. So, as the datum of a complex structure on a manifold $\MM$ induces a real analytic structure on it, that of a nilpotent structure induces some ``polynomial structure'' along the leaves of $\II$. If $N^2=0$, this structure is of degree one {\em i.e.\@} is a flat affine structure. Let us see it directly, without Theorem \ref{development}. Take any $U=NV\in T_m\II$. There is a unique way, modulo $\ker N$, to extend $V$ in a basic vector field along this leaf $\II_m$ of $\II$. This induces a canonical way to extend $U$ along $\II_m$ {\em i.e.\@} $\II_m$  is endowed with a flat affine structure, preserved by $\Diff(\MM, N)$. This amounts to a flat affine connection $\nabla$ on $T\II={\rm J}^1\II$.

In the case where all the Jordan blocks of $N$ have the same size $n$, using the point of view developed in Example \ref{jets}, we see that the $\Diff(\MM, N)$-invariant structure on the leaves of $\II$, providing their ``polynomial structure'' is a flat connection $\nabla$ on its bundle ${\rm J}^{n-1}\II$.
\end{rem}

\section{The germs of metrics making parallel some self adjoint nilpotent endomorphism}\label{germsN}

Here appears the link between what precedes and our subject:

\begin{prop}\label{lien} Let $(\MM,g)$ be a pseudo-Riemannian manifold admitting a parallel field $N\in\Gamma(\End(T\MM))$ of self adjoint endomorphisms, nilpotent of index $n$. Then $N$ is integrable and $g$ is pre-nilomorphic for $N$, in the sense of \ref{defprenilomorphic}. After \ref{defbilinearassociated}, $g$ is the real metric associated with the nilomorphic metric $h:=\sum_{a=0}^{n-1}\nu^{a}g(\,\cdot\,,N^{n-1-a}\,\cdot\,)$.

Conversely, suppose that $h=\sum_{a=0}^{n-1}\nu^ah_{a}$, with $h_a$ real, is some non degenerate symmetric $\R[\nu]$-bilinear  nilomorphic form on a manifold $\MM$ with a nilpotent structure $N$. Set $g:=h_{n-1}$. Then $(\MM,g)$ is pseudo-Riemannian, $N$ is self adjoint and parallel on it ---~and, according to Definition \ref{defbilinearassociated}, $h$ is the nilomorphic metric associated with $g$.
\end{prop}

We need the following lemma, to our knowledge (through \cite{bolsinov-tsonev}) first proven by \cite{shirokov}, then independently by \cite{lehmann-lejeune}; \cite{thompson} is a short recent proof.
\begin{lem}\label{integration_N}A field of nilpotent endomorphisms $N$ on $\R^d$ with constant invariant factors is integrable if and only if its Nijenhuis tensor ${\mathcal N}_N$ vanishes and the distributions $\ker N^k$ are involutive for all $k$.
\end{lem}

\noindent{\bf Proof of the proposition.} $D$ is torsion free and $DN=0$, so ${\mathcal N}_N=0$, and the distributions $\ker N^a$ are involutive: by Lemma \ref{integration_N}, $N$ is integrable. Hence $g$ satisfies Definition \ref{defprenilomorphic}. Indeed, {\bf(i)} is the fact that $N$ is self adjoint. For {\bf(ii)} take $V$ any nilomorphic vector field and check that $(\LL_{NV}g)(A,B)=(\LL_{V}g^1)(A,B)$, where $g^1:=g(\,\cdot\,,\,N\,\cdot\,)$, for any $A$, $B$. By Prop.\@ \ref{defnilomorphicvectorfield} we may suppose, without loss of generality, that the field $V$, and some fields extending $A$ and $B$, are coordinate vector fields of some integral coordinate system for $N$, so that $V$, $A$, $B$, $NV$, $NA$ and $NB$ commute. Then we must check: $(NV).\bigl(g(A,B)\bigr)=V.\bigl(g(A,NB)\bigr)$.
\begin{align*}
&(NV).\bigl(g(A,B)\bigr)\\
=&g({D}_{NV}A,B)+g(A,{D}_{NV}B)\\
=&g({D}_{A}(NV),B)+g(A,{D}_{B}(NV))\quad\text{as }[V,NA]=[V,NB]=0,\\
=&g(N{D}_{A}V,B)+g(A,N{D}_{B}V)\quad\text{as }{D} N=0,\\
=&g({D}_{A}V,NB)+g(A,N{D}_{B}V)\quad\text{as }N^\ast=N,\\
=&g({D}_{V}A,NB)+g(A,N{D}_{V}B)\quad\text{as }[V,A]=[V,B]=0,\\
=&g({D}_{V}A,NB)+g(A,{D}_{V}(NB))\quad\text{as }{D} N=0,\\
=&V.\bigl(g(A,NB)\bigr).
\end{align*}
For the converse part, first $g$ is non degenerate: if $g(V,\,\cdot\,)=0$ then for any $a$, $h_a(V,\,\cdot\,)=h_{n-1}(V,N^{n-1-a}\,\cdot\,)=0$ so $V=0$. Then $N^\ast=N$ is immediate. To ensure ${D} N=0$, it is sufficient to prove that, for any $N$-integral coordinate vector fields $X_i$, $X_j$, $X_k$ and any $a$, $b$, $c$ in $\N$, $g({D}_{N^aX_i}N^bX_j,N^cX_k)=g(N^b{D}_{N^aX_i}X_j,N^cX_k)$.
\begin{align*}
&\ 2g({D}_{N^aX_i}N^bX_j,N^cX_k)\\
=&\ (N^aX_i).\bigl(g(N^bX_j,N^cX_k)\bigr)+(N^bX_j).\bigl(g(N^aX_i,N^cX_k)\bigr)
%JDG
\\&\qquad\qquad\qquad\qquad\qquad\qquad\qquad-(N^cX_k).\bigl(g(N^aX_i,N^bX_j)\bigr)\\
=&\ X_i.\bigl(g(X_j,N^{a+b+c}X_k)\bigr)+X_j.\bigl(g(X_i,N^{a+b+c}X_k)\bigr)%JDG
\\&\qquad\qquad\qquad\qquad\qquad\qquad\qquad-X_k.\bigl(g(X_i,N^{a+b+c}X_j)\bigr)\\
&\qquad\text{as $h$ is nilomorphic, so $g$ pre-nilomorphic, see \ref{defprenilomorphic},}\\
=&\ 2g({D}_{X_i}X_j,N^{a+b+c}X_k),
\end{align*}
which gives in particular the wanted equality.\hfill{\rm q.e.d.}\medskip

We are done: combining Propositions \ref{lien} and \ref{development_form} provides exactly Theorem \ref{theoremenilpotentadj}. In the statement, if needed, see Definitions \ref{defnilomorphicform}, \ref{defbilinearassociated}, \ref{def_coordonnees_nilomorphes}, \ref{defadapted} and \ref{defsignature} for ``nilomorphic'', ``associated with'', ``nilomorphic coordinates'', ``adapted function'' and ``characteristic signatures''.

\begin{te}\label{theoremenilpotentadj}
A nilpotent structure $N$ of nilpotence index $n$ on $\MM$ is self-adjoint and parallel for a pseudo-Riemannian metric $g$ if and only if $g$ is the real metric associated with a metric $h$ nilomorphic for $N$.

In nilomorphic local coordinates $(z_i)_{i=1}^D:=\bigl((x_i+(\nu y_i))\bigr)_{i=1}^D$, $h$ is an $\R[\nu]$-valued, $\R[\nu]$-bilinear metric of the form:
\begin{align*}h&=\sum_{i,j=1}^Dh_{i,j}\dd z_i\otimes\dd z_j,\ \begin{array}[t]{l}\text{with nilomorphic functions $h_{i,j}=h_{j,i}$}\\
\text{applying in $\nu^{n-max(n(i),n(j))}\R[\nu]$}\\\text{and $(h_{i,j})_{i,j=1}^D$ non degenerate,}\end{array}\\
&=\sum_{i,j=1}^D\sum_\alpha\frac1{\alpha!}\frac{\partial^{|\alpha|}\ch_{i,j}}{\partial x^\alpha}(\nu y)^\alpha\dd z_i\otimes\dd z_j,
\end{align*}
with $\alpha$ a multi-index $(\alpha_i)_{i=1}^D$ and $\ch_{i,j}=\ch_{j,i}$ adapted functions of the $(x_i)_i$ giving the properties of $(h_{i,j})_{i,j=1}^D$ above. The characteristic signatures $(r_a,s_a)_{a=1}^n$ of $(N,g)$ are those of $h$.
\end{te}

In the theorem, recall that $h=\sum_{a=0}^{n-1}\nu^{a}g(\,\cdot\,,N^{n-1-a}\,\cdot\,)$, as stated in Proposition \ref{lien}; in particular, $g$ is the coefficient of $\nu^{n-1}$. See also the important Remark \ref{remarqueparametrage}, and a matricial formulation in Remark \ref{remarquenilpotentadj_matrice}.

\begin{rem}Set $\ch=\sum_{i,j=1}^D\ch_{i,j}\dd z_i\otimes\dd z_j$. Then $\ch=\sum_{a=0}^{n-1}\ch_a\nu^a$ where each $\ch_a$ is the value of $g(\cdot\,,N^{n-1-a}\,\cdot\,)$ along the transversal $\{(\nu y)=0\}$ to $\II$. This gives more explicitly the link between $g$ and $h$.
\end{rem}

Using Notation \ref{notation_multi_index}, Corollary \ref{corollairenilpotentadj} translates Theorem \ref{theoremenilpotentadj} into purely real terms.

\begin{notation}\label{notation_multi_index} {\bf (i)} If $((\alpha_{i,a})_{a=1}^{n(i)-1})_{i=1}^D$ is a multi-index designed so that: $\displaystyle y^\alpha:=\prod_{i,a}y_{i,a}^{\alpha_{i,a}}$, then $x^\alpha$ and ${\mathcal L}_{X^\alpha}$ denote:\vspace{-2.5ex}
$$x^\alpha:=\prod_{i,a}x_i^{\alpha_{i,a}}=\prod_ix_i^{\sum_a\alpha_{i,a}}\quad\text{and}\quad{\mathcal L}_{X^\alpha}:=\frac{\partial^{|\alpha|}}{\partial x^\alpha}:=\prod_{i,a}\frac{\partial^{\alpha_{i,a}}}{\partial x_i^{\alpha_{i,a}}}.$$

{\bf (ii)} Using point {\bf (i)}, if $(x_i,(y_{i,a})_{a=1}^{n(i)-1})_{i=1}^D$ are $N$-integral coordinates, if $\eta$ is a (multi)\-linear form on $\UU/\II$ (hence, depending only on the coordinates $x_i$) and if $b\in\llbracket0,n-1\rrbracket$, we set:
$$\eta^{(b)}:=\sum_{\text{\scriptsize\begin{tabular}{c}$\alpha$ such that\\$\sum_{i,a}a\alpha_{i,a}=b$\end{tabular}}}\frac1{\alpha!}\left({\mathcal L}_{X^\alpha}\eta\right) y^\alpha.$$
\end{notation}

\begin{cor}\label{corollairenilpotentadj} In the framework of Theorem \ref{theoremenilpotentadj}, $\alpha$ being a multi-index $((\alpha_{i,a})_{a=1}^{n(i)-1})_{i=1}^D$, $g$ is a metric defined by (on the right side, we make use of Notation \ref{notation_multi_index} {\bf (ii)}):
\begin{align*}
g(X_i,N^cX_j)&=\sum_{b=c}^{n-1}\hspace{-.2em}
\sum_{\text{\scriptsize\begin{tabular}{c}$\alpha$ such that\\$\sum_{i,a}a\alpha_{i,a}=n\!-\!1\!-\!b$\end{tabular}}}\hspace{-1.2em}
\frac1{\alpha!}\frac{\partial^{|\alpha|}}{\partial x^\alpha}\Bigl(\widetilde{B}^{b-c}(X_i,X_j)\Bigr)y^\alpha\\
&=\sum_{b=c}^{n-1}(\widetilde{B}^{b-c})^{(n-1-b)},
\end{align*}
where each $B^a$, for $a\in\llbracket0,n-1\rrbracket$, is a $\pi(\KK^{n-1-a})$-basic symmetric bilinear form on $\pi(\UU)$, non degenerate along $\pi(\KK^{n-a})/\pi(\KK^{n-1-a})$, and $\widetilde B^\star$ denotes $\pi^\ast B^\star$. The characteristic signatures $(r_a,s_a)_{a=1}^n$ of the couple $(N,g)$ are the signatures of the   $B^{n-a}$ restricted to $\pi(\ker N^{n-1-a})$.
\end{cor}

\noindent{\bf Proof.} By definition, the form $h$ of Theorem \ref{theoremenilpotentadj} is equal to \linebreak[4]%JDG
$\sum_{a=0}^{n-1}\nu^ag(\,\cdot\,,N^{n-1-a}\,\cdot\,)$. So $g(X_i,N^cX_j)$ is the coefficient of $\nu^{n-1-c}$ in $h$. Defining the $\widetilde B^a$ by $\ch=\sum_{a=0}^{n-1}\nu^a\widetilde B^{a}$, we let the reader expand $h$ as given in Theorem \ref{theoremenilpotentadj} in powers of $\nu$, proving the Corollary.\hfill{\rm q.e.d.}

\begin{rem}\label{remarqueBh}$B^{a}$ is the matrix of $g(\,\cdot\,,N^{n-1-a}\,\cdot\,)=h_{a}$ on the transversal $\{(\nu y)=0\}$ to $\II$.
\end{rem}

\begin{importantrem}\label{remarqueparametrage}By Lemma \ref{integration_N}, any parallel field of endomorphisms is integrable, so Theorem \ref{theoremenilpotentadj} and Corollary \ref{corollairenilpotentadj} give a parametrization of the set of pseudo-Riem\-an\-nian metrics on $\UU$, with a holonomy representation preserving some fixed (arbitrary) self adjoint nilpotent endomorphism $N$. The parameters are the $(B^a)_{a=0}^{n-1}$, up to an action of $C^\infty(\R^{d-\operatorname{rk} N},\R^{\operatorname{rk} N})$. Indeed, the $B^a$ are chosen freely, and characterize $g$ once the level $\{(\nu y)=0\}$, {\em i.e.\@} some section of $\pi:\MM\rightarrow\MM/\II$, is chosen. As $\dim \II=\operatorname{rk}N$, $C^\infty(\R^{d-\operatorname{rk} N},\R^{\operatorname{rk} N})$ acts simply transitively on those sections. We give here only an idea of this action, in Remark \ref{dependanceT}. Some part of the $B^a$ is invariant under it, see Remark \ref{annulation_impossible}.
\end{importantrem}

\begin{rem}\label{annulation_impossible}The adapted bilinear form $\ch$ in Theorem \ref{theoremenilpotentadj} depends on the chosen trans\-ver\-sal $\TT=\{(\nu y)=0\}$; see a similar remark for nilomorphic functions in \ref{fchecknoncanonique}. Yet the restriction of each $\ch_a=g(\,\cdot\,,N^{n-1-a}\,\cdot\,)$ to the leaves of $\KK^{n-a}$, ---~encoded below by the matrix $\cH^0_a$ in Remark \ref{remarquenilpotentadj_matrice} and by the matrix $B^a_0$ in Example \ref{exemple_trois}~---, is canonical {\em i.e.\@} does not depend on the choice of $\TT$. Indeed for each $a$, $g(\,\cdot\,,N^{a}\,\cdot\,)$ does not pass on the quotient $\pi(\UU/\KK^a)$, but its restriction to the leaves of $\pi(\KK^{a+1})$ does, see Remark \ref{formes_quotient}. This invariant shall be noticed.

Therefore, it is natural to choose coordinates $\bigl(x_i; n(i)=a+1\bigr)=(x_i)_{i=1+D_a}^{D_{a+1}}$ satisfying some property with respect to it, if they exist: {\em e.\@~g.\@} if $h_{n-1-a}=g(\,\cdot\,,N^{a}\,\cdot\,)$ is flat along each leaf of $\pi(\KK^{a+1})$, coordinates such that its matrix $B^{n-1-a}_0$ is $I_{r_{a+1},s_{a+1}}$. In particular if $d_{a+1}=1$ we can take $B^{n-1-a}_0\equiv\pm1$, after the corresponding characteristic signature.

As the $\ch_a$, except for their restriction to $\ker N^{n-a}$, depend on $\TT$, may $\TT$ be chosen such that they satisfy some specific property ? The answer is given in Remarks \ref{remarque_orthogonalite} and \ref{dependanceT}.\end{rem}

The complicated expression of Corollary \ref{corollairenilpotentadj} is simpler in some cases. We present two of them. Case ({\bf A}): all invariant factors of $N$ have the same degree {\em i.e.\@} $\im N^p=\ker N^{n-p}$ for $p\leqslant n$ {\em i.e.\@} $N$ is conjugated to:
\begin{center}{\small$\left(\begin{array}{cccc}0&I_{d_n}\\&\ddots&\ddots\\&&0&I_{d_n}\\&&&0\end{array}\right)$}, \begin{tabular}{l}with $n$ null blocks on the diagonal,\\$d_n$ being the $n^{\text{th}}$ characteristic\\ dimension of $N$, the other ones being null.\end{tabular}
\end{center}
\noindent In other terms, the $\R[\nu]$-module $E=(\R^d,N)%=(\R^{nd_n},N)
$ is free: $E\simeq\R[\nu]^{d_n}$. Case ({\bf B}): the nilpotence index of $N$ is small, namely we took $N^3=0$.

\begin{example}{$\!$({\bf A})}\label{exempleA} After Corollary \ref{corollairenilpotentadj}, a metric $g$ makes $N$ self adjoint and parallel if and only if, in $N$-integral coordinates giving $N$ the block form displayed just above, its matrix reads:
\begin{center}$\Mat(g)=${\small$\left(\begin{array}{cccc}0&\cdots&0&G^{0}\\
\vdots&\udots&\udots&G^{1}\\
0&\udots&\udots&\vdots\vspace*{.15cm}\\
\!G^{0}\!&\!\!\!G^{1}\!\!\!&\cdots&G^{n-1}\end{array}\right)$},
\end{center}
\noindent with the following $G^a$. For each $a$, we denote also by ${B}^a$ the matrix of the form ${B}^a$ introduced in Corollary \ref{corollairenilpotentadj}. Using Notation \ref{notation_multi_index} {\bf (ii)}:
$$G^a=\sum_{b=0}^{n-1-a}(B^{a-b})^{(b)}.$$
The ${B}^a$ are symmetric matrices, function of the coordinates $x_i$, with $ B^{0}$ non degenerate. For each $a$, $B^a$ represents $g(\,\cdot\,,N^{n-1-a}\,\cdot\,)$ along $\{(\nu y)=0\}$. The couple $(N,g)$ has only one characteristic signature, namely $(r_n,s_n)=\sign( B^{0})$. So $G^{0}= B^{0}$, $G^{1}= B^{1}+\sum_{i}\bigl(\frac{\partial B^{0}}{\partial x_i}\bigr)y_{i,1}$ {\em etc.\@}; as an example, let us expand $G^{3}$:
\begin{align*}G^{3}&=B^{3}+{B^{2}}^{(1)}+{B^{1}}^{(2)}+{B^{0}}^{(3)}\\
&=B^{3}+\sum_{i}\Bigl(\frac{\partial B^{2}}{\partial x_i}\Bigr)y_{i,1}+\frac1{2!}\sum_{i,j}\Bigl(\frac{\partial^2 B^{1}}{\partial x_i\partial x_j}\Bigr)y_{i,1}y_{j,1}+\sum_{i}\Bigl(\frac{\partial B^{1}}{\partial x_i}\Bigr)y_{i,2}\\
&+\frac1{3!}\sum_{i,j,k}\Bigl(\frac{\partial^3 B^{0}}{\partial x_i\partial x_j\partial x_k}\Bigr)y_{i,1}y_{j,1}y_{k,1}+\sum_{i,j}\Bigl(\frac{\partial^2 B^{0}}{\partial x_i\partial x_j}\Bigr)y_{i,2}y_{j,1}\\[-2ex]&\hspace*{21.5em}+\sum_{i}\Bigl(\frac{\partial B^{0}}{\partial x_i}\Bigr)y_{i,3}.
\end{align*}
If $n=2$, $\Mat(g)$ is an affine function of the $y_{i,1}$:
$$\Mat(g)=\left(\begin{array}{cc}0&B^0\\B^0&B^1+\sum_i\frac{\partial B^{0}}{\partial x_i}y_{i,1}\end{array}\right).$$
To be totally explicit, on $\R^4$ with coordinates $(x,x',y,y')$,  this means:
$$\Mat(g)=\left(\begin{array}{cc}0&B^0(x,x')\\B^0(x,x')&B^1(x,x')+\frac{\partial B^{0}}{\partial x}y+\frac{\partial B^{0}}{\partial x'}y'\end{array}\right)$$
with $B^0$ and $B^1$ symmetric 2-2 matrices, $B^0$ everywhere nondegenerate.
\end{example}

\begin{example}\label{exemple_trois}{$\!$({\bf B})} Recall that, after Notation \ref{notationNpointwise}:
\begin{center}$\displaystyle d_1=\dim(\pi(\ker N))$, $\displaystyle d_2=\dim(\pi(\ker N^2)/\pi(\ker N))$\\[.5ex]and $\displaystyle d_3=\dim(\pi(\ker N^3)/\pi(\ker N^2))=\dim(\pi(T\MM)/\pi(\ker N^2))$,
\end{center}
Order the coordinates as: $\displaystyle ((y_{i,2})_{i=d_1+d_2+1}^{d_1+d_2+d_3}, (y_{i,1})_{i=d_1+1}^{d_1+d_2+d_3}, (x_i)_{i=1}^{d_1+d_2+d_3})$.
$$\text{Then, if $N^3=0\neq N^2$, }\Mat(N)=\left(\begin{array}{cccccc}
0&0&I_{d_3}&0&0&0\\
0&0&0&0&I_{d_2}&0\\
0&0&0&0&0&I_{d_3}\\
0&0&0&0&0&0\\
0&0&0&0&0&0\\
0&0&0&0&0&0
\end{array}\right),$$
with columns and lines of respective sizes $d_3$, $d_2$, $d_3$, $d_1$, $d_2$, $d_3$. Columns 1, 2--3, 4--6 correspond respectively to the $y_{i,2}$, $y_{i,1}$, $x_i$. Then a metric $g$ makes $N$ self adjoint and parallel if and only if its matrix reads:
$$\Mat(g)=\left(\begin{array}{cccccc}
0&0&0&0&\hspace*{.5\arraycolsep}0&\hspace*{-.5\arraycolsep}G^0\\
\begin{nnarray}{c}0\\0\end{nnarray}&\begin{nnarray}{c}0\\0\end{nnarray}&\begin{nnarray}{c}0\\G^0\end{nnarray}&\begin{nnarray}{c}0\\0\end{nnarray}&\multicolumn{2}{c}{\hspace*{-\arraycolsep}\biggl(\begin{array}{c}G^1\end{array}\biggr)\hspace*{-\arraycolsep}}\\
\begin{nnarray}{c}0\\0\\G^0\end{nnarray}&\multicolumn{2}{c}{\hspace*{-\arraycolsep}
\begin{array}{cc}0&0\\\multicolumn{2}{c}{\hspace*{-\arraycolsep}\biggl(\begin{nnarray}{c}\,\,\,G^1\,\,\,\end{nnarray}\biggr)\hspace*{-\arraycolsep}}\end{array}\hspace*{-\arraycolsep}}&\multicolumn{3}{c}{\hspace*{-\arraycolsep}\Biggl(\begin{array}{ccc}&G^2&\end{array}\Biggr)\hspace*{-\arraycolsep}}
\end{array}\right)$$
 with, using again Notation \ref{notation_multi_index}, {\bf (ii)}
, $G^0=B^0$ and:
$$
\addtolength{\arraycolsep}{-.2ex}G^1=B^1+\left(\begin{array}{cc}0&0\\0&{B^0 }^{(1)}\end{array}\right)
\text{, }G^2=B^2+
\left(\begin{array}{ccc}\begin{nnarray}{c}0\\0\\0\end{nnarray}&\multicolumn{2}{c}{\hspace*{-\arraycolsep}
      \begin{array}{cc}0&0\\\multicolumn{2}{c}{\hspace*{-\arraycolsep}\biggl(\begin{nnarray}{c}{B^1}^{(1)}\end{nnarray}\biggr)\hspace*{-\arraycolsep}}
      \end{array}
\hspace*{-\arraycolsep}}
\end{array}\right)+
\left(\begin{array}{ccc}0&0&0\\0&0&0\\0&0&{B^0}^{(2)}\end{array}\right),\addtolength{\arraycolsep}{.2ex}$$
where the $B^a$ are symmetric matrices:
$\displaystyle B^0$ depending on the $\displaystyle (x_i)_{i>d_1+d_2}$, %\\[.5ex]
$\displaystyle B^1$ on the $\displaystyle (x_i)_{i>d_1}$ %\\[.5ex]
and $\displaystyle B^2$ on all the $\displaystyle (x_i)_{i=1}^{d_1+d_2+d_3}$.
We recall that:
$${B^a}^{(1)}=\sum_{i}\bigl(\frac{\partial B^{a}}{\partial x_i}\bigr)y_{i,1},\quad{B^a}^{(2)}=\frac1{2!}\sum_{i,j}\bigl(\frac{\partial^2 B^{a}}{\partial x_i\partial x_j}\bigr)y_{i,1}y_{j,1}+\sum_{i}\bigl(\frac{\partial B^{a}}{\partial x_i}\bigr)y_{i,2}.$$

Ensuring the non degeneracy of $g$ is ensuring the non degeneracy condition stated in Corollary \ref{corollairenilpotentadj}, {\em i.e.\@} here, $B^0$ is non degenerate and:
\begin{center}$B^1=${\small$\,\left(\!\begin{array}{cc}B_0^1&\ast\\\ast&\ast\end{array}\right)$}  and $B^2=${\small$\,\left(\!\begin{array}{ccc}B_0^2&\ast&\ast\\\ast&\ast&\ast\\\ast&\ast&\ast\end{array}\right)$}, \begin{tabular}[m]{l}with $B^1_0$ and $B^2_0$\\non degenerate.\end{tabular}\end{center}
The characteristic signatures of $(N,g)$, {\em cf.\@} Remark \ref{formes_quotient}, are:
\begin{gather*}
\bigl(\sign(B^2_0), \sign(B^1_0), \sign(B^0_0)\bigr)=\bigl((r_1,s_1), (r_2,s_2), (r_3,s_3)\bigr),\\
\text{so: }\ \sign(g)=(d_2+d_3+r_1+r_3,d_2+d_3+s_1+s_3).
\end{gather*}
If $N^2=0$, relabelling the $G^a$ and $B^a$ we find, setting \addtolength{\arraycolsep}{-.3ex}$B^1=${\small$\left(\begin{array}{cc}
B^1_0&B^{1\prime}\\
^t\!B^{1\prime}&B^{1\prime\prime}
\end{array}\right)$\addtolength{\arraycolsep}{.3ex}}:
%JDG
$$\Mat(N)=\text{\small$\displaystyle\left(\begin{array}{ccc}
0&0&I_{d_2}\\
0&0&0\\
0&0&0
\end{array}\right)$}\ \text{and}$$
$$
\Mat(g)=\text{\small$\displaystyle\left(\begin{array}{ccc}
0&\hspace*{.5\arraycolsep}0&\hspace*{-.5\arraycolsep}G^0\\
\begin{nnarray}{c}0\\G^0\end{nnarray}&\multicolumn{2}{c}{\hspace*{-\arraycolsep}\biggl(\begin{array}{c}G^1\end{array}\biggr)\hspace*{-\arraycolsep}}
\end{array}\right)$}=
\text{\small$\displaystyle\left(\begin{array}{ccc}0&0&B^0\\
0&B^1_0&B^{1\prime}\\
B^0&^t\!B^{1\prime}&B^{1\prime\prime}+{B^0}^{(1)}\end{array}\right)$},
$$
with $B^0$ and $B^1_0$ non degenerate, $B^0$ depending on the $(x_i)_{i>d_1}$ and $B^1$ on all the $(x_i)_{i=1}^{d_1+d_2}$. We recall that ${B^0}^{(1)}=\sum_{i}\bigl(\frac{\partial B^{0}}{\partial x_i}\bigr)y_{i,1}$. In case $d_1=0$ we re-find the end of Example A above.
\end{example}

\begin{rem}\label{remarquenilpotentadj_matrice}For any nilpotence index $n$, ordering the coordinates as $((y_{i,n-1})_i,\ldots,\linebreak[1](y_{i,1})_i,\linebreak[1](x_{i})_i)$, one may build similarly $\Mat(g)$. The principle is the same as when $N^3=0$ and no new phenomenon appears, but the matrix becomes rapidly very cumbersome. So it seems that the use of real coordinates, forms and matrices, if it may be avoided, should be, and replaced by the use of $\R[\nu]$-linear ones like in Theorem \ref{theoremenilpotentadj}, just as complex expressions replace real ones in complex geometry.

The matrix of $h$ given in Theorem \ref{theoremenilpotentadj} reads as follows. Take $(z_i)_{i=1}^D$ nilomorphic coordinates of $\UU$, ordered by increasing values of $n(i)$. So, $(z_i)_{i=1}^{D_a}=(z_i; n(i)\leqslant a)$ parametrize the leaves of $\KK^{a}$. Using Prop.\@ \ref{matrice_h_en_general},
$$\Mat(h)=\sum_\alpha\frac1{\alpha!}\Bigl(\frac{\partial^{|\alpha|}}{\partial x^\alpha}\cH\Bigr)(\nu y)^\alpha$$
with $\cH$ a sum of symmetric matrices $\nu^a\text{\small $\left(\begin{array}{cc}0&0\\0&\cH^a\end{array}\right)$}\in\nu^a{\operatorname M}_D(\R)$ satisfying, for each $a\in\llbracket0,n-1\rrbracket$:\medskip

{\bf (i)} {\em (Adaptation condition, see \ref{defadapted})} the block $\cH^a$, corresponding to the coordinates $\bigl(z_i;\linebreak[1] n(i)\geqslant n-a\bigr)$, depends only on the $\bigl(x_i; n(i)\geqslant n-a\bigr)$,\medskip

{\bf (ii)} {\em (Non degeneracy condition, see  \ref{matrice_h_en_general})} cutting the $(z_i; n(i)\geqslant n-a)$ into $\bigl((z_i; n(i)=n-a),(z_i; n(i)>n-a)\bigr)$, $\cH^a$ splits into {\small $\left(\!\begin{array}{cc}\cH^a_0&\ast\\\ast&\ast\end{array}\right)$} with $\cH^a_0$ non degenerate.\medskip

Finally, notice that $\cH^a=\nu^a B^{a}$, with $B^a=\Mat(g(\,\cdot\,,N^{n-1-a}\,\cdot\,))$ along $\{(\nu y)=0\}$, as appearing in Corollary \ref{corollairenilpotentadj}.
\end{rem}

\begin{rem}\label{remarque_orthogonalite}Another choice of ``preferred'' coordinates $(x_i)_{i=D_a+1}^{D_{a+1}}=\bigl(x_i;\linebreak[1]n(i)=a+1\bigr)$ for some $a$ could be to try to get $B^a=\text{\small$\left(\begin{array}{cc}B^a_0&0\\0&B^{a\prime\prime}\end{array}\right)$}$, that is to say to choose them so that the $\bigl\{\frac{\partial}{\partial x_i}; n(i)>a+1\bigr\}$ be $B^a$-orthogonal to $\KK^{a+1}\cap\{(\nu y)=0\}$. It is impossible, as the orthogonal distribution to $\KK^{n-a}\cap\{(\nu y)=0\}$ is not integrable in general, even by seeking an ``adequate'' transversal $\TT=\{(\nu y)=0\}$ to $\II$. See Remark \ref{dependanceT} for a proof. Notice here that on the contrary, this is possible for {\em alternate} nilomorphic bilinear forms, for which we can even achieve $B^a=\text{\small$\left(\begin{array}{cc}B^a_0&0\\0&0\end{array}\right)$}$, see Lemma \ref{darboux-poincare} and Rem.\@ \ref{onefunctionofDvariables}.

It works however in one case. If $\sum_{c=1}^{b}d_{n-a+c}=1$ for some $b\geqslant c+1$ (all $d_{n-a+c}$ null except one), the $B^{a}$-orthogonal to $\KK^{n-a}\cap\TT$ is 1-dimensional inside of the leaves of $\KK^{n-a+b}\cap\TT$, so may be integrated, within these leaves. This puts a null line under $B^a_0$ and a null column on its right. Moreover, \cite{boubel2007} provides a transversal $\TT$ such that the unique coordinate vector $\frac{\partial}{\partial x_i}$ transverse to $\KK^{n-a}$ in $\KK^{a+b}\cap\TT$ has a constant $B^a$-square, {\em e.g.\@} is isotropic. This puts a constant coefficient $B^{a\prime\prime}$ below $B^a_0$ on the right of it. See Example \ref{example_lorentz_sym}. If moreover $d_{n-a}=1$, you must choose between achieving all this, and achieving $B^{a}_0\equiv\pm1$ as in \ref{annulation_impossible}.
\end{rem}

\begin{example}\label{example_lorentz_sym}With Remarks \ref{remarquenilpotentadj_matrice} and \ref{remarque_orthogonalite}, we treat the Lorentzian case. The indecomposable germs of Lorentzian metric making a non trivial self adjoint endomorphism $N$ parallel are those with holonomy algebra included in:
$$\left\{\text{\small$\left(\begin{array}{ccc}0&L&0\\0&A&-^t\!L\\0&0&0\end{array}\right)$};A\in\goth{so}(n-2),\ L\in\R^{n-2}\right\},
$$
in a basis where $g=\text{\small$\left(\begin{array}{ccc}1&0&0\\0&I_{n-2}&0\\0&0&1\end{array}\right)$}$. Then $N=\text{%JDG\scriptsize
\small$\left(\begin{array}{ccc}0&0&1\\0&0&0\\0&0&0\end{array}\right)$}$ is parallel. Writing $g$ on the form given at the end of \ref{exemple_trois} and setting $B^0=1$ by Remark \ref{annulation_impossible} and $B^{1\prime}=0$ by Remark \ref{remarque_orthogonalite} 
we get:
$$\Mat(g)=\text{\small$\left(\begin{array}{ccc}0&0&1\\0&B^1_0&0\\1&0&b^{1\prime\prime}\end{array}\right)$},
$$
with $B^1_0$ and $b^{1\prime\prime}$ independent of the first coordinate. Again by \ref{remarque_orthogonalite}, one can get $b^{1\prime\prime}$ constant (zero or not), getting a classical form for $g$.
\end{example}

\begin{rem}\label{dependanceT}The adapted bilinear form $\ch$ appearing in Theorem \ref{theoremenilpotentadj} depends on the choice of $\TT=\{(\nu y)=0\}$. Let us illustrate here how this dependence works through the 
case $N^2=0$. In $N$-adapted coordinates:
\begin{center}
$\displaystyle\bigl((y_i)_{i=1}^{d_2},(x_i)_{i=1}^{d_1},(x_i)_{i=d_1+1}^{d_2}\bigr)$, we get: $\Mat(N)=\text{\small$\left(\begin{array}{ccc}0&0&I\\0&0&0\\0&0&0\end{array}\right)$}$
\end{center}
and $\Mat(g)$ is as given at the end of Example \ref{exemple_trois}.
Then $B^0$ is the matrix of $h_0=g(\,\cdot\,,N\,\cdot\,)$, which is well defined on $\MM/\KK$, so does not depend on $\TT$; $B^1$ is the matrix of $\ch_1=g_{|T\TT}$. It depends on $\TT$, which is the image of a section $\sigma$ of $\pi:\MM\rightarrow\MM/\II$. Changing $\TT$ amounts to add to $\sigma$ a vector field $U$ defined along $\TT$, and tangent to $\II$. Indeed, $\II$ is endowed with a flat affine connection $\nabla$, see Remark \ref{connexion}, so identifies with its tangent space. In turn, this field $U$ is equal to $NV$ with $V$ a section of $\pi(T\MM)/\pi(\ker N)$ defined on $\pi(\MM)=\MM/\II$. If $\sum_{i=d_1+1}^{d_2} v_i \frac{\partial}{\partial x_i}$ represents $V$, so that $U=\sum_{i=1}^{d_2} v_i \frac{\partial}{\partial y_i}$, it follows from the expression of $\Mat(g)$ that $B^1$ becomes $B^1_V$ given by:
\begin{align*}
\textstyle B^1_V(\frac{\partial}{\partial x_j},\frac{\partial}{\partial x_k})&=B^1(\frac{\partial}{\partial x_j},\frac{\partial}{\partial x_k})+V.\bigl(B^0(\frac{\partial}{\partial x_j},\frac{\partial}{\partial x_k})\bigr)%JDG
\\&\qquad\qquad\qquad+\sum_{i=1}^D\frac{\partial v_i}{\partial x_j}B^0(V,\frac{\partial}{\partial x_k})+\frac{\partial v_i}{\partial x_k}B^0(\frac{\partial}{\partial x_j},V)%[.5ex]
\end{align*}
{\em i.e.\@} $B^1_V=B^1+{\mathcal L}_VB^0$. The Lie derivative ${\mathcal L}_VB^0$ is well defined, even if $V$ is defined modulo $\ker N$, as $\ker N=\ker B^0$. So $\ch_1$ is defined, through the choice of $\TT$, up to addition of an infinitesimal deformation of $h_0=g(\,\cdot\,,N\,\cdot\,)$ by a diffeomorphism.

So, the orbit $\bigl\{B^1_V, V\in\Gamma(T\pi(\MM))\bigr\}$ of the matrix $B^1$ when $\TT$ varies does not contain in general a matrix such that $B^{1\prime}=0$ and/or $B^{1\prime\prime}=0$, with the notation of Example \ref{exemple_trois}. Getting $B^{1\prime}=0$ is impossible in general if $\dim\bigl(\KK\cap\TT\bigr)$ is greater than 1. Indeed for all $i>d_1$ the 1-form $({\mathcal L}_VB^0)(X_i,\,\cdot\,)$ is closed along the leaves of $\KK\cap\TT$. So $\dd\bigl(\iota_{X_i}B^1\bigr)$ is canonical, independent of $\TT$. This proves the impossibility announced in Remark \ref{remarque_orthogonalite}. Getting $B^{1\prime\prime}=0$ is also impossible in general: see below. 

To simplify, we suppose now that $d_1=0$ {\em i.e.\@} $\ker N=\im N$; now $h_0$ is a non degenerate metric on $\pi(\MM)=\MM/\II$. After the infinitesimal part of Ebin's slice theorem (see {\em e.g.\@} \cite{besse2} 4.2--5 for a short explanation):
$${\rm Sym}^2T^\ast\MM=\bigl\{{\mathcal L}_Vh_0;V\in\Gamma(T\pi(\MM))\bigr\}\oplus\ker\delta_{h_0}$$
with $\delta_{h_0}$ the divergence operator with respect to $h_0$. I thank S.\@ Gallot who let me think to this. So there is a local privileged choice of $\TT$, making $\ch_1$ divergence free; this divergence free $\ch_1$ is canonical; besides it is given by $\frac{D(D-1)}2$ functions of $D$ variables. The choice of $\TT$, depending on $D$ functions of $D$ variables, acts ``with a slice'' on the value of $\ch_1$, and each orbit of this action consists of forms depending also on $D$ functions of $D$ variables. Finally, up to diffeomorphism, $g$ depends on ${D(D-1)}$ functions of $D$ variables: $2\frac{D(D+1)}{2}$ for the couple $(B^0,B^1)$, minus the choice of $\TT$ and of a chart of $\MM/\II$ {\em i.e.\@} twice $D$ functions of $D$ variables.
\end{rem}

\begin{importantrem}{\mathversion{bold}\bf [The resulting algebra $\End(T\MM)^{\goth{h}}$ --- Indecomposability of the obtained metrics]} The form of the resulting algebra $\End(T\MM)^{\goth{h}}$ for metrics built in Theorem \ref{theoremenilpotentadj} is given by Corollary \ref{cor_forme_e} in the next section. This corollary also discusses, together with Remark \ref{indecomposabilite}, if the metrics of Theorem \ref{theoremenilpotentadj} are decomposable.
\end{importantrem}

\begin{rem}In the next section, Corollary \ref{cor_forme_e} shows that in cases {\bf (1)}--{\bf (1$^\C$)}, the holonomy group of a metric making $N$ parallel acts trivially on $\im N^{n_2}$, which is non zero as soon as $(n=)n_1>n_2$. So one may take $N$-adapted coordinates such that the $(Y_{1,a})_{a=n_2}^{n-1}:=\bigl(N^{a}\frac{\partial}{\partial x_1}\bigr)_{a=n_2}^{n-1}$ span $\im N^{n_2}$ and are parallel. Here is a direct way to see it:\medskip

-- Take $N$-adapted coordinates such that $N^{n-1}\frac{\partial}{\partial x_1}\neq 0$; the metric $g_{n-1}:=g(\,\cdot\,,N^{n-1}\,\cdot\,)$ is well-defined on $\MM/\KK^{n-1}$, which is 1-dimensional, so you may modify the coordinate $x_1$ so that $g(\frac{\partial}{\partial x_1},N^{n-1}\frac{\partial}{\partial x_1})\equiv\pm1$ {\em i.e.\@} $B^0=\pm1$. If $n-n_2=1$ we are done.\medskip

-- If $n-n_2>1$, then $\pi(\MM/\KK^{n-2})=\pi(\MM/\KK^{n-1})$ is also 1-dimensional so the matrix $B^1$ is also a scalar. Modifying the level $\TT=\{(\nu y)=0\}$ as in Rem.\@ \ref{dependanceT} to get: $B^1_V=B^1+{\mathcal L}_VB^0$ with $V$ a vector field on $\MM/\KK^{n-1}$, enables to modify arbitrarily $B^1$. An adequate $V$ gives $B^1_V\equiv 0$. The new level $\{(\nu y)=0\}$ is obtained by addition of the field $N^{n-1}V$ {\em i.e.\@}, setting $V=v.\pi\bigl(\frac{\partial}{\partial x_1}\bigr)$, by replacing the coordinate $y_{1,n-1}$ by $y_{1,n-1}-v$.\medskip

-- We modify then similarly the coordinates $y_{1,n-a}$, for $a\in\llbracket1,n-n_2-1\rrbracket$, by induction on $a$: at each step, we choose $V=v.\pi\bigl(\frac{\partial}{\partial x_1}\bigr)$ such that $B^a+{\mathcal L}_VB^0\equiv0$ and replace $y_{1,n-a}$ by $y_{1,n-a}-v$. We let the reader check that this makes the scalars $B^1$, \ldots, $B^{n-n_2-1}$ null.\medskip

\noindent Because of the form of the metric given by Corollary \ref{corollairenilpotentadj}, no coefficient of the metric $g$ depends on the coordinates $(y_{a})_{a=n_2}^{n-1}$. At $m$, we may take the $\bigl(N^a\frac{\partial}{\partial x_1}\bigr)_{a=n-1}^{0}$ such that, on the subspace they span, $\Mat(N)=N_n$ and $\Mat(g)=\pm K_n$, see Notation \ref{notation_forme_e2} {\bf (a)}. The $\pm$ sign is given by the signature of the metric $g_{n-1}$ introduced just above in this remark.
\end{rem}

\begin{example}We end with an example of a pseudo-Riemannian metric where a parallel $N$ naturally arises. I thank L.\@ B\'erard Bergery for having told me that this example works with non locally symmetric metrics. Example \ref{jets} showed that jet bundles carry naturally a nilpotent structure; in particular, tangent bundles carry a 2-step nilpotent structure $N$, see the {\em Note} at its end. Now if $g$ is a (pseudo\nobreakdash-)Riemannian metric on $\MM$ of dimension $d$, $T\MM$ carries a natural metric $\widehat g$, of signature $(d,d)$,  making $N$ self adjoint and parallel. In $TT\MM$, we call ``vertical'' the tangent space $\bf V$ to the fibre, and $\bf H$ its ``horizontal'' complement given by the Levi-Civita connection $D$ of $g$. If $X\in T_{(m,V)}T\MM$, as  the fibre of $T\MM$ is a vector space, $X$ is naturally identified with a vector of $T_m\MM$, that we denote by $X'$. We set $\pi:T\MM\rightarrow\MM$. We define $\widehat g$ by:
\begin{center}{\bf V} and {\bf H} are $\widehat g$-totally isotropic and, if $(X,Y)\in{\bf V}\times{\bf H}$, $\widehat g(X,Y)=g(X',\dd\pi(Y))$.
\end{center}
Notice that $N$ is $\widehat g$-self adjoint. Moreover $\widehat DN=0$, with $\Dc$ the Levi-Civita connection of $\widehat g$.\medskip

{\em Proof.} We must show that, if $X, Y, Z$ are vector fields tangent to $T\MM$, $\widehat g(\Dc_XNY,Z)=\widehat g(\Dc_XY,NZ)$. Let us show it at an arbitrary point $m\in T\MM$. Around $m$, we build a frame field $\beta$ on $T\MM$ as follows. We take normal coordinate vector fields $(U_i)_{i=1}^d$ at $\pi(m)$ on $\MM$. We lift them horizontally on $T\MM$, getting a frame field $\beta_{\bf H}$ of $\bf  H$. Besides we denote by $\beta_{\bf V}$ the frame field of $\bf V$ such that $\beta'_{\bf V}=(U_i)_{i=1}^d$. In other words, $\beta_{\bf V}=N.\beta_{\bf H}$. We set $\beta:=(\beta_{\bf V},\beta_{\bf H})$. In the following, $X$, $Y$, $Z$ denote vectors of $\beta$. Notice that, by construction, $\beta_{\bf V}$ is constant along each fibre of $T\MM$, so $[X,Y]=0$ if $X,Y\in\beta_{\bf V}$; besides $[X,Y]\in{\bf V}$ if $X,Y\in\beta_{\bf H}$, as they are lifts of commuting fields on $\MM$. Finally, by definition of $\bf H$, if $X\in{\bf H}$ and $Y\in{\bf V}$, $[X,Y]\in{\bf V}$ and $[X,Y]'=D_{\dd\pi(X)}Y'$. As the $(U_i)_{i=1}^d$ are normal coordinate vectors at $\pi(m)$, for all $(i,j,k)$, $U_i.(g(U_j,U_k))$ is null at $\pi(m)$; combined with the definition of $\widehat g$, it gives that for any $X,Y,Z$ in $\beta$, $X.\widehat g(Y,Z)$ is also null at $m$. Then, if $X,Y,Z\in\beta$:\begin{align*}
2\widehat g(\Dc_XNY,Z)&=X.\widehat g(NY,Z)+NY.\widehat g(NY,Z)-Z.\widehat g(X,NY)\\
&\qquad\quad-\widehat g(X,[NY,Z])-\widehat g(NY,[X,Z])+\widehat g(Z,[X,NY])\\
&\underset{\text{\tiny(at $m$)}}{=}\!\!\!-\widehat g(X,[NY,Z])-\widehat g(NY,[X,Z])+\widehat g(Z,[X,NY]).
\end{align*}
If $Y\in\beta_{\bf V}$, $NY=0$ so this vanishes. If $Z\in\beta_{\bf V}$ (or if $X\in\beta_{\bf V}$, left to the reader) then:\medskip

-- as $NY\in\beta_{\bf V}$, $[NY,Z]=0$ so $\widehat g(X,[NY,Z])=0$,\medskip

-- $[X,Z]\in{\bf V}$ so $\widehat g(NY,[X,Z])\in\widehat g({\bf V},{\bf V})=\{0\}$,\medskip

-- as $NY\in\beta_{\bf V}$,$[X,NY]\in{\bf V}$ so $\widehat g(Z,[X,NY])\in\widehat g({\bf V},{\bf V})=\{0\}$.\medskip

\noindent So ($X\in\beta_{\bf V}$ or $Y\in\beta_{\bf V}$ or $Z\in\beta_{\bf V}$) $\Rightarrow$ $g_{|m}(\Dc_XNY,Z)=0$. If $X,Y,Z\in\beta_{\bf H}$, at $m$%JDG
, $\widehat g(\Dc_XNY,Z)$ is equal to:
$${\textstyle\frac12}\left(g(\dd\pi(X),D_{\dd\pi(Z)}\dd\pi(Y))-0+g(\dd\pi(Z),D_{\dd\pi(X)}\dd\pi(Y))\right)=0$$ as the $(U_i)_i$ are normal at $m$. Now similarly, at $m$, $2\widehat g(\Dc_XY,NZ)=-\widehat g(X,[Y,NZ])-\widehat g(Y,[X,NZ])+\widehat g(NZ,[X,Y])$, which also vanishes. So at any $m$, $\Dc_XNY=N\Dc_XY$ {\em i.e.\@} $\Dc N=0$, q.e.d.

We could build a similar $\widehat g$ on any jet bundle ${\rm J}^n\MM$, making the $N$ of Ex.\@ \ref{jets} parallel.
\end{example}

\section{\mathversion{bold}Metrics such that $\End(T\MM)^\goth{h}$ has an arbitrary semi-simple part and a non trivial radical\mathversion{normal}}

In \cite{boubel2013a} Section 2 we built metrics such that the semi-simple part $\goth s$ of $\End(T\MM)^\goth{h}$ is in each of the eight possibilities listed by Theorem \ref{structure_s}. Here we built in section \ref{germsN} metrics admitting an arbitrary self adjoint nilpotent structure as a parallel endomorphism; in particular, this makes in general the radical $\goth n$ of $\End(T\MM)^\goth{h}$ non trivial.
Mixing here both arguments, we build metrics with $\goth s$ arbitrary in the list of Theorem \ref{structure_s}, and whose holonomy commutes moreover with some arbitrary self adjoint nilpotent endomorphism $N$. We considered only the case where $N$ is in the commutant of $\goth s$. It is natural: it means that $N$ is a complex endomorphism if $\goth s=\langle J\rangle$ induces a complex structure {\em etc.\@}, see Remark \ref{matrices_commutant}. Besides by Proposition 1.8 of \cite{boubel2013a} and as we suppose $N$ self adjoint, the non commuting case is strongly constrained.

This builds metrics whose holonomy group is the commutant of $N$ in each of the holonomy groups given in Remark \ref{groupes_dholonomie}.

The principle is roughly the following. We repeat the constructions of \cite{boubel2013a} Section 2, on a manifold $\MM$ with a nilpotent structure $N$, replacing everywhere real differentiable or complex holomorphic functions by nilomorphic $\R[\nu]$-valued, or nilomorphic+holomorphic $\C[\nu]$-valued ones. This gives Theorem \ref{realisation_nilomorphe}. To state it, if $J$ (or $L$) and $N$ are commuting (para)complex and nilpotent structures, we need integral coordinates for both $J$ (or $L$) and $N$. They exist, this is Lemma \ref{integrabilite_amelioree}, proven on p.\@ \pageref{preuve_integrabilite_amelioree}.

\begin{lem}\label{integrabilite_amelioree}(small enhancement of Lemma \ref{integration_N}) Suppose that $\R^{2d}$ is endowed with a (para)complex structure $J$ (or $L$) and a nilpotent structure $N$ commuting with it. Coordinates are called here adapted to $J$ (or $L$) if they make $\Mat(J)$ (or $\Mat(L)$) constant $J_1$- (or $L_1$-) block diagonal. Then there is a coordinate system simultaneously adapted to $J$ (or $L$) and integral for $N$ if and only if ${\mathcal N}_N=0$ and the $\ker N^k$ are involutive for all $k$. This holds also with complex coordinates if all those endomorphisms are $\C$-linear on $\C^{2d}$.
\end{lem}

\begin{te}\label{realisation_nilomorphe}{\bf (Corollary and generalization of Theorem \ref{theoremenilpotentadj})} 
 Let $H_{\goth s}$ be the generic holonomy group corresponding to an algebra $\goth s$ in any of the eight cases of Theorem \ref{structure_s}, and $N$ any self adjoint nilpotent endomorphism in the commutant of $\goth s$ {\em i.e.\@} the bicommutant of $H_{\goth{s}}$. We denote by $\GG_{H_{\goth{s}}^N}$ the set of germs of metrics whose holonomy group $H$ is included in the commutant $H_{\goth{s}}^N$ of $N$ in $H_{\goth{s}}$. If $g\in\GG_{H_{\goth{s}}^N}$, then $N$ extends as a parallel endomorphism, in particular as a nilpotent structure. Besides $\End(T\MM)^{{\goth h}}\supset\langle\goth s\cup\{ N\}\rangle$.\medskip
 
 \noindent{\bf (a)} For $\goth s$ in each case of Theorem \ref{structure_s}, $g\in\GG_{H_{\goth{s}}^N}$ if and only if:\medskip

-- In case {\bf (1$^\C$)}, $g$ is the real metic associated with the real part of a non degenerate, $\C[\nu]$-valued, $\C[\nu]$-bilinear metric $h$, which is both holomorphic and nilomorphic on $(\MM,\underline J,N)$.\medskip

-- In case {\bf (2)}, and in $N$-adapted coordinates $\bigl(x_i,(y_{i,a})_{a=1}^{n(i)-1}\bigr)_{i=1}^D$ such that $J\frac{\partial}{\partial x_{i}}=\frac{\partial}{\partial x_{i+1}}$ for any odd $i$, given by Lemma \ref{integrabilite_amelioree}, $g$ is the real metric associated with a non degenerate nilomorphic metric $h$ given by:
\begin{align*}
h\left({\textstyle\frac{\partial}{\partial w_i},\frac{\partial}{\partial \overline w_j}}\right)&=\frac{\partial^2u}{\partial w_i\partial \overline w_j},\begin{array}[t]{l}\text{with $u$ an $\R[\nu]$-valued nilomorphic function,}\\\text{and $\frac{\partial}{\partial w_{(j+1)/2}}:=\frac{\partial}{\partial x_j}-{\rm i}\frac{\partial}{\partial x_{j+1}}\in\T^{1,0}\R^d$}\\
\text{for all odd $j$,}\end{array}\\
&=\sum_\alpha\frac{\partial^{|\alpha|}}{\partial x_{\alpha}}\left(\frac{\partial^2\cu}{\partial w_i\partial \overline w_j}\right)(\nu y)^\alpha\text{\begin{tabular}[t]{l}with $\cu$ some $\R[\nu]$-valued\\adapted function of\\the coordinates $(x_i)_i$.\end{tabular}}
\end{align*}

-- In case {\bf (2')}, and in $N$-adapted coordinates such that $L\frac{\partial}{\partial x_{i}}=\frac{\partial}{\partial x_{i+1}}$ for any odd $i$, given by Lemma \ref{integrabilite_amelioree}, $g$ is the real metric associated with a non degenerate nilomorphic metric $h$ given by:
\begin{align*}
h\left({\textstyle\frac{\partial}{\partial x_i},\frac{\partial}{\partial x_{j+1}}}\right)&=\frac{\partial^2u}{\partial x_i\partial x_{j+1}}\begin{array}[t]{l}\text{\begin{tabular}[t]{l}for $i$, $j$ odd, with $u$ some $\R[\nu]$-valued\\nilomorphic function,\end{tabular}}\end{array}\\
&=\sum_\alpha\frac{\partial^{|\alpha|}}{\partial x_{\alpha}}\left(\frac{\partial^2\cu}{\partial x_i\partial x_{j+1}}\right)(\nu y)^\alpha\text{\begin{tabular}[t]{l}with $\cu$ some $\R[\nu]$-valued\\adapted function of\\the coordinates $(x_i)_i$.\end{tabular}}
\end{align*}

-- In case {\bf (2$^\C$)}, and in complex $N$-adapted coordinates such that $J\frac{\partial}{\partial x_{i}}=\frac{\partial}{\partial x_{i+1}}$, or $L\frac{\partial}{\partial x_{i}}=\frac{\partial}{\partial x_{i+1}}$, for $i$ odd, $g$ is the real part of the complex metric associated with a non degenerate complex nilomorphic metric $h$ given by the complexifications, which are equivalent, of any of the two formulas given in the previous cases. The potential $u$ is then a $\C[\nu]$-valued holomorphic and nilomorphic function.\medskip

-- In cases {\bf (3)}, {\bf (3')} and {\bf (3$^\C$)}, and in the real analytic category, it is the solution of the exterior differential system given in \cite{boubel2013a} Section 2, formulated in nilomorphic coordinates and with $\R[\nu]$- or $\C[\nu]$-valued nilomorphic functions instead of real or complex ones. In particular, in cases {\bf (3)} and {\bf (3')}, the elements of $\GG_{H_{\goth s}^N}$, considered up to diffeomorphism, are parametrized by $b\frac {D}2$ real analytic functions of $\frac D2+1$ variables, where $b:=\min\{a\in\llbracket1,n\rrbracket;d_a\neq 0\}$, with $D$ the number of invariant factors of $N$ and $d_a$ the number of repetitions of  $X^a$ among them. In case {\bf (3$^\C$)}, they are by $b\frac {D}4$ holomorphic functions of $\frac D4+1$ complex variables.\medskip 

\noindent{\bf (b)} In a dense open subset for the $C^2$ topology in ${\mathcal G}_{H_{\goth{s}}^N}$ (inside of the real analytic category in cases {\bf (3)}--{\bf (3')}--{\bf (3$^\C$)}), the holonomy group $H$ of the metric is exactly $H_{\goth s}^N$. So those commutants are holonomy groups.
\end{te}

\begin{rem} For the meaning of $D$ and the $d_a$, see also Notation \ref{notationNpointwise}.
\end{rem}

\begin{rem}\label{matrices_commutant} When does a nilpotent endomorphism $N$ commute with $\goth s$? Exactly when $\Id+N$ does. Now by definition, the automorphisms commuting with $\goth s$ are those preserving the $G$-structure defined by the commutant of $\goth s$. In cases {\bf (1$^\C$)} and {\bf (2)}, they are exactly the $J$- or $\underline J$-complex automorphisms, and in case {\bf (2')}, the paracomplex ones {\em i.e.\@} {\small$\biggl\{\left(\begin{array}{cc}U&0\\0&V\end{array}\right)\!\!$} with $U,V\!\!\in\!\!{\rm GL}_{d/2}(\R)\biggr\}$. In case {\bf (3)}, they are the $(J_1,J_2,J_3)$-quaternionic automorphisms. In case {\bf (3')},  if $\Mat(L)=I_{d/2,d/2}$ and $\Mat(J)=J_{d/2}$, they are {\small$\biggl\{\left(\begin{array}{cc}U&0\\0&-JUJ\end{array}\right)$} with $U\in{\rm GL}_{d/2}(\R)\biggr\}$. Cases  {\bf (2$^\C$)} and {\bf (3$^\C$)} are the complexification of  cases {\bf (\mathversion{bold}$2$\mathversion{normal})}--{\bf (\mathversion{bold}2$^{\prime}$\mathversion{normal})} and {\bf (\mathversion{bold}$3$\mathversion{normal})}--{\bf (\mathversion{bold}3$^{\prime}$\mathversion{normal})}.
\end{rem}

Now, in the different cases of Theorem \ref{realisation_nilomorphe}, Corollary \ref{cor_forme_e} gives the form of $\goth{e}:=\End(T\MM)^{\goth h}$, and notably that of its radical $\goth n$, asked for in the introduction. Here I thank the referee for his noticing a mistake in my first version. The commutant ${\goth c}:=\End(T\MM)^{\left(H_{\goth{s}}^N\right)}$ of $H_{\goth s}^N$ contains $\goth s$ and the bicommutant $N^{\text{\sf cc}}:=\End(T\MM)^{\left(\End(T\MM)^N\right)}$ of $N$ (classically equal to $\R[N]$). In fact $\goth{c}=\langle{\goth s}\cup\{N\}\rangle$, except a bit surprisingly in some sub-cases of {\bf (1)} and {\bf (1$^\C$)}. This could be expected: look at the lowest possible dimension in the case $N=0$. For $\K\in\{\R,\C\}$, $\goth{o}(g)^{\goth s}=\goth{o}_1(\K)=\{0\}$ is trivial in cases {\bf (1)}--{\bf (1$^{\C}$)}, so $\End_\K(\T\MM)=\goth{c}\supsetneq{\goth o}(g)$, whereas in the other cases $\goth{o}(g)^{\goth s}\neq\{0\}$, {\em e.g.\@} in case {\bf (2)}, $\End(\T\MM)\supsetneq{\goth c}=\goth{u}_1=\langle J\rangle$.

Corollary \ref{cor_forme_e} gives the details. It follows from a sequence of results in linear algebra, gathered at the end of this section, pp.\@ \pageref{alg_lin} {\em sq.} To know the general form of the matrices at stake in it, or related to algebras involved in this article ---~commutants {\em etc.\@}~---, notably in the delicate cases {\bf (1)} and {\bf (1$^\C$)}, see Lemma \ref{lemme_commutants}. Replace in it the real coefficients by complex ones to get cases {\bf (1\mathversion{bold}$^\C$\mathversion{normal})}--{\bf (2$^\C$)}--{\bf (3$^\C$)}.

\begin{cor}\label{cor_forme_e}Suppose that $H$ is the holonomy group of a metric $g$ and that  $H\subset H_{\goth s}^N$ (see Theorem \ref{realisation_nilomorphe}). Then immediately, $\End(T\MM)^H\supset\End(T\MM)^{H_{\goth s}^N}$, with equality if $H=H_{\goth s}^N$ (which holds generically). Now:\medskip

{\bf (a)} If $\goth s$ is not in case {\bf (1)} or {\bf (1$^\C$)} of Theorem \ref{structure_s} {\em i.e.\@} ${\goth s}\neq\R.\Id$ and ${\goth s}\neq\langle \underline J\rangle$ with $\underline J$ a self adjoint complex structure then:
\begin{gather*}
\End(T\MM)^{H_{\goth s}^N}=\langle{\goth s}\cup\{N\}\rangle=\goth{s}\oplus(N),
\end{gather*}
so if $H=H_{\goth s}^N$, $(\MM,g)$ is indecomposable. The sum $\goth{s}\oplus(N)$ is the decomposition of $\langle{\goth s}\cup\{N\}\rangle$ into a semi-simple part and its radical; $(N)$ is the radical spanned by $N$ in $\End(T\MM)^{H_{\goth s}^N}$.\medskip

{\bf (b)} If $\goth s$ is in case {\bf (1)} or {\bf (1$^\C$)} of Theorem \ref{structure_s}, set $n=n_1\geqslant n_2\geqslant\ldots\geqslant n_D$ the sizes of the Jordan blocks of $N$ ---~viewed as a $\underline J$-complex endomorphism in case {\bf (1$^\C$)}. By convention, we set $n_i=0$ for $i>D$.\medskip

{\bf (b1)} If $n_1>2n_2$, locally, $(\MM,g)\simeq(\MM',g')\times(\R^{n_1-2n_2},g_{\text{flat}})$ is decomposable. The sizes of the Jordan blocks of $N':=N_{|T\MM'}$ are $(2n_2,\linebreak[3]n_2,\ldots,\linebreak[2]n_D)$. So the situation on $(\MM',g')$ is described by {\bf (b2)}.\medskip

{\bf (b2)} Else, set $\K:=\R$ or $\K:=\R[\underline J]\simeq \C$, so that $\K\simeq{\goth s}$. Set $\pi':T\MM\rightarrow T\MM/\ker N^{n_2}$, $\pi'':T\MM\rightarrow T\MM/\ker N^{n_3}$ and $g_{n_3}:=g(\,\cdot\,,N^{n_3}\,\cdot\,)$ [or its complexification in case $\K=\C$], defined on $\pi''(T\MM)$. Then:
\begin{gather*}\End(T\MM)^{H_{\goth s}^N}=
\left(\K[N]+(N^{n_2})_\ast\pi^{\prime\ast}\End_\K\bigl(T\MM/\ker N^{n_2}\bigr)\right)\\
\qquad\qquad\qquad\qquad\oplus (N^{n_3})_\ast\pi^{\prime\prime\ast}\goth{o}_\K(g_{n_3})^{N}
\end{gather*}
where $\varphi_\ast\psi^\ast A\!:=\!\{\varphi\circ f \circ\psi,f\in A\}$. So $(\MM,g)$ is indecomposable if $H=H_{\goth s}^N$. A semi-simple part of $\End(T\MM)^{H_{\goth s}^N}$ is $\goth{s}=\K.\Id\subset\K[N]$ and its radical is
$\left(N\K[N]+(N^{n_2})_\ast\pi^{\prime\ast}\End_\K\bigl(T\MM/\ker N^{n_2}\bigr)\right)\oplus (N^{n_3})_\ast\pi^{\prime\prime\ast}\goth{o}_\K(g_{n_3})^{N}$.

If $n_1=n_2$ and $n_3=0$, $\pi'=0$ and $\pi''=\Id$ so: $\End(T\MM)^{H_{\goth s}^N}=\K[N]\oplus \goth{o}_\K(g)^{N}$. In this case, $\goth{o}_\K(g)^{N}$ contains $J$ or $L$ a skew adjoint (para)complex structure commuting with $N$, so
semi-simple part $\goth{s}'$ of $\End(T\MM)^{H_{\goth s}^N}$ is {\em not} $\goth{s}\simeq\K\Id$, but $\goth{s}'=\K.J$ or $\goth{s}'=\K.L$, $\End(T\MM)^{H_{\goth s}^N}=\goth s'\oplus(N)$ and we are respectively in case {\bf (2)} or {\bf (2')} of Theorem \ref{structure_s} if $\K=\R$, and in case {\bf (2$^\C$)} if $\K=\C$. 

Moreover, the signature of the flat metric $g_{\text{flat}}$ appearing in the first point is$ \left(\frac{d''}2,\frac{d''}2\right)$ or $\left(\frac{d''\pm 1}2,\frac{d''\mp1}2\right)$ with $d'':=n_1-2n_2$, according to the cases given in Proposition \ref{forme_e}.
\end{cor}

\noindent{\bf Proof.} Apply Proposition \ref{forme_e}, complexified for cases {\bf (1$^\C$)}--{\bf (2$^\C$)}--{\bf (3$^\C$)}. For the last, exceptional case, see the second point of Remark \ref{cas_elementaires}. \hfill{\rm q.e.d.}

\begin{importantrem}\label{indecomposabilite}Corollary \ref{cor_forme_e} shows that in Theorem \ref{realisation_nilomorphe}, gen\-eric metrics are indecomposable, except in cases {\bf (1)}--{\bf (1$^\C$)} when $n_1>2n_2$.
\end{importantrem}

\begin{rem}In Corollary \ref{cor_forme_e} {\bf (b2)}, the situation is like in {\bf (a)} {\em i.e.\@} $\End(T\MM)^{H_{\goth s}^N}=\K[N]=\K\Id\oplus(N)$ if and only if $\pi'=\pi''=0$ {\em i.e.\@} $n_1=n_2=n_3$ {\em i.e.\@} $N$ has at least three Jordan blocks of maximal size.
\end{rem}

\begin{rem}\label{action_triviale}
The cases where $T\MM/\ker N^{n_2}\neq\{0\}$ {\em i.e.\@} $(n=)n_1>n_2$ (considering, in case {\bf (1$^\C$)}, the complex Jordan blocks) are exactly those where the commutation with $N$ (or with $\{N,\underline J\}$ in case {\bf (1$^\C$)}) forces the holonomy group to act trivially on a non null subspace, namely $\im N^{n_2}$. 
\end{rem}

\begin{rem}\label{cas_elementaires}Some very elementary cases turn out to be ``exceptional'' in the classification of Corollary \ref{cor_forme_e}.\medskip

-- If $N$ consists of only one Jordan block, so of order $n=d>0$, then immediately $\goth{h}\subset\goth{o}(g)^N=\{0\}$, do the calculation or use Lemma \ref{lemme_commutants} {\bf (b)} and {\bf (c)}. So the metric is flat and $\End(T\MM)^{\goth{h}}=\End(T\MM)$. This is case {\bf (1)} with $n_2=0$ so $2n_1>n_2$, and $\End_\R\bigl(T\MM/\ker N^{n_2}\bigr)=\End(T\MM)$.\medskip

-- If $N$ consists of exactly two Jordan blocks, of the same size $n=\frac d2>0$, then a skew adjoint complex structure $J$, if $\sign(g_{n-1})\in\{(2,0),(0,2)\}$, or paracomplex structure $L$, if $\sign(g_{n-1})=(1,1)$, commuting with $N$, is also parallel, and we are in case {\bf\mathversion{bold}($2$)\mathversion{normal}} or {\bf\mathversion{bold}($2^{\prime}$)\mathversion{normal}}. Then $N$, as a (para)complex endomorphism, has one single Jordan block and $\End(T\MM)^{\goth{h}}=\goth{s}\R[N]$ with $\goth{s}=\langle J\rangle$ or $\goth{s}=\langle L\rangle$. Apply case {\bf (1)} of Corollary \ref{cor_forme_e} with $n_1=n_2>0$ and $n_3=0$; then $g_{n_3}=g$ and $\End(T\MM)^{\goth{h}}=\R[N]\oplus\goth{o}(g)^{N}$.\medskip

-- This does not go on: if $N$ consists of $k>3$ Jordan blocks of the same size $n=\frac dk>0$, $\End(T\MM)^{\goth{h}}=\K[N]$ is the bicommutant of $N$.\medskip

-- All this appears in the case $N=0$. For $g$ generic, $\End(T\MM)^{\goth{h}}=\End(T\MM)^{\goth{o}(g)}$. If $d\geqslant 2$, ${\goth{o}(g)}$ is abelian and $\End(T\MM)^{\goth{h}}$ is the (Euclidian or Lorentzian) conformal group. If $d>2$, the commutant of ${\goth{o}(g)}$ is trivial, $\End(T\MM)^{\goth{h}}=\R.\Id$. You may see this difference appear in case {\bf (1)}, $(n_1,n_2,n_3)=(1,1,0)$ if $d=2$ and $(n_1,n_2,n_3)=(1,1,1)$ if $d>2$.
\end{rem}

\noindent{\bf Proof of Lemma \ref{integrabilite_amelioree}.}\label{preuve_integrabilite_amelioree} We do it in the real case. The ``only if'' is immediate. The converse is immediate in the case of $J$: repeat the  proof of \cite{thompson} with the field $\C$ replacing $\R$: $(\R^{2d},J)\simeq\C^d$, so work with complex coordinates. For the case of $L$, or of $J$ for an alternative proof, we have to check that the proof works with $L$-adapted coordinates at each step (on p.\@ 610 of \cite{thompson}). The author builds the coordinates by induction on the nilpotence index $n$ of $N$. For $n=1$ the result is empty hence true. If it holds for index $n-1$ and if $N^{n-1}\neq N^n=0$, as $\ker N$ is $L$-invariant we may find coordinates $((x_i)_i,(y_i)_i)$ that are:\medskip

-- $L$-adapted {\em i.e.\@}: $L\frac{\partial}{\partial x_{i}}=\frac{\partial}{\partial x_{i+1}}$ and $L\frac{\partial}{\partial y_{i}}=\frac{\partial}{\partial y_{i+1}}$ for any odd $i$,\medskip

-- such that the $(y_i)_i$ parametrize the leaves of $\KK$.\medskip

\noindent Then the induction assumption applies on $\R^{2d}/\KK$, providing coordinates $(x_i)_i$ of the wished type on $\R^{2d}/\KK$. As the fields $\left(N\frac{\partial}{\partial x_{i}}\right)_i$ commute with each other, \cite{thompson} extends these coordinates to the whole $\R^{2d}$, obtaining $N$-adapted coordinates $((x_i)_i,(\ybar_i)_i)$. Here, we need moreover to check that they are also $L$-adapted {\em i.e.} {\bf (i)} $L\left(\frac{\partial}{\partial \xbar_{i}}\right)=\frac{\partial}{\partial\xbar_{i+1}}$ and {\bf (ii)} $L\left(\frac{\partial}{\partial \ybar_{i}}\right)=\frac{\partial}{\partial\ybar_{i+1}}$ for $i$ odd. We follow \cite{thompson}: the $\frac{\partial}{\partial \xbar_{i}}$ are equal to some $N\left(\frac{\partial}{\partial x_{j}}\right)$ if they are in $\im N$, else they are equal to $\frac{\partial}{\partial x_{i}}$. As $LN=NL$, this gives {\bf (i)}. The  $\frac{\partial}{\partial \ybar_{i}}$ are equal to some $N\left(\frac{\partial}{\partial x_{k}}\right)$ if they are in $\im N$, else they are chosen freely. As $LN=NL$, this gives {\bf (ii)}.\hfill{\rm q.e.d.}

\begin{rem}\label{integrabilite_amelioree_bis} {\em (This will be used in Part \ref{nonprincipal})} The key properties used in the proof are that $J$, or $L$, has a constant matrix in the basis $\left(\frac{\partial}{\partial x_{i}}\right)_i$, and that it commutes with $N$. So the same proof, and result, hold for any integrable field of endomorphism playing the role of $J$ or $L$.
\end{rem}

To show the theorem we need the following remarks.

\begin{rem}\label{omega_ponctuelle}{\bf(\mathversion{bold}natural matricial form for an $\R[\nu]$-bilinear alternate 2-form)\mathversion{normal}} Using Notation \ref{notationNpointwise} and \ref{remNpointwise}, we adapt  \ref{formes_quotient}, \ref{defbilinearassociated} and \ref{matrice_h_en_general} for $\omega$ an $\R[\nu]$-bilinear {\em alternate} 2-form on $E=(\R^d,N)$.\vspace{-.5ex}
\begin{gather*}
\omega=\sum_{a=0}^{n-1}\nu^a\omega_a\ \begin{array}[t]{l}\text{with $\omega_{n-1}$ such that }\omega_{n-1}(\,\cdot\,,N\,\cdot\,)=\omega_{n-1}(N\,\cdot\,,\,\cdot\,)\\
\text{and for any $a$, }\omega_{a}=\omega_{n-1}(\,\cdot\,,N^{n-1-a}\,\cdot\,).\end{array}
\end{gather*}
For $a\in\llbracket1,n\rrbracket$ we denote by $r_a\in2\N$ the rank of $\omega_{n-a}:=\omega_{n-1}(\,\cdot\,,N^{a-1}\,\cdot\,)$ defined on $\ker N^a/\ker N^{a-1}$. It is standard that the couple $(N,\omega_{n-1})$ on $\R^d$ is characterized up to conjugation by the $(r_a)_{a=1}^n$, called here the {\em ranks of $\omega$}. See {\em e.g.\@} \cite{leep-schueller} like in \ref{formes_quotient}.

If $\beta=(X_i)_{i=1}^D$ is an adapted spanning family (see \ref{remNpointwise}) of $E$, $\Mat_\beta(\omega)=\sum_{a=0}^{n-1}\nu^a\Omega_a\in{\rm M}_D(\R[\nu])$ where:\medskip

{\bf (i)} $\Omega_a=\text{\small $\left(\begin{array}{cc}0&0\\0&\cOmega^a\end{array}\right)$}$, the upper left null square block, of size $D_{n-1-a}$, corresponding to $\Span_{\R[\nu]}\bigl\{X_i;N^{n-1-a}X_i=0\bigr\}$,\medskip

{\bf (ii)}  the upper left block $\cOmega^a_0$ of $\cOmega^a$ of size $d_{n-a}$, corresponding to $\Span_{\R[\nu]} \bigl\{X_i;N^{n-1-a}X_i\neq N^{n-a}X_i=0 \bigr\}$ is of rank $r_{n-a}$. So if $S\oplus \im N=E$, $r_a$ is the rank of the (well defined) from $\omega_{n-a}$ on the quotient $(S\cap\ker N^a)/(S\cap\ker N^{a-1})$.\medskip

\noindent For $r\leqslant \delta$ and $r$ even, we denote by $J_{\delta,r/2}$ the matrix $\diag(J_{r/2},0)\in{\rm M}_\delta(\R[\nu])$. There are adapted spanning families $\beta=(X_i)_{i=1}^D$ of $E$ such that $\cOmega_a$ is null except $\cOmega^a_0=J_{d_a,r_a}$, for all $a$ {\em i.e.\@}:
\begin{align*}
\Mat_\beta(\omega)&=\diag\left(\nu^{n-a}J_{d_a,r_a/2}\right)_{a=1}^n\\&=\left(\begin{array}{ccc}\nu^{n-1}J_{d_{1},r_{1}/2}\\&\ddots\\&&\nu^{0}J_{d_{n},r_{n}/2}\end{array}\right)\in{\rm M}_D(\R[\nu]).
\end{align*}
\noindent Each block $\nu^{n-a}J_{d_a,r_a/2}$ corresponds to the factor $(\nu^{n-a}\R[\nu])^{d_{a}}=%JDG
\linebreak[4]\Span\bigl((X_i)_{D_{a-1}<i\leqslant D_a}\bigr)$ of $E$. The form $\omega$ is non degenerate if and only if each $\omega_a$ is {\em i.e.\@} $r_a=d_a$ for all $a$.
\end{rem}

\begin{lem}\label{darboux-poincare}The Poincar\'e and Darboux lemmas admit a natural ``nilomorphic'' version. For example, if $B$ is some ball in $(\R^d,N)$ with $N$ nilpotent, in constant Jordan form:\medskip

{\bf (a)} If $\lambda\in\Lambda^k_{\R[\nu]}(B)$ with $p>0$ is a closed nilomorphic k-form on $B$ then there is a nilomorphic $\alpha\in\Lambda^{k-1}_{\R[\nu]}(B)$ such that $\lambda=\dd\alpha$.\medskip

{\bf (b)} If $\omega\in\Lambda^2_{\R[\nu]}(B)$ is closed and has constant ranks, there exist nilomorphic coordinates on $B$ in which $\Mat(\omega)$ has the (constant) form given at the  end of Remark \ref{omega_ponctuelle}.
\end{lem}

\noindent{\bf Proof.} {\bf(a)} Classically, if $\lambda$ is a closed real form and $X$ a vector field on $B$ whose flow $(\varphi^t)_{t\in[0,+\infty[}$ is a retraction of $B$ on a point, then:
$$\alpha:=\int_0^\infty\iota_X\left(\varphi^{t\ast}\lambda\right)\dd t$$
fits. In our case, use the retraction $(e^{-t}\Id_{B})_{t\in[0,+\infty[}$ generated by $X=-\Id_{\R^d}$; $\varphi^t$ and $X$ being nilomorphic, so is the obtained integral form $\alpha$.

{\bf (b)} As $\omega$ is closed and has constant ranks, its kernel integrates in some foliation $\FF$ and $N$ acts on the quotient $B/\FF$. So we may suppose that $\omega$ is non degenerate. Then we use Moser's path method. We set $\omega_0$ the constant nondegenerate 2-form on $\R^d$, given in Remark \ref{omega_ponctuelle}. To simplify, we suppose that $\omega_t:=\omega_0+t(\omega-\omega_0)$ never degenerates. Else, iterate the method along an adequate piecewise affine path form $\omega_0$ to $\omega$. We want to build a {\em nilomorphic} homotopy $(\varphi^t)_{t\in[0,1]}$ such that:\medskip

\noindent\phantom{$(\ast)$}\hspace*{\fill}$\varphi_t^\ast\omega_t=\omega_0$.\hspace*{\fill}$(\ast)$\medskip

\noindent Once this is done, $\varphi_1^\ast\omega=\omega_0$ and $\varphi_1^\ast N=N$ as we want. Let $X_t$ be the field such that $X_t(\varphi^t(p))=\frac{\dd}{\dd t}\varphi^t(p)$, then $(\ast)$ amounts to:
$$\varphi^{t\ast}\left({\mathcal L}_{X_t}\omega_t+\frac{\dd}{\dd t}\omega_t\right)=0$$
{\em i.e.\@}, as $\dd\omega_t=0$, $\dd(\iota_{X_t}\omega_t)+\omega=0$. By {\bf (a)}, there is a {\em nilomorphic} 1-form $\lambda$ on $B$ with $\omega=\dd\lambda$. Then it is sufficient to find $X_t$ with $\iota_{X_t}\omega_t+\lambda=0$. As $\omega_t$ is nondegenerate, this defines indeed $X_t$; as $\omega_t$ and $\lambda$ are nilomorphic, so is $X_t$.\hfill{\rm q.e.d.}

\begin{rem}\label{onefunctionofDvariables}An infinitesimal deformation of nilomorphic Darboux coordinates for a non degenerate $\omega$ is a nilomorphic field $X$ such that $0={\mathcal L}_X\omega=\dd(\iota_X\omega)$. So classically, $X$ is the symplectic gradient of some (nilomorphic) potential $f$. By Prop.\@ \ref{development_form}, $f$ is given by its ``adapted'' restriction $f=\sum_af_a\nu^a$ to $\TT=\{(\nu y)=0\}$ {\em i.e.\@} by one real function 
 ($f_{n-1}$) of $D-D_0=D$ variables, one ($f_{n-2}$) of $D-D_{1}$ variables {\em etc.\@}, and one ($f_{0}$) of $D-D_{n-1}$ variables. Set $b:=\min\{a\in\llbracket1,n\rrbracket;d_a\neq 0\}$. Then $0=D_0=\ldots=D_{b-1}$. So $f_{n-1}$ up to $f_{n-b}$ depend, each, on $D-D_{0}=D$ variables ---~strictly more variables than any other $f_c$. So the Darboux coordinates depend on  $b$ functions of $D$ variables. 

In passing, we add the following. After Rem.\@ \ref{formes_quotient}, on each leaf of the (quotient) foliation $\pi(\KK^{a+1})/\pi(\KK^a)$, the (real) symplectic form $\omega(\,\cdot\,,N^a\,\cdot\,)$ is well defined. Choosing nilomorphic Darboux coordinates means in particular choosing Darboux coordinates for them, but also choosing a transversal $\TT=\{(\nu y)=0\}$ to $\II$ such that for all $a$:\medskip

-- the orthogonal distribution to $\TT\cap\KK^a$ with respect to $\omega(\,\cdot\,,N^{a-1}\,\cdot\,)$ is integrable,\medskip

-- its intersection with $\TT$ is totally isotropic.\medskip

\noindent In fact, it amounts exactly to both these choices.
\end{rem}

\noindent{\bf Proof of Theorem \ref{realisation_nilomorphe}.}  The fact that $N$ extends as a nilpotent structure was given by Lemma \ref{integration_N} in \S\ref{germsN}. To prove {\bf (a)} and {\bf (b)}, we begin with cases {\bf (3)} and {\bf (3')} ---~case {\bf (3$^\C$)} is only their complexification.\medskip

{\bf Part (a)}  We follow the line of \cite{boubel2013a} Section 2. The triple $(J,U,N)$ is given and we look for a quadruple $(g,J,U,N)$. This is equivalent to a quadruple $(J,N,\widetilde\omega_0,\widetilde\omega)$ where:\vspace{-1ex}
$$\widetilde\omega_0:=\sum_{a=0}^{n-1}\omega_0(\,\cdot\,,N^{n-1-a}\,\cdot\,)\nu^a$$
is the ($J$)-complex nilomorphic symplectic 2-form associated with the complex symplectic 2-form $\omega_0:=g(\,\cdot\,,U\,\cdot\,)+{\rm i}g(\,\cdot\,,JU\,\cdot\,)$, and where:
$$\widetilde\omega:=\sum_{a=0}^{n-1}\omega(\,\cdot\,,N^{n-1-a}\,\cdot\,)\nu^a$$
is the symplectic nilomorphic (1,1)-form associated with $\omega:=\omega_0(\,\cdot\,,U\,\cdot\,)$. As $UN=NU$, $\widetilde\omega=\widetilde\omega_0(\,\cdot\,,U\,\cdot\,)$. By Lemma \ref{darboux-poincare} {\bf (b)}, we may use local coordinates $(x_j,(y_{j,a})_a)_{j=1}^{D}$ adapted to $N$ and such that for all odd $j$, $J\frac{\partial}{\partial x_j}=\frac{\partial}{\partial x_{j+1}}$. As usual, we set $z_j:=x_j+(\nu y_j)\in\R[\nu]$ for all $j$ and introduce the ``complex and nilomorphic'' local coordinates $(w_j)_{j=1}^{D/2}$ by:
$$w_{(j+1)/2}:=z_j+{\rm i}z_{j+1}\in\C[\nu].$$
Through those coordinates, we consider that we are in some ball $\mathcal B$ of $(\C^{d/2},N)$ with $N$ a complex endomorphism field. The matrix $\Omega_0$ of $\omega_0$ as a $\C[\nu]$-bilinear form is as at the end of Remark \ref{omega_ponctuelle}, with non degenerate $J_{d_a,r_a/2}=J_{d_a}$.  In the complexification $T^\C\BB=T{\mathcal B}\otimes\C$ of the tangent bundle, we introduce also, for $j$ odd, the vector fields $\frac{\partial}{\partial w_{(j+1)/2}}:=\frac{\partial}{\partial x_j}-{\rm i}\frac{\partial}{\partial x_{j+1}}$ and $\frac{\partial}{\partial \overline w_{(j+1)/2}}:=\frac{\partial}{\partial x_j}+{\rm i}\frac{\partial}{\partial x_{j+1}}$. We introduce the matrix $V=(v_{i,j})_{i,j=1}^{D/2}\in{\rm M}_{D/2}(\C[\nu])$ of the $\nu$-linear, $J$-antilinear morphism $U$ by:
\begin{gather*}
U\left(\frac{\partial}{\partial w_i}\right)=\sum_{j=1}^{D/2}v_{i,j}\frac{\partial}{\partial \overline w_j}\quad\text{{\em i.e.\@}, if $v_{i,j}=\sum_{a=0}^{n-1} v_{i,j,a}\nu^a$, then:}\\
\forall b\in\llbracket0,n-1\rrbracket, U\left(N^b\frac{\partial}{\partial w_i}\right)=\sum_{j,a}v_{i,j,a} N^{b+a}\frac{\partial}{\partial \overline w_j}.
\end{gather*}
Notice that, as $N^a\frac{\partial}{\partial \overline w_j}=0$ if $a\geqslant n(j)$, $\deg_\nu v_{i,j}<n(j)$ for all $i,j$. In other terms, in the block decomposition of $V$ corresponding to the flag of the $\pi(\ker N^a)$, {\em i.e.\@} to the blocks of $\Mat(\omega)$ given in Remark \ref{omega_ponctuelle}, the $\nu$-degree of the $a^{\text{th}}$ line of blocks is strictly less than $a$.

We keep following \cite{boubel2013a} Section 2. Let $\dd w$ be the column $(\dd w_j)_{j=1}^{D/2}$, then:
$$\omega_0={}^t\!\dd w\wedge\Omega_0\wedge\dd w\quad\text{and, setting }H:=-\Omega_0V,\ \ \omega=\frac{\rm i}2{}^t\!\dd w\wedge H\wedge\dd \overline w.$$
Here, notice that $\Phi:V\mapsto-\Omega_0V=H$ is injective. Indeed, it associates with the matrix $V$ of $U$, the matrix $H$ of the form $\omega_0(\,\cdot\,,U\,\cdot\,)$, and $\omega_0$ is non degenerate. An alternative, computational argument is the fact that, as $\deg_\nu v_{i,j}<n(j)$ for all coefficient $v_{i,j}$ of $V$, no product $\nu^a\nu^b$ with $a+b\geqslant n$ appears when computing $\Omega_0V$. This shows also the following. Introduce:
$$\Omega_0^{-1}:=-\diag\left(\nu^{-n+1}J_{d_1/2},\linebreak[1]\nu^{-n+2}J_{d_2/2},\ldots,\linebreak[1]\nu^{-1}J_{d_{n-1}/2},\linebreak[1]J_{d_n/2}\right),$$
$\nu^{-a}$ standing for the application $\sum_{b\geqslant a}f_b\nu^b\mapsto\sum_{b\geqslant a}f_b\nu^{b-a}$, from $\nu^a\C[\nu]$ to $\C[\nu]$  {\em i.e.\@} for the left inverse of the multiplication by $\nu^a$. Then $H\mapsto-\Omega_0^{-1}H$ is well-defined on the space of the matrices of $\R[\nu]$-bilinear forms, as they satisfy property {\bf (i)} of Remark \ref{omega_ponctuelle}. More precisely, $\Phi$ is a bijection:
\begin{gather*}\hspace*{-3em}\Phi:\{V=(v_{i,j})\in M_D(\R[\nu]);\deg_\nu v_{i,j}<n(j)\ \text{for all}\ i,j,\\\qquad\qquad\qquad\qquad\qquad\qquad^t\overline V\Omega_0=-\Omega_0V\ \text{and}\ V\overline V=\varepsilon I\}\rightarrow\\
\hspace*{-3em}\widetilde\HHH_\varepsilon:=\{H\in M_D(\R[\nu]);H\ \text{satisfy property {\bf (i)} of \ref{omega_ponctuelle},}\\\qquad\qquad\qquad\qquad\qquad\qquad^t\overline H=H\ \text{and}\ H\Omega_0^{-1}\overline H=\varepsilon\Omega_0\}.
\end{gather*}
Notice that in the set above, $H$ is necessarily non degenerate, so satisfies property {\bf (ii)} of \ref{omega_ponctuelle} with $r_a=d_a$ for all $a$. Now that this adaptation to the nilomorphic case is done, we may perform Cartan's test (see \cite{bcggg,ivey-landsberg}). In fact, it works just like in \cite{boubel2013a} Section 2. Indeed:\medskip

-- We look for an $N$-stable integral manifold of the exterior differential equation ${\bf I}:\mbox{$^t\!\dd w$}\wedge\dd H\wedge \dd \overline w=0$ {\em i.e.\@} for a nilomorphic function $H:(\C^{d/2},N)\rightarrow\widetilde\HHH_\varepsilon$ around the origin, whose graph is an integral manifold of ${\bf I}$. By Theorem \ref{development}, it amounts to find an {\em adapted} function $H:\TT\rightarrow\widetilde\HHH_\varepsilon$ (see \ref{defadapted}), with $\TT$ the transversal $\{(\nu y)=0\}$ to the foliation $\II$. So in the following we work on $\TT$, identified with $ \C^{D/2}$ by the coordinates.\medskip

-- Then, Cartan's test rests on the equation of the tangent space $\widetilde W_\varepsilon:=T_{H(0)}\widetilde\HHH_\varepsilon$, which we will see to be nearly the same as in \cite{boubel2013a} Section 2. More precisely, along $\TT$, the ${\rm M}_D(\C[\nu])$-valued function $H$ we look for reads $H=\sum_{a=0}^{n-1}H_a\nu^a$, each $H_a$ being the pull back of some complex valued matrix function on $\TT/(\TT\cap\KK^{n-1-a})$. We will see that the coefficient of  $\nu^a$ in {\bf I} is an exterior differential equation involving only the function $H_a$, thus is an exterior differential equation defined on $\TT/(\TT\cap\KK^{n-1-a})$. Each of those equations is like that of \cite{boubel2013a} Section 2.\medskip

To alleviate the formulas, we write them in the case $\varepsilon=-1$ and $g$ positive definite. The other cases work alike. We introduce $\widehat \Omega_0:=\diag(J_{d_1/2},\linebreak[1]J_{d_2/2},\ldots,\linebreak[1]J_{d_n/2})$. As in \cite{boubel2013a} Section 2, at the origin, we may take $V=\widehat \Omega_0$ {\em i.e.\@}, at the origin:
$$H=\widehat I:=\diag\left(\nu^{n-1}I_{d_1/2},\linebreak[1]\nu^{n-2}I_{d_2/2},\ldots,\linebreak[1]\nu I_{d_{n-1}/2},\linebreak[1]I_{d_n/2}\right),\vspace{-2ex}$$
\begin{align*}
\text{then: }W\in\widetilde W_\varepsilon&\Leftrightarrow W\Omega_0^{-1}\widehat I+\widehat I\Omega_0^{-1}W=0\\
&\Leftrightarrow W\widehat\Omega_0+\widehat \Omega_0W=0\\
&\Leftrightarrow \sum_{a=0}^{n-1}\nu^a(W_a\widehat\Omega_0+\widehat \Omega_0W_a)=0.
\end{align*}
 The coefficient of  $\nu^a$ involves only $W_a$. So, as announced, the coefficient of  $\nu^a$ in {\bf I} is an equation involving only $H_a$. This equation is stated on the $(D-D_{n-1-a})/2$-dimensional quotient $\TT/(\TT\cap\KK^{n-1-a})$ as we look for a $\KK^{n-1-a}$-basic function $H_a$. (This is consistent with the fact that only the bottom right square of $W_a$, appearing in {\bf (i)} of Remark \ref{omega_ponctuelle} and corresponding  to the quotient by $\ker N^{n-1-a}$, is non vanishing.) Now on this quotient, the vectors may be reordered so that $\widehat\Omega_0$ reads {\scriptsize$\left(\begin{array}{cc}0&I\\-I&0\end{array}\right)$},  and $W_a\widehat\Omega_0+\widehat \Omega_0W_a=0$ is the same equation as that defining $W_\varepsilon$ in \cite{boubel2013a} Section 2. So Cartan's criterion is fulfilled and the solutions $H_a$ depend on $(D-D_{n-1-a})/2$ of $(D-D_{n-1-a})/2+1$ variables. Finally, set $b:=\min\{a\in\llbracket1,n\rrbracket;d_a\neq 0\}$. Then $0=D_0=\ldots=D_{b-1}$. So $H_{n-1}$ up to $H_{n-b}$ depend, each, on $(D-D_{0})/2=\frac D2$ functions of $(D-D_{0})/2+1=\frac D2+1$ variables ---~strictly more variables than any other $H_c$. So the whole function $H$ depends on  $b\frac D2$ functions of $(D-D_{0})/2+1=\frac D2+1$ variables. By Remark \ref{onefunctionofDvariables}, the choice of the complex Darboux coordinates for $\omega$ amounts to that of $b$ function of $D/2$ variables, so this does not interfer.\medskip

\noindent{\em Note.}\label{nilomorphisation} The above technique may be used as a standard reasoning to adapt arguments about an exterior differential system to the nilomorphic framework.\medskip

\noindent{\bf Part (b) for cases {\bf (3)}--{\bf (3')}--{\bf (3$^\C$)}.} Following the note just above, {\bf (b)} is given by the reasoning referred to in Remark 2.11 of \cite{boubel2013a}, adapted to the nilomorphic case {\em i.e.\@} applied to adapted functions, or jets, defined on $\TT$.\medskip

\noindent {\bf Part (a) for case {\bf (1$^\C$)}.} Use the first point of the Reminder at the beginning of \cite{boubel2013a} \S2, and repeat the very proof of Theorem \ref{theoremenilpotentadj}, with the field $\C$ and holomorphic functions replacing $\R$ and smooth functions.\medskip

\noindent {\bf Part (a) for cases {\bf (2)} and {\bf (2')}}, hence {\bf (2$^\C$)}, which is their complexification. We are now quicker. Repeat the classic proofs (see respectively {\em e.g.\@} \cite{moroianu} \S11.2 and \S8.3, and \cite{BB-ikemakhen1997} \S2) with nilomorphic coordinates and functions replacing real ones. In other words:\medskip

-- Take a nilomorphic coordinate system which is also integral for $J$ or $L$ {\em i.e.\@} such that $J\frac{\partial}{\partial x_i}=\frac{\partial}{\partial x_{i+1}}$ or $L\frac{\partial}{\partial x_i}=\frac{\partial}{\partial x_{i+1}}$ for all odd $i$. This is given by Lemma \ref{integrabilite_amelioree}.\medskip

-- Set $\TT:=\{(\nu y)=0\}=\{\forall i,(\nu y_i)=0\}$, then apply these proofs along $\TT$, to each of the (real) coefficients $f_a$, factor of $\nu^a$, of the nilomorphic functions $f$ that appear. As the latter are adapted along $\TT$ (see \ref{defadapted}), this means applying the proofs on $\TT/(\TT\cap\KK^{n-1-a})$ for each $f_a$. Then extend the value of all functions along the leaves of $\II$ by the formula of Theorem \ref{development}. So, you get the announced potentials $u$.\medskip

{\bf Part {\bf (b)} for cases {\bf (1)}, {\bf (2)} and {\bf (2')}} ---~thus also for their complexifications {\bf (1$^\C$)} and {\bf (2$^\C$)}. We adapt the standard arguments given in Proposition 2.4 of \cite{boubel2013a}. First, take nilomorphic coordinates being, at the origin in $\TT=\{(\nu y)=0\}$, tangent to normal coordinates {\em i.e.\@} such that the $X_i.g(X_j,X_k)$ vanish; in particular all $D_{X_i}X_j$ are null. As $DN=0$, so are the $D_{N^bX_i}N^aX_j$, and we still get, for any vectors $A,B,U,V$ among the $N^aX_i$ at the origin, that $g(R(A,B)U,V)$ is equal to:
$${\textstyle \frac12}\bigl(A.U.(g(B,V))-B.U.(g(A,V))-A.V.(g(B,U))+B.V.(g(A,U))\bigr).$$As $N$ is also self ajoint, $R(N^a A,B)=R(A,N^aB)$ for all $a$. The $R(X_i,N^aX_j)$ are determined by their restriction on $T\TT$ {\em i.e.\@} by the $R(X_i,N^aX_j)X_k$. We denote by $\widehat R$ this restriction. 

In case {\bf (1)}, for each $a$, $\widehat R(X_i,N^aX_j)$, which is defined at the origin on $T\TT/\ker N^a$, is the alternate part of the bilinear form:
$$\beta_{i,j,a}:(U,V)\mapsto X_i.U.g(N^aX_j,V)-X_j.U.g(N^aX_i,V),$$
also defined on $T\TT/\ker N^a$. For each $a$, the $\beta_{i,j,a}$ depend on the second derivatives at $0$ of the coefficients of $g(\,\cdot\,,N^a\cdot\,)$, defined on $\TT/(\TT\cap\KK^{a})$. Those derivatives are free in normal coordinates. Indeed $g(\,\cdot\,,N^a\cdot\,)$ is the coefficient of $\nu^{n-1-a}$ of the nilomorphic metric $h$ given in Theorem \ref{theoremenilpotentadj}, hence is chosen freely. So, on a dense open subset of the 2-jets of metrics, the alternate parts of the $(\beta_{i,j,a})_{i,j=1}^D$ are linearly independent and hence span a $\frac{\delta(\delta-1)}2$-dimensional space, with $\delta=\dim(\TT/(\TT\cap\KK^{a}))=\sharp\{i;N^aX_i\neq0\}$. So, the sum for all $a$ of those dimensions is the number $K$ of triples $(i,j,a)$ with $i<j$ and $N^aX_j\neq0$, and generically, the holonomy algebra is $K$-dimensional. Now an element $\gamma$ of $\goth{o}_d(\R)^N$ is precisely given by the $g(\gamma(N^aX_j),X_i)$ for those triples $(i,j,a)$ {\em i.e.\@} $\dim\goth{o}_d(\R)^N=K$. We are done.

The adaptation for case {\bf (2)} is similar. The forms $\beta_{i,j,a}$ are:
$$\beta_{i,j,a}:(Z_k,\Zbar_l)\mapsto {\textstyle \frac12}\bigl(-\Zbar_j.Z_k.(g(Z_i,N^a\Zbar_l))-Z_i.\Zbar_l.(g(\Zbar_j,N^aZ_k))\bigr),$$
defined for $i,j$ in $\llbracket1,\frac D2\rrbracket$ and $a\in\llbracket 0,n-1\rrbracket$, see the proof of Proposition 2.4 
of \cite{boubel2013a}. 
Each $\beta_{i,j,a}$ is given by the derivatives of $u_{n-1-a}$, the coeffcient of $\nu^{n-1-a}$ in the potential $u$ given by part {\bf (a)} of the theorem. This $u_{n-1-a}$ is defined on $\TT/(\TT\cap\KK^{a})$, so the $(\beta_{i,j,a})_{i,j}$ may be chosen freely and span a $\left(\frac\delta2\right)^2$-dimensional space, with $\delta=\dim(\TT/(\TT\cap\KK^{a}))=\sharp\{i;N^aX_i\neq0\}$. The sum for all $a$ of those dimensions is the number $K'$ of triples $(i,j,a)$ with $i<j\leqslant \frac D2$ and $N^aZ_j\neq0$. So generically, the holonomy algebra is $K'$-dimensional. Now an element $\gamma$ of $\goth{u}_{d/2}^N$ is precisely given by the $g(\gamma(N^aZ_j),\Zbar_j)$ for those triples $(i,j,a)$ {\em i.e.\@} $\dim\goth{u}_{d/2}^N=K'$. We are done. Case {\bf (2')} is entirely similar and left to the reader.\hfill{\rm q.e.d.}

\begin{importantrem}The note p.\@ \pageref{nilomorphisation} points out that the technique used in the proof of Theorem \ref{realisation_nilomorphe} may be used as a standard way to generalize reasonings on germs of real functions to germs of nilomorphic ones. In particular, by this means, we might show similar statements as Theorem \ref{realisation_nilomorphe} for $H$ any semi-simple classical pseudo-Riemannian holonomy group.
\end{importantrem}

\label{alg_lin}We finally state the sequence of results in linear algebra leading to Corollary \ref{cor_forme_e}. This provides the matrix of the elements of all algebras appearing here. To let clearly appear the reason the final Proposition \ref{forme_e} works, we give the statements as a sequence of steps; this makes the work to prove each step quite clear. So we let the proofs to the reader.

\begin{notation}\label{notation_forme_e1} The two cases are here labelled in reference to Theorem \ref{structure_s}. We recall standard facts in point denoted by {\bf (2)}, to introduce in {\bf (2')} a ``paracomplex'' counterpart of them.\medskip

{\bf (i)} If $J$ is a $g$-skew adjoint morphism of $\R^{2d}$ with $J^2=-\Id$, its commutant ${\rm U}(g)={\rm O}(g)^J$ in ${\rm O}(g)$ may be seen as a subgroup of GL$_d(\C)$. More precisely, in a basis such that $J=\diag(J_1,\ldots,J_1)$, ${\rm U}(g)$ is a group of real matrices consisting of square subblocks in:
$${\rm M}_2(\R)^{J_1}=\left\{\text{\small$\left(\begin{array}{cc}a&-b\\b&a\end{array}\right)$};a,b\in\R\right\}\simeq\{a+{\rm i}b;a,b\in\R\}=\C,$$
so in such a basis, these matrices may be considered as complex.\medskip

{\bf (ii)} Similarly, if $L$ is a $g$-skew adjoint morphism of $\R^{2d}$ with $L^2=\Id$, we denote here by ${\rm U}_L(g)={\rm O}(g)^L$ its commutant in ${\rm O}(g)$. It may be seen as a subgroup of GL$_d(\R\oplus\R)$, so as a group of matrices with entries in $\R\oplus\R$. More precisely, in a basis such that $L=\diag(I_{1,1},\ldots,I_{1,1})$, ${\rm U}_L(g)$ is a group of real matrices consisting of square subblocks in:
$${\rm M}_2(\R)^{I_{1,1}}=\left\{\text{\small$\left(\begin{array}{cc}a&0\\0&b\end{array}\right)$};a,b\in\R\right\}\simeq\{(a,b);a,b\in\R\}=\R\oplus\R,$$
with the natural product $(a,b).(a',b')=(aa',bb')$ on the ring $\R\oplus\R$. These subblocks may be considered as elements of M$_{d/2}(\R\oplus\R)$.\medskip

The inclusion $\R.\Id\subset{\rm M}_2(\R)^{J_1}\simeq\C$ sends $\R$ onto $\{a+{\rm i}b;b=0\}$, stabilized by the involution $a+{\rm i}b\mapsto \overline{a+{\rm i}b}=a-{\rm i}b$. Similarly $\R.\Id\subset{\rm M}_2(\R)^{I_{1,1}}\simeq\R\oplus\R$ sends $\R$ onto $\{(a,b);a=b\}$, stabilized by the natural involution $(a,b)\mapsto (b,a)$, called here ``$L$-'' or  ``paracomplex'' conjugation and denoted by $(a,b)\mapsto\overline{(a,b)}$. Notice this involution at line {\bf (2')} of Table  1 of \cite{boubel2013a}. To sum up, in terms of real 2-2 matrices or submatrices:
$$\text{in {\bf (i)}, }\overline{\text{\small$\left(\begin{array}{cc}a&-b\\b&a\end{array}\right)$}}=\text{\small$\left(\begin{array}{cc}a&b\\-b&a\end{array}\right)$}\ \text{and in {\bf (ii)}, }\overline{\text{\small$\left(\begin{array}{cc}a&0\\0&b\end{array}\right)$}}=\text{\small$\left(\begin{array}{cc}b&0\\0&a\end{array}\right)$}.$$
\end{notation}

\begin{notation}\label{notation_forme_e1bis} The two cases are here labelled in reference to Theorem \ref{structure_s}. We recall standard facts in point denoted by {\bf (3)} about the group Sp$(p,q)$, to introduce in {\bf (3')} a ``paraquaternionic'' counterpart of them.\medskip

{\bf (i)} If $\langle J,J',J''\rangle$ is a $g$-skew adjoint quaternionic structure on $\R^{4d}$, {\em i.e.\@} $(J^{(\prime\prime)})^2=-\Id$ and $JJ'=-J'J=J''$, its commutant ${\rm Sp}(g)={\rm O}(g)^{\{J,J',J''\}}$ in ${\rm O}(g)$ may be seen as a subgroup of GL$_d(\HH)$. More precisely, identifying $(\R^{4d},J)$ with $(\C^{2d},{\rm i}I_{1,1})$, in a basis such that $J'=-\diag(J_1,\ldots,J_1)$, ${\rm Sp}(g)$  is a group of real matrices consisting of 4-4 square subblocks in ${\rm M}_4(\R)^{\{J,J'\}}$ identified with:
$$\left\{\text{\small$\left(\begin{array}{cc}\alpha&-\overline\beta\\\beta&\overline\alpha\end{array}\right)$};\alpha,\beta\in\C\right\}\simeq\{a+{\rm i}a'+{\rm j}b+{\rm k}b';\alpha=a+{\rm i}a', \beta=b+{\rm i}b'\}=\HH,$$
so in such a basis, these matrices may be considered as quaternionic.\medskip

{\bf (ii)} Similarly, if $\langle L,L',J\rangle$ is a $g$-skew adjoint paraquaternionic structure on $\R^{4d}$, {\em i.e.\@} $(L^{(\prime)})^2=\Id$ and $LL'=-L'L=J$, we denote here by ${\rm Sp}_{L,L'}(g)$ its commutant ${\rm O}(g)^{\{L,L',J\}}$ in ${\rm O}(g)$. It may be seen as a subgroup of GL$_d({\rm M}_2(\R))$. More precisely,  in a basis such that $L=\diag(I_{2,2},\ldots,I_{2,2})$ and $L'=\diag(L'_2,\ldots,L'_2)$ where $L'_2=\text{\scriptsize$\left(\begin{array}{cc}0&J_1\\-J_1&0\end{array}\right)$}$, ${\rm Sp}_{L,L'}(g)$ is a group of real matrices consisting of square subblocks in:
$${\rm M}_4(\R)^{\{I_{2,2},L'_2\}}=\left\{h_P:=\left(\begin{array}{cc}P&0\\0&^t\!\widetilde P\end{array}\right);P\in{\rm M_2}(\R)\right\}\simeq{\rm M}_2(\R),$$
where $\widetilde P$ stands for the comatrix of $P$.\medskip

In {\bf (i)}, the inclusion $\R.\Id\subset{\rm M}_4(\R)^{\{J,J'\}}\simeq\HH$ sends $\R$ to $\{a+{\rm i}a'+{\rm j}b+{\rm k}b'; a'=b=b'=0\}$, stabilized by the quaternionic conjugation. The inclusion $\R.\Id\subset{\rm M}_4(\R)^{\{I_{2,2},L'_2\}}\simeq {\rm M}_2(\R)$, sends $\R$ to $\{h_{\lambda\Id}\}$, stabilized by $h_P\mapsto h_{^t\!\widetilde P}$ called here ``paraquaternionic conjugation'' and denoted by $h_P\mapsto \overline{h_P}$. This involution appears at line {\bf (3')} of Table 1 of \cite{boubel2013a}.
\end{notation}

\begin{notation}\label{notation_forme_e2} {\bf (a)} For each $p\in\N$ and for $\delta\in\{1,2,4\}$ we introduce:
$$%\begin{array}{l}
N_p^{(\delta)}:=\text{\small$\left(\begin{array}{cccc}0&I_\delta\\&\ddots&\ddots\\&&\ddots&I_\delta\\&&&0\end{array}\right)$}%\in {\rm M}_{\delta p}(\R)\text{,}\\
\ \text{and:}\ K_p^{(\delta)}:=\text{\small$\left(\begin{array}{ccc}&&I_\delta\\&\iddots\\I_\delta&&\end{array}\right)$}%\in {\rm M}_{\delta p}(\R).\end{array}
,$$
both in ${\rm M}_{\delta p}(\R)$. We will also denote $N_p^{(1)}$ by $N_p$ and $K_p^{(1)}$ by $K_p$.

{\bf (b)} If $M=(m_{i,j})_{i=1,}^p{}_{j=1}^q$ is a matrix with $p$ lines and $q$ columns, we denote here by $^{\slash\hspace*{-.8ex} t}M$ its ``transpose with respect to the anti diagonal'' $^{\slash\hspace*{-.8ex} t}M=(m_{q-j,p-i})_{j=1,}^q{}_{i=1}^p=K_q.^t\!M.K_p$.
\end{notation}

\begin{lem-notation}\label{bases_privilegiees} Let $g$ be a pseudo-Riemannian metric on $\R^d$, and $N$ a $g$-self adjoint nilpotent endomorphism of $\R^d$, of nilpotence index $n$. In five different cases, labelled as in Theorem \ref{structure_s}, let us introduce a Lie subgroup $Q$ of GL$_d(\R)$, its Lie algebra $\goth q$, and an induced set of bases of $\R^d$, called here ``priviledged'', on which $Q^N$ acts simply transitively. The $N_i$, $N_i^{(\delta)}$, $K_i$, $K_i^{(\delta)}$ are as in Notation \ref{notation_forme_e2} {\bf (a)} and $\varepsilon_i\in\{-1,1\}$.\medskip

{\bf (1)} Here $Q:={\rm O}(g)$. There are bases of $\R^d$ in which:
$$\left\{\begin{array}{l}\Mat(N)=\diag(N_{n_1},\ldots,N_{n_D})\ \text{with }n=n_1\geqslant\ldots\geqslant n_D\\\Mat(g)=\diag(\varepsilon_{n_1} K_{n_1},\ldots,\varepsilon_{n_D} K_{n_D}).\end{array}\right.$$

{\bf (2)} Here $J\in\End(\R^d)$ with $J^\ast=-J$ and $J^2=-\Id$; $Q:=U(g)={\rm O}(g)^J$. There are bases of $\R^d$ in which $\Mat(J)=\diag(J_{1},\ldots,J_{1})$ and: $$\left\{\begin{array}{l}\Mat(N)=\diag(N_{n_1}^{(2)},\ldots,N_{n_D}^{(2)})\ \text{with }n=n_1\geqslant\ldots\geqslant n_D\\\Mat(g)=\diag(\varepsilon_{n_1} K_{n_1}^{(2)},\ldots,\varepsilon_{n_D} K_{n_D}^{(2)}).\end{array}\right.$$

{\bf (2')} Here $L\in\End(\R^d)$ with $L^\ast=-J$ and $L^2=-\Id$; $Q:=U_L(g)={\rm O}(g)^L$, see Notation \ref{notation_forme_e1} {\bf (ii)}. There are bases of $\R^d$ in which $\Mat(L)=\diag(I_{1,1},\ldots,I_{1,1})$ and:
$$\left\{\begin{array}{l}\Mat(N)=\diag(N_{n_1}^{(2)},\ldots,N_{n_D}^{(2)})\ \text{with }n=n_1\geqslant\ldots\geqslant n_D\\\Mat(g)=\diag(K_{2n_1},\ldots,K_{2n_D}).\end{array}\right.$$

{\bf (3)} Here $(J,J',J'')$ is a $g$-skew adjoint quaternionic structure on $\R^{d}$, see \ref{notation_forme_e1bis} {\bf (i)}; $Q:={\rm Sp}(g)={\rm O}(g)^{\{J,J'\}}$. There are bases of $\R^d\simeq\C^{d/2}$ in which $\Mat(J)={\rm i}\diag(I_{1,1},\ldots,I_{1,1})$, $\Mat(J)=-\diag(J_{1},\ldots,J_{1})$ and:
$$\left\{\begin{array}{l}\Mat(N)=\diag(N_{n_1}^{(4)},\ldots,N_{n_D}^{(4)})\ \text{with }n=n_1\geqslant\ldots\geqslant n_D\\\Mat(g)=\diag(\varepsilon_{n_1} K_{n_1}^{(4)},\ldots,\varepsilon_{n_D} K_{n_D}^{(4)}).\end{array}\right.$$

{\bf (3')} Here $(L,L',J)$ is a $g$-skew adjoint paraquaternionic structure on $\R^{d}$, see \ref{notation_forme_e1bis} {\bf (ii)}; $Q={\rm Sp}_{L,L'}(g)={\rm O}(g)^{\{L,L'\}}$. There are bases of $\R^d$ in which $\Mat(L)=\diag(I_{2,2},\ldots,I_{2,2})$, $\Mat(L')=\diag(L'_{2},\ldots,L'_{2})$ with $L'_2$ as in \ref{notation_forme_e1bis} {\bf (ii)} and:
$$\left\{\begin{array}{l}\Mat(N)=\diag(N_{n_1}^{(4)},\ldots,N_{n_D}^{(4)})\ \text{with }n=n_1\geqslant\ldots\geqslant n_D\\\Mat(g)=\diag(K_{2n_1}^{(2)},\ldots,K_{2n_D}^{(2)}).\end{array}\right.$$
\end{lem-notation}

\noindent{\em Reference for the proof.} Priviledged bases are provided by \cite{leep-schueller}.

\begin{lem}\label{lemme_commutants} Take $\gamma\in\End(\R^d)$. Using a priviledged basis given by Lemma \ref{bases_privilegiees}, we consider its matrix $M$ as a matrix with coefficients in $\Abb:=\R$, $\Abb:=\C$, $\Abb:=\R\oplus\R$, $\Abb:=\HH$ or $\Abb:={\rm M}_2(\R)$ in cases {\bf (1)}, {\bf (2)}, {\bf (2')}, {\bf (3)} or {\bf (3')} respectively, see Notation \ref{notation_forme_e1} and \ref{notation_forme_e1bis}. Consider $(M_{i,j})_{i,j=1}^D$ the block-decomposition of $M$ corresponding to the Jordan blocks of $N$, and the Lie algebra $\goth q$ introduced in Lemma \ref{bases_privilegiees}.\medskip

{\bf (a)} $\gamma\in\goth q$ if and only if, with the notation $^{\slash\hspace*{-.8ex} t}$ set in \ref{notation_forme_e2} {\bf (b)}:\medskip

-- in case {\bf (1)}, $M_{j,i}=-\varepsilon_i\varepsilon_j{}^{\slash\hspace*{-.8ex} t}\!M_{i,j}$ for all $i,j\in \llbracket 1,D\rrbracket$,\medskip

-- in cases {\bf (2)} and {\bf (3)}, $M_{j,i}=-\varepsilon_i\varepsilon_j{}^{\slash\hspace*{-.8ex} t}\overline{M_{i,j}}$ for all $i,j\in \llbracket 1,D\rrbracket$, with the complex or quaternionic conjugation in cases {\bf (2)} and {\bf (3)} respectively,\medskip

-- in case {\bf (2')} and {\bf (3')}, $M_{j,i}=-{}^{\slash\hspace*{-.8ex} t}\overline{M_{i,j}}$ for all $i,j\in \llbracket 1,D\rrbracket$, with the paracomplex or paraquaternionic conjugation, in cases {\bf (2')} and {\bf (3')} respectively, introduced in Notation \ref{notation_forme_e1} {\bf (2')} and \ref{notation_forme_e1bis} {\bf (3')}.\medskip

{\bf (b)} $\gamma$ commutes with $N$, and additionally with $J$  in case {\bf (2)}, $L$ in case {\bf (2')}, $\{J,J',J''\}$ in case {\bf (3)} or $\{L,L',J\}$ in case {\bf (3')} if and only if, for all $i\leqslant j$:
$$(\pentagram)\left\{\begin{array}{l}M_{i,j}=\text{\small$\left(\!\!\begin{array}{c}M'_{i,j}\\0_{n_i-n_j,n_j}\end{array}\!\!\right)$},\ M_{j,i}=\left(\!\begin{array}{cc}0_{n_j,n_i-n_j}&\!M'_{j,i}\end{array}\!\right)\!,\ \text{\em i.e.\@ }M_{i,i}=M'_{i,i}\\[.3cm]
\text{where: }M'_{i,j}=\text{\small$\displaystyle\left(\begin{array}{cccc}m_{i,j}^1&m_{i,j}^2&\ldots &m_{i,j}^{n_j}\\
                             0 &m_{i,j}^1&\ddots&\vdots\\
\vdots&\ddots&\ddots&m_{i,j}^2\\
0&\cdots&0&m_{i,j}^1\\
                             \end{array}\right)$}\in{\rm M}_{n_k}(\Abb).\end{array}\right.$$
\smallskip

{\bf (c)} Hence, $\gamma\in\goth q^N$ if and only if the $M_{i,j}$ are as above and:\medskip

-- in case {\bf (1)}, $M'_{j,i}=-\varepsilon_i\varepsilon_jM'_{i,j}$ for all $i,j\in \llbracket 1,D\rrbracket$,\medskip

-- in cases {\bf (2)} and {\bf (3)}, $M'_{j,i}=-\varepsilon_i\varepsilon_j\overline{M'_{i,j}}$ for all $i,j\in \llbracket 1,D\rrbracket$,\medskip

-- in cases {\bf (2')} and {\bf (3')}, $M'_{j,i}=-\overline{M'_{i,j}}$ for all $i,j\in \llbracket 1,D\rrbracket$, with the conjugation introduced in Notation \ref{notation_forme_e1} {\bf (2')} and \ref{notation_forme_e1bis} {\bf (3')}.\medskip

Indeed, the $M'_{i,j}$ being as in $(\pentagram)$, ${}^{\slash\hspace*{-.8ex} t}\!M'_{i,j}=M'_{i,j}$. So finally, for $i<j$, $M_{j,i}$ is given by $M_{i,j}$, which is freely chosen as in $(\pentagram)$, and:\medskip

-- in case {\bf (1)}, the $M'_{i,i}$ are null,\medskip

-- in cases {\bf (2)}--{\bf (3)}, the $M'_{i,i}$ are as in $(\pentagram)$, with purely imaginary coefficients {\em i.e.\@} in ${\rm i}\R\subset\C$ for case {\bf (2)} and in $\Span\{{\rm i},{\rm j},{\rm k}\}\subset\HH$ for {\bf (3)}.\medskip

-- in cases {\bf (2')}--{\bf (3')}, the $M'_{i,i}$ are as in $(\pentagram)$, with coefficients of the type $(a,-a)\in\R\oplus\R$ for {\bf (2')}, and $\diag(P,{}^t\!\widetilde P)$ with $\tr P=0$ for {\bf (3')}.\medskip

{\bf (d)} Therefore, $\gamma\in\End(\R^d)^{(\goth q^N)}$ if and only if:\medskip

-- in the other  cases than {\bf (1)}, all the $M_{i,j}$ are null except when $i=j$, where $M_{i,i}$ is as in $(\pentagram)$.\medskip

-- in case {\bf (1)}, the situation is the same except for some similarity types of $N$. All the $M_{i,j}$ are null except when $i=j$, and except possibly $M_{1,2}$ and $M_{2,1}$. Remind that $n=n_1\geqslant n_2\geqslant\ldots\geqslant n_D$ are the sizes of the Jordan blocks of $N$; conventionally, we set $n_i=0$ for $i>D$, for example $n_3=0$ if $D=2$ {\em i.e.\@} $N$ has two Jordan blocks. Then:
\begin{itemize}
\item[$\bullet$] All $M_{i,i}$ for $i\geqslant2$ are as in $(\pentagram)$, and $M_{1,1}=M'_{1,1}+M''_{1,1}$ where $M'_{1,1}$ is as in $(\pentagram)$ and:
$$M''_{1,1}=\text{\small$\left(\begin{array}{cc}0_{n_1-n_2,n_2}&\ast\\0_{n_2,n_2}&0_{n_2,n_1-n_2}\end{array}\right)$},\ \text{with }\ast\text{ arbitrary in }{\rm M}_{n_1-n_2}(\R).$$
\item[$\bullet$] $\displaystyle M_{1,2}=\text{\small$\left(\begin{array}{cc}0_{n_2-n_3,n_3}&M''_{1,2}\\0_{n_3,n_3}&0_{n_3,n_2-n_3}\end{array}\right)$}$ and $M_{2,1}=-^{\slash\hspace*{-.8ex} t}M_{1,2}$, with:
$$M''_{1,2}=\text{\small$\displaystyle\left(\begin{array}{cccc}m^1&m^2&\ldots &m^{n_2-n_3}\\
                             0 &m^1&\ddots&\vdots\\
\vdots&\ddots&\ddots&m^2\\
0&\cdots&0&m^1\\
                             \end{array}\right)$}\in{\rm M}_{n_2-n_3}(\R).$$
\end{itemize}
Thus, $M''_{1,1}=0$ {\em i.e.\@} $M_{1,1}$ is as in $(\pentagram)$ if and only if $n_1=n_2$ {\em i.e.\@} $d_n>1$ {\em i.e.\@} $N$ has several Jordan blocks of maximal size, and $M_{1,2}=-^{\slash\hspace*{-.8ex} t}M_{1,2}=0$ if and only if $n_2=n_3$. So  {\bf (1)} is like the other cases if and only if $n_1=n_2=n_3$.
\end{lem}

\noindent{\em Hint for the proof.} In {\bf (d)}, for cases {\bf\mathversion{bold}($2$)\mathversion{normal}}--{\bf\mathversion{bold}($2^{\prime}$)\mathversion{normal}} and {\bf\mathversion{bold}($3$)\mathversion{normal}}--{\bf\mathversion{bold}($3^{\prime}$)\mathversion{normal}}, to get $M_{i,j}=0$ if $i\neq j$, it is sufficient to involve two Jordan blocks. For case {\bf (1)}, involving three blocks is necessary to get the result, and four or more give no additional constraint.

\begin{rem}
The reason case {\bf (1)} behaves differently is that $M_{i,i}=-{}^{\slash\hspace*{-.8ex} t}M_{i,i}\Rightarrow M_{i,i}=0$, whereas $M_{i,i}=-{}^{\slash\hspace*{-.8ex} t}\overline{M_{i,i}}\not\Rightarrow M_{i,i}=0$.
\end{rem}

\begin{prop}\label{forme_e}Take $\goth s$ any of the subalgebras of $\End(\R^d)$ appearing in Theorem \ref{structure_s}. With the assumptions of Lemma \ref{lemme_commutants}:\medskip

-- In cases {\bf (2)}, {\bf (2')}, {\bf (3)} and {\bf (3')}, the action of $\goth{o}(g)^{\goth{s} \cup\{N\}}$ on $(\R^d,g)$ is indecomposable and: $\End(T\MM)^{\left(\goth{o}(g)^{\goth{s} \cup\{N\}}\right)}=\goth{s}.\R[N].$\medskip

-- In case {\bf (1)}, the action of $\goth{o}(g)^{\goth{s} \cup\{N\}}=\goth{o}(g)^{N}$ on $(\R^d,g)$ is indecomposable if and only if $2n_2<n_1$ {\em i.e.\@} $N$ has only one Jordan block of maximal size, the second largest size being more than twice smaller.\medskip

$\bullet$ If $2n_2<n_1$, take $F$ any supplement of $\ker N^{n_2}$ in $\im N^{n_2}$. By Lemma \ref{lemme_commutants} {\bf (c)}, $\goth{o}(g)^{N}$ acts trivially on $\im N^{n_2}$ thus on $F$; besides $\ker N^{n_2}=\ker(g_{|\im N^{n_2}})$, so $g'':=g_{|F}$ is non degenerate. So $\goth{o}(g)^{N}$ acts trivially on $F$ and stabilizes $F^\perp$. In particular, $\End(F)^{\left(\goth{o}(g)^{N}\right)}=\End(F)$. Moreover, $d'':=\dim F=n_1-2n_2$ and the signature of $g''$ is:\medskip

$\ast$ $\textstyle(\frac{d''}2,\frac{d''}2)$ if $d''$ is even {\em i.e.\@} $n$ is,\medskip

$\ast$ $(\frac{d''+1}2,\frac{d''-1}2)$, respectively $(\frac{d''-1}2,\frac{d''+1}2)$, if $d''$ is odd {\em i.e.\@} $n$ is, and $g(\,\cdot\,,N^{n-1}\,\cdot\,)$ is positive, respectively negative on the line $\R^d/\ker N^{n-1}$.\medskip

On $F^\perp$, the operator $N'=N_{|F^{\perp}}$ is of nilpotence index $n'=2n_2$ and has Jordan blocks of sizes $(2n_2,n_2,n_3,\ldots,n_D)$. So the action of $\goth{o}(g)^{N}$ on $F^\perp$ is described by the second item below.\medskip

$\bullet$ If $2n_2\geqslant n_1$, set $\pi':\R^d\rightarrow\R^d/\ker N^{n_2}$ and $\pi'':\R^d\rightarrow\R^d/\ker N^{n_3}$. The (pseudo\nobreakdash-)euclidian product $g_{n_3}:=g(\,\cdot\,,N^{n_3}\,\cdot\,)$ is well-defined and non degenerate on $\R^d/\ker N^{n_3}$, on which $N$ acts. This defines the group ${\rm O}(g_{n_3})$ and the commutant ${\rm O}(g_{n_3})^N$ of $N$ in it. Besides, we may see the $N^{n_i}$ as morphisms $\R^d/\ker N^{n_i}\rightarrow\im N^{n_i}$. Then:
\begin{gather*}
\End(\R^d)^{\left(\goth{o}(g)^{N}\right)}\!=\!
\left(\R[N]+N^{n_2}_\ast\pi^{\prime\ast}\End\bigl(\R^d/\ker N^{n_2}\bigr)\right)\!\oplus\! N^{n_3}_\ast\pi^{\prime\prime\ast}\goth{o}(g_{n_3})^{N}.
\end{gather*}
\end{prop}

\noindent{\em Hint for the proof.} For the (in)decomposability, you see in Lemma \ref{lemme_commutants} {\bf (c)} that, in cases {\bf\mathversion{bold}($2$)\mathversion{normal}}--{\bf\mathversion{bold}($2'$)\mathversion{normal}} and {\bf\mathversion{bold}($3$)\mathversion{normal}}--{\bf\mathversion{bold}($3^{\prime}$)\mathversion{normal}}, no proper subspace of $\R^d$ is stable by $\goth{o}(g)^{\goth{s}\cup\{N\}}$, except possibly isotropic lines in some very exceptional subcases of {\bf\mathversion{bold}($2^{\prime}$)\mathversion{normal}}. In case {\bf (1)}, the subspace on which $\goth{o}(g)^N$ acts trivially is exactly $\{0\}$ if $n_1=n_2$ and $\im N^{n_2}$ if $n_1>n_2$.

\section{\mathversion{bold}A glimpse on metrics whose holonomy group is the commutant of {\em several} algebraically independent nilpotent endomorphisms\mathversion{normal}}\label{nonprincipal}

We investigate here the simplest example where $\goth{h}=\goth{o}(g)^{\{N,N'\}}$ with $(N,N')$ algebraically independent. This will show a phenomenon appearing when $\goth n$ is abelian, non principal, see Comment \ref{commentaire}. We will see that the results of both Parts 1 and 2 are needed to describe it.
Here is it. We parametrize the set ${\mathcal G}$ of germs of metrics such that $\goth n=(N,N')$ with $N$ and $N'$ self adjoint and $N^2=N'^2=NN'=N'N=0$. To simplify a bit more, we assume $\im N=\ker N$ ---~our goal here is no kind of general theory. Besides, taking $\dim\MM\geqslant 6$ ensures $n_1=n_2=n_3$ so we are not in the exceptional cases of Corollary \ref{cor_forme_e}.

\begin{rem}\label{u}{\bf (i)} Then, there is a $U\in\End(T\MM/\im N)$ such that $N'=NU$, which makes sense as the argument of $N$ may be defined only modulo $\ker N$. Indeed, $N$ gives an isomorphism $\theta_N:T\MM/\ker N\rightarrow\ker N$. As $NN'=N'N$, $N'$ gives also a morphism $\theta_{N'}:T\MM/\ker N\rightarrow\ker N$, and is determined by it. Set $U:=\theta_N^{-1}\circ\theta_{N'}\in\End(T\MM/\ker N)$.

{\bf (ii)} Identify $U$ and its matrix. In coordinates $(y_i,x_i)_{i=1}^{d/2}$ adapted to $N$:
$$\Mat(N)=\text{\small$\left(\begin{array}{cc}0&I\\0&0\end{array}\right)$}\ \text{ and: }\ \Mat(N')=\text{\small$\left(\begin{array}{cc}0&U\\0&0\end{array}\right)$}.$$
Then, using Lemma \ref{integrabilite_amelioree} complemented by Remark \ref{integrabilite_amelioree_bis}, we take coordinates such that Mat$(U)$ is also constant ---~and for instance in Jordan form. Explicitly, $U$ is well defined on $\MM/\II$ (in other words it is $\II$-basic) and, if $\prod_\alpha P_\alpha^{n_\alpha}$ is the decomposition of its minimal polynomial in powers of irreducible polynomials, we may take coordinates $(x_i)_i$ which are:\medskip

-- product coordinates for $\MM/\II\simeq\prod_k\MMbar_k$, the integration of the decomposition $\overset{\perp}\oplus_\alpha\ker P_\alpha^{n_\alpha}$,\medskip

-- on each factor, adapted to the nilpotent part of $U$ on it and, if $\deg P_\alpha=2$ {\em i.e.\@} if the semi-simple part of $U$ induces a complex structure $\underline J_\alpha$ on it, also complex coordinates for it.
\end{rem}

\begin{prop}\label{nn'} With the $U\in\End(T\MM/\im N)$ introduced in Remark \ref{u} {\bf (i)}, a metric $g$ makes $N$ and $N'$ parallel if and only if:\medskip

{\bf (a)} it is the real metric $h_1$ associated with a $(\nu,N)$-nilomorphic metric $h=h_0+\nu h_1$ with value in $\R[\nu]=\R[X]/(X^2)$,\medskip

{\bf (b)} the bilinear form $h_0$, defined on $\MM/\II$, makes $U$ parallel (recall that $h_0=g(\,\cdot\,,N\,\cdot\,)$ so is nondegenerate, hence is a metric, on $\MM/\II$).\medskip

Then there exist coordinates $(x_i,y_i)_i$ that are simultaneously $N$- and $N'$-adapted, as given in Remark \ref{u} {\bf (ii)}. In such coordinates, on each factor $\MMbar_k$, $h_0$ is itself of the form given by Theorem \ref{theoremenilpotentadj}, for the nilpotent part $N_k$ of $U$ on $\MMbar_k$ (\/{\em i.e.\@} $h_0$ is the real metric associated with some $(\mu,N_k)$-nilomorphic metric, with $\R[\mu]=\R[X]/(X^{n_k})$\/), and also complex Riemannian for $\underline J_k$ if $\deg P_k=2$ (then $\R[\mu]$ is replaced by $\C[\mu]$).
\end{prop}

\noindent{\bf Proof.} After Th.\@ \ref{theoremenilpotentadj}, $g$ makes $N$ parallel if and only if it satisfies {\bf (a)}. It makes also $N'$ parallel if and only if for any $(i,j,k)$, $g(D_{X_i}N'X_j,X_k)=g(D_{X_i}X_j,N'X_k)$, that is to say $g(D_{X_i}NUX_j,X_k)=g(D_{X_i}X_j,NUX_k)$, or $g(D_{X_i}UX_j,NX_k)=g(D_{X_i}X_j,NUX_k)$, as $DN=0$. This means $h_0(D_{X_i}UX_j,X_k)=h_0(D_{X_i}X_j,UX_k)$ {\em i.e.\@} $U$ is parallel for $h_0$.\hfill{\rm q.e.d.}\medskip

In real terms, following Example \ref{exempleA}, and in the basis $((Y_i)_i,(X_i)_i)$, $\Mat(g)=\text{\small$\left(\begin{array}{cc}0&G^0\\G^0&G^1\end{array}\right)$}$ where% $G^a=\Mat(h_{1-a})$ and
:\medskip

-- $G^0$, depending only on the $(x_i)_i$, is the matrix of a metric on $\MM/\II$ making $U$ parallel,\medskip

-- $G^1=B^1+\sum_{i}\bigl(\frac{\partial G^{1}}{\partial x_i}\bigr)y_{i,1}$, with $B^1$ the matrix of any bilinear symmetric form, depending only on the $(x_i)_i$.\medskip

\noindent So building a metric making $N$ and $N'$ parallel means taking a metric such that $N$ is, and adding constraints on it by repeating the general story of this article, on the quotient $\MM/\II$, endowed with the metric $g(\,\cdot\,,N\,\cdot\,)$ and of a parallel endomorphism $U$ induced by $N'$.\medskip

\begin{comment} \label{commentaire}If $(N,N')$ is a general pair of commuting nilpotent endomorphisms, their characteristic flags $F^{a,b}$ introduced in Table \ref{table_F} p.\@ \pageref{table_F} induce, each alone  and with each other, many quotient spaces $E_\alpha$ on each of which the metric $g$ and the pair $(N,N')$ induce an endomorphism $U_\alpha$ and a metric $g_\alpha$. Take a metric $g$ making $N$ parallel, and $\mathcal T$ a transversal to $\II$. If $g$ makes also $N'$ parallel, it makes each $U_\alpha$ parallel for $g_\alpha$. This means that some coefficients of $g$ on $\mathcal T$ must be themselves real parts of $N_\alpha$-nilomorphic functions as in Theorem \ref{realisation_nilomorphe}, with $N_\alpha$ the nilpotent part of $U_\alpha$. In fact, this is not sufficient, as a pair of commuting nilpotent endomorphisms is not something simple ---~{\em e.g.\@} consider simply the case of Proposition \ref{nn'}, without the assumption $\im N=\ker N$. If more than two endomorphisms $N$ and $N'$ are involved, it may appear quotients of the type of the $E_\alpha$ on which any of the situations of Theorem \ref{realisation_nilomorphe} appears, even if $\goth s=\R\Id$ on $T\MM$ itself.

So the general situation, though not entirely new with respect to the case where $\goth n$ is principal, seems to be complicated. What the good next questions are it is still unclear. The case where $\goth n$ is not abelian may induce new phenomena, but is strongly constrained, when $\goth n$ consists of self adjoint elements, by Proposition 1.8 of \cite{boubel2013a}.\end{comment}

\noindent Charles Boubel\\Institut de Recherche Math\'ematique Avanc\'ee, UMR 7501 -- 
Universit\'e de Strasbourg et CNRS,\\
7 rue Ren\'e Descartes\\
67084 STRASBOURG CEDEX, FRANCE

\end{document}